\documentclass{amsart}
\usepackage[utf8]{inputenc}

\usepackage{amsthm, amsmath, amsfonts, amssymb, mathrsfs}
\usepackage{mathtools,thmtools}

\usepackage{enumerate}
\usepackage{hyperref}
\usepackage[capitalise]{cleveref}
\usepackage[margin=1.25in]{geometry}

\usepackage{float}

\numberwithin{equation}{section}

\newtheorem{theorem}{Theorem}[section]
\newtheorem{proposition}[theorem]{Proposition}
\newtheorem{lemma}[theorem]{Lemma}

\newtheorem{fact}[theorem]{Fact}
\newtheorem{corollary}[theorem]{Corollary}
\newtheorem{conjecture}[theorem]{Conjecture}

\theoremstyle{definition}
\newtheorem{definition}[theorem]{Definition}

\newtheorem{question}[theorem]{Question}
\newtheorem{remark}[theorem]{Remark}
\newtheorem{construction}[theorem]{Construction}

\newcommand{\CC}{\mathbb C}

\newcommand{\FF}{\mathbb F}
\newcommand{\RR}{\mathbb R}
\newcommand{\ZZ}{\mathbb Z}
\newcommand{\cB}{\mathcal B}
\newcommand{\cF}{\mathcal F}
\newcommand{\cG}{\mathcal G}

\newcommand{\bP}{\mathbf P}
\newcommand{\bU}{\mathbf U}
\newcommand{\bV}{\mathbf V}
\newcommand{\bW}{\mathbf W}
\newcommand{\bX}{\mathbf X}
\newcommand{\bZ}{\mathbf Z}

\newcommand{\eigm}{\mathsf m}
\newcommand{\eps}{\varepsilon}

\newcommand{\dspan}[1]{\langle #1 \rangle_{\pm}}

\DeclareMathOperator{\im}{im}
\DeclareMathOperator{\codim}{codim}
\DeclareMathOperator{\pdim}{pdim}
\DeclareMathOperator{\ord}{ord}
\DeclareMathOperator{\spn}{span}
\DeclareMathOperator{\A}{A}
\DeclareMathOperator{\Syl}{Syl}
\DeclareMathOperator{\SL}{SL}
\DeclareMathOperator{\PGL}{PGL}
\DeclareMathOperator{\GL}{GL}

\DeclareMathOperator{\Cay}{Cay}
\DeclareMathOperator{\Lmul}{L*}
\DeclareMathOperator{\Rmul}{R*}
\DeclareMathOperator{\Mat}{Mat}
\DeclareMathOperator{\End}{End}
\DeclareMathOperator{\Conj}{Conj}
\DeclareMathOperator{\lcm}{lcm}
\DeclareMathOperator{\rk}{rk}

\DeclareMathOperator{\tr}{tr}
\DeclareMathOperator{\vtr}{vtr}

\newcommand{\nimplies}{\;\not\nobreak\!\!\!\!\implies} 

\DeclarePairedDelimiter\abs{\lvert}{\rvert}
\DeclarePairedDelimiter\floor{\lfloor}{\rfloor}
\DeclarePairedDelimiter\ceil{\lceil}{\rceil}

\DeclarePairedDelimiter{\pr}{\lparen}{\rparen}
\DeclarePairedDelimiter{\br}{\lbrack}{\rbrack}
\DeclarePairedDelimiter{\set}{\lbrace}{\rbrace}

\usepackage{etoolbox}
\AtBeginEnvironment{quote}{\par\small}

\title[Abelian structure and Alon's conjecture]{Abelian structure in approximate groups and \\ Alon's conjecture on Ramsey Cayley graphs}
\author{Carl Schildkraut}
\address{Department of Mathematics, Stanford University}
\email{\texttt {carlsch@stanford.edu}}
\date\today

\begin{document}

\begin{abstract} 
A result of Pyber states that every finite group $G$ contains an abelian subgroup whose order is quasi-polynomially large in $\lvert G\rvert$. We prove a similar result for $K$-approximate subgroups of solvable groups under only modest restrictions on $K$. We show that, if $A$ is a finite $K$-approximate group contained in some solvable group, then some abelian group intersects $A^4$ in at least $\exp(\Omega(\log^{1/6}\lvert A\rvert/\log 2K))$ elements. We also prove a similar result for approximate subgroups of finite groups with no large alternating subquotients. Along the way, we obtain polynomial (instead of quasi-polynomial) bounds for the same statement of approximate subgroups of linear groups.

We give two applications. Firstly, we consider the conjecture of Alon that every finite group $G$ admits a Cayley graph with clique number and independence number $O(\log\lvert G\rvert)$. Conlon, Fox, Pham, and Yepremyan have recently proven that, for almost all positive integers $N$, every abelian group of order $N$ satisfies Alon's conjecture. Extending their result, we verify Alon's conjecture for all (not necessarily abelian) groups of almost all orders. Secondly, we prove a ``local'' version of Roth's theorem in (many) non-abelian settings with quasi-polynomial bounds, using the recent breakthroughs of Kelley and Meka on Roth's theorem and of Jaber, Liu, Lovett, Ostuni, and Sawhney on the corners problem.
\end{abstract}

\maketitle

\section{Introduction}\label{sec:intro}

Many questions in additive combinatorics can be asked in non-abelian settings, wherein several tools from the abelian world are unavailable. Progress can often be made in such settings by looking within smaller abelian structures, when they exist. The simplest type of such arguments involves partitioning a group into cosets of an abelian subgroup and working within each coset. This technique, introduced by Bergelson, McCutcheon, and Zhang \cite[Corollary~6.4]{BMZ}, has been applied to many problems in arithmetic combinatorics, including non-abelian versions of Roth's theorem \cite{BMZ,Solymosi,SandersRoth} and the corners problem \cite{Austin,JLLOS}, as well as the Brown--Erd\H{o}s--S\'os conjecture \cite{SolymosiBES, SWBES, LongBES, NSTBES, WongBES}.

To apply such an argument, one needs to find large abelian subgroups of groups. In the setting of general (finite) groups, this can be done with the following result of Pyber.

\begin{theorem}[{\cite[Theorem~1.1]{Pyber}}]\label{thm:Pyber} There exists an absolute constant $c>0$ for which every finite group $G$ has an abelian subgroup of order at least $e^{c\sqrt{\log\abs{G}}}$.
\end{theorem}

This strategy of working within cosets of abelian groups is most useful when the sets one is concerned with are dense in the ambient group. It is often necessary, however, to understand the multiplicative structure of sparse subsets of finite groups, or of subsets of infinite groups. From the perspective of non-abelian arithmetic combinatorics, it is particularly important to consider sets $A$ with \emph{small tripling}, i.e., those sets for which the product set $A^3:=\{a_1a_2a_3:a_1,a_2,a_3\in A\}$ is not too much larger than $A$ itself. (These sets are closely related to the well-studied \emph{approximate groups}.) We refer the reader to \cite{Tointon,ShkredovSurvey} for the history behind this definition, and for the myriad applications of understanding such sets in non-abelian settings.

To this end, our main results are the following two extensions of \cref{thm:Pyber} to sets with small tripling (for example, to approximate groups). A set $A$ is \emph{symmetric} if $x^{-1}\in A$ whenever $x\in A$.

\begin{theorem}\label{thm:solv-abelian-struct} There exists an absolute constant $c>0$ for which the following holds. Let $G$ be a (not necessarily finite) solvable group, and let $A\subset G$ be a finite symmetric subset satisfying $\abs{A^3}\leq K\abs{A}$. Then there exists some commuting set $T\subset A^4$ satisfying
\[\abs{T}\geq\exp\pr*{c\frac{\log^{1/6}\abs{A}}{\log 2K}}.\]
\end{theorem}

\begin{theorem}\label{thm:gen-abelian-struct} There exist absolute constants $C,c>0$ for which the following holds. Let $G$ be a finite group and let $d$ be a positive integer. Suppose that $G$ has no subquotient\footnote{A \emph{subquotient} of $G$ is a quotient of a subgroup of $G$.} isomorphic to the alternating group $\A_d$. Let $A\subset G$ be a finite symmetric subset satisfying $\abs{A^3}\leq K\abs{A}$. Then there exists some commuting set $T\subset A^4$ satisfying
\[\abs{T}\geq\exp\pr*{\frac{c\log^{1/8}\abs{A}}{\exp\pr[\big]{C\log^2d\log(1+\log d)}\cdot\log 2K}}.\]
\end{theorem}

Like \cref{thm:Pyber}, the proof of \cref{thm:gen-abelian-struct} relies on the classification of finite simple groups. We state the particular corollary of the classification that we use as \cref{lem:comp-rank-alternating}.

\cref{thm:solv-abelian-struct,thm:gen-abelian-struct} can be viewed as generalizations of \cref{thm:Pyber} from the setting of groups to that ``near-groups.''\footnote{The appearance of the exponent $4$ in the conclusions of \cref{thm:solv-abelian-struct,thm:gen-abelian-struct} should not be cause for concern. If $A$ has small tripling, then $A^4$ is not much larger than $A$ itself, and we shall see later (in \cref{lem:abelian-struct-near}) that the commuting sets $T$ we obtain from these two theorems are related to $A$ closely enough for our purposes.} In the most structured setting where $K=O(1)$ (and, in \cref{thm:gen-abelian-struct}, $d=O(1)$ as well), both results also give quasi-polynomial dependence on the size $\abs{T}$ of the commuting subset in $\abs{A}$, albeit with different exponents. In this sense, \cref{thm:solv-abelian-struct,thm:gen-abelian-struct} are qualitatively optimal: since \cref{thm:Pyber} is tight even in the solvable setting \cite{Olshanskii,BGH}, one cannot necessary find a commuting set of size larger than quasi-polynomial in $\abs{A}$ even when $A$ is a genuine solvable group.

To demonstrate the utility of \cref{thm:solv-abelian-struct,thm:gen-abelian-struct}, we give two applications to non-abelian arithmetic combinatorics. Our first application, which is our primary motivation for this work, is to settle most (in a certain sense) cases of a conjecture of Alon on the existence of Ramsey Cayley graphs. Alon \cite{AlonConjStatement,AlonConjPaper} conjectured that, for some absolute constant $C$, every finite group $G$ admits a Cayley graph with no clique or independent set of size exceeding $C\log\abs{G}$, and we confirm this (see \cref{cor:ramsey-cayley-almost-all}) for all groups of almost all orders. Our work builds on and generalizes recent work of Conlon, Fox, Pham, and Yepremyan \cite{CFPY}: we remove the abelian condition from many of their results while maintaining a similar quantitative strength. Our second application is to give a local version of Roth's theorem in many non-abelian groups, including settings where no modeling lemma exists to transform small subsets of groups to dense subsets of other groups (see \cite{GreenModeling,HegyvariHennecart}). For example, we show (see \cref{cor:roth-odd}) that, in odd-order groups, sets with no nontrivial three-term arithmetic progression have at least quasi-polynomial doubling. This strengthens a result of Freiman \cite{Freiman} in the integer setting. 

\subsection{Application: Ramsey Cayley graphs}\label{sec:applications}

A classical result of Erd\H{o}s and Szekeres \cite{ErdosSzekeres}, strengthening Ramsey's theorem, states that any graph on $n$ vertices must have a clique or independent set of size at least $\frac12\log_2n$. Motivated by this, a graph on $n$ vertices is termed \emph{$C$-Ramsey} for some $C>1/2$ if it has no clique or independent set of size larger than $C\log_2n$. Erd\H{o}s \cite{ErdosRamseyNum} showed that almost all graphs are $2$-Ramsey. However, no explicit constructions of $C$-Ramsey graphs are known for any constant $C$.

One class of graphs within which to search for $C$-Ramsey graphs are \emph{Cayley graphs}, graphs arising from finite groups. Given a group $G$ and a symmetric subset $S\subset G$, we let $\Cay(G,S)$ denote the graph with vertex set $G$ and edge set $\{(x,y):x^{-1}y\in S\}$. Cayley graphs form a rich source of lower-bound examples in Ramsey theory. Following the idea that Cayley graphs have good Ramsey-theoretic properties, and motivated by applications to coding theory, Alon made the following conjecture.

\begin{conjecture}[\cite{AlonConjStatement}, {\cite[Conjecture~4.1]{AlonConjPaper}}]\label{conj:alon} There exists some absolute constant $C>0$ for which every finite group has a $C$-Ramsey Cayley graph. 
\end{conjecture}

Towards \cref{conj:alon}, Alon and Orlitsky \cite{AlonOrlitsky} showed that every group of order $n$ has a Cayley graph with clique and independence number $O(\log^2n)$. This was recently improved by Conlon, Fox, Pham, and Yepremyan \cite{CFPY} to $O(\log n\log\log n)$. 

When one restricts focus to special classes of groups, more can be said. The earliest such result is due to Green \cite{Green}, who proved \cref{conj:alon} for cyclic groups. Conlon, Fox, Pham, and Yepremyan proved Alon's conjecture for a much larger family of groups. Their most general result along these lines is the following. 

\begin{theorem}[{\cite[Theorem~4.6]{CFPY}}]\label{thm:CFPY-conclusion} There exists an absolute constant $C>0$ for which the following holds: if $G$ is a group of order $N$ with an abelian subgroup $H$ whose order is (i) at least $N/(\log N)^{1/1000}$ and (ii) relatively prime to $6$, then $G$ has a $C$-Ramsey Cayley graph.    
\end{theorem}

While \cref{thm:CFPY-conclusion} covers \emph{abelian} groups in quite wide generality, most non-abelian groups (even those of order relatively prime to $6$) fail to have abelian subgroups large enough for \cref{thm:CFPY-conclusion} to apply. Our main result towards \cref{conj:alon} is the following extension of \cref{thm:CFPY-conclusion} to a fully non-abelian setting. 

\begin{theorem}\label{thm:ramsey-cayley} Let $N$ be a sufficiently large positive integer, and suppose that the largest factor of $N$ which is a product of powers of $2$ and $3$ is at most $\exp(\exp((\log\log\log N)^{1/3})))$. Then every group $G$ of order $N$ has a $24000$-Ramsey Cayley graph.
\end{theorem}

We point out the following two corollaries of \cref{thm:ramsey-cayley}, the second of which generalizes \cite[Corollary~1.11]{CFPY} by removing the abelian condition.

\begin{corollary}\label{cor:ramsey-cayley} \cref{conj:alon} holds for the family of groups of order relatively prime to $6$.
\end{corollary}

\begin{corollary}\label{cor:ramsey-cayley-almost-all} For almost all\footnote{``Almost all'' here is in the sense of natural density.} integers $N$, every group of order $N$ has a $24000$-Ramsey Cayley graph.
\end{corollary}

We now say a few words about how our proof differs from the work \cite{CFPY}. The construction given in \cite{CFPY} to prove \cref{thm:CFPY-conclusion} is random but not uniformly random. The authors choose a Cayley graph $\Gamma:=\Cay(G,S)$ from a carefully constructed probability distribution such that (1) if $\abs{AA^{-1}}$ is large compared to $\abs{A}$, then $\Gamma$ is not too likely to contain $A$ as a clique or independent set, and (2) if $\abs{AA^{-1}}$ is small compared to $\abs{A}$, then $\Gamma$ \emph{deterministically} cannot have $A$ as a clique or independent set. It is only part (2) that, in their work, uses the abelian assumption. In their construction, sets which may be cliques or independent sets avoid nontrivial instances of the pattern\footnote{We are simplifying the situation slightly for the purposes of exposition; the actual pattern is $x_1y_1^{-1}=(x_2y_2^{-1})^2=(x_3y_3^{-1})^4=(x_4y_4^{-1})^8$. The pattern $xy^{-1}=(zw^{-1})^2$ satisfies many of the same properties as the more complicated pattern at hand, and the problem can be reduced to studying this simpler pattern in model cases like $(\ZZ/5\ZZ)^n$.} $xy^{-1}=(zw^{-1})^2$. In the abelian setting (written additively), this equation becomes $x+2w=y+2z$, and solutions can be found using the Pl\"unnecke--Ruzsa inequality and the pigeonhole principle. In the non-abelian setting, the manipulations used to reduce $xy^{-1}=(zw^{-1})^2$ to such a pattern are invalid.

Our construction of Cayley graphs in \cref{thm:ramsey-cayley} is a natural non-abelian analogue of that in \cite{CFPY} (in particular, it is also a randomized construction). Using \cref{thm:solv-abelian-struct,thm:gen-abelian-struct}, we are able to obtain (a slightly weaker version of) their part (2) by simply restricting our focus to solutions to $xy^{-1}=(zw^{-1})^2$ where $\{x,y,z,w\}$ all lie in some coset of an abelian subgroup. Within such a coset, the argument in \cite{CFPY} applies and the relevant pattern can be easily found.

By pushing our argument somewhat further, we are able to obtain the following strengthening of \cref{thm:ramsey-cayley}, which pins down precisely which subgroups form obstructions to finding Ramsey Cayley graphs.

\begin{theorem}\label{thm:ramsey-cayley-general} The following holds for all sufficiently large $N$. Let $G$ be a group of order $N$ and let $r=\floor*{\exp((\log\log\log N)^{1/3})}$. Suppose that $G$ possesses neither
\begin{itemize}
    \item a subgroup isomorphic to $(\ZZ/2\ZZ)^r$ or $(\ZZ/3\ZZ)^r$, nor
    \item a subgroup isomorphic to $\ZZ/L\ZZ$, where
    \[L=\lcm\pr*{1,2,\ldots,\floor*{\exp\pr*{(\log\log N)^{1/27}}}}.\]
\end{itemize}
Then $G$ has a $24000$-Ramsey Cayley graph.
\end{theorem}

\subsection{Application: Growth of \texorpdfstring{$3$}{3}-AP-free sets}

In his monograph on set addition, Freiman \cite[Theorem~2.30]{Freiman} gave the following generalization of Roth's celebrated theorem on $3$-term progressions: if $A\subset\ZZ$ avoids nontrivial $3$-term arithmetic progressions, then $\abs{A+A}/\abs{A}$ must grow without bound as $\abs{A}$ grows. Ruzsa \cite{Ruzsa} strengthened Freiman's result to the explicit asymptotic $\abs{A+A}\geq c\abs{A}(\log\abs{A})^c$ for some constant $c>0$. Using the recent breakthrough of Kelley and Meka \cite{KelleyMeka}, This result can be strengthened to the statement that $\abs{A+A}/\abs{A}$ grows at least quasi-polynomially in $\abs{A}$. (This was observed by Schoen \cite{Schoen}, who also strengthened Ruzsa's results in the case of longer progressions). Such a statement has applications to the study of restricted sumsets \cite{SchoenRest,HHP}, as well as to various problems in discrete geometry \cite{Pach,PW,Dumitrescu}. 

There are two different ways to generalize $3$-term arithmetic progressions to non-abelian groups. The first is the pattern $\{x,xt,xt^2\}$, which we call the \emph{translational} $3$-AP; the second is a solution to $xz=y^2$, which we call the \emph{averaging} $3$-AP. The translational $3$-AP, as it is invariant under left-translation, tends to be easier to handle. In particular, the ``working within cosets of an abelian group'' strategy reduces the translational $3$-AP problem to finding $3$-APs in abelian groups (see \cite[Theorem~1.4]{SandersRoth}, where this observation is credited to Croot), while it reduces the averaging $3$-AP problem instead to finding \emph{corners} in abelian groups (see \cite[Corollary~1.6]{JLLOS}, where this observation is credited to Fox). Conversely, arguments of Petrov \cite{Petrov} suggest that the averaging $3$-AP may be more natural in groups of bounded exponent, from the perspective of the polynomial method.

As our second application of \cref{thm:solv-abelian-struct,thm:gen-abelian-struct}, we provide an extension of Freiman's result to non-abelian groups, obtaining quasi-polynomial bounds.

\begin{theorem}\label{thm:local-3-AP} There is an absolute constant $c>0$ for which the following holds. Let $G$ be a group, and let $A\subset G$ be a set with no nontrivial translational three-term arithmetic progression, i.e.~that if $x,xy,xy^2\in A$, then $y=1$.\footnote{The literature is not consistent on the definition of ``nontrivial'' for three-term arithmetic progressions in groups which may have $2$-torsion. Following \cite{SandersRoth} (but not \cite{Solymosi}, for example), we take a minimal approach, classifying an arithmetic progression as nontrivial as long as its elements are not all equal. In particular, we allow $xy^2=x$ in our definition of a ``nontrivial'' translational $3$-AP $x,xy,xy^2$.} Suppose that either
\begin{enumerate}[(a)]
    \item $G$ is solvable, or
    \item $G$ is a finite group which does not possess $(\ZZ/2\ZZ)^d$ as a subgroup, where $d=\floor*{\exp\pr*{c\sqrt{\frac{\log\log\abs{A}}{\log\log\log\abs{A}}}}}$.
\end{enumerate}
Then
\[\frac{\abs{AA^{-1}}}{\abs{A}}\geq\exp\pr*{c(\log\abs{A})^{1/1305}}.\]
If $A$ instead avoids averaging three-term arithmetic progressions, i.e.~solutions to $xz=y^2$ with $x\neq y$, then under the same conditions on $G$ we obtain
\[\frac{\abs{A^2}}{\abs{A}}\geq\exp\pr*{c(\log\abs{A})^{1/64809}}.\]
\end{theorem}

Unlike in the abelian setting, our arguments do not go directly by finding a Freiman-isomorphic copy of $A$ which is dense in an ambient finite group. Indeed, no such ``dense model'' can exist in general for non-abelian sets of small doubling; see \cite[Theorem~1.2]{GreenModeling} and \cite[Theorem~1]{HegyvariHennecart}. Instead, we first use \cref{thm:solv-abelian-struct,thm:gen-abelian-struct} to reduce our problem to a subset of $A$ which lives in a coset (or, in the averaging case, a double coset) of an abelian group, and then appeal to results from the abelian world which do require modeling lemmas. 

We remark that the conclusion in the translational case is (qualitatively, not just quantitatively) stronger than that of the averaging case: lower bounds on $\abs{AA^{-1}}/\abs{A}$ imply lower bounds on $\abs{A^2}/\abs{A}$ by Ruzsa's triangle inequality, but the converse does not hold. We are able to obtain a lower bound on $\abs{AA^{-1}}$ in the translational case, but not in the averaging case, precisely because translational $3$-APs are preserved under left-translation. In any case, Ruzsa's triangle inequality furnishes the following attractive corollary.

\begin{corollary}\label{cor:roth-odd} There exists an absolute constant $c$ for which the following holds. Let $G$ be a finite group of odd order,\footnote{To read this corollary from \cref{thm:local-3-AP}, one either needs to use statement (a) combined with the Feit--Thompson theorem or to use statement (b), which relies on the classification of finite simple groups. The former is thus more elementary, and in fact gives slightly better quantitative bounds than what we have written here if worked out in full.} and let $A\subset G$. Suppose either that $A$ contains no translational $3$-APs or that $A$ contains no averaging $3$-APs. 
\[\frac{\abs{A^2}}{\abs{A}}\geq\exp\pr*{c(\log\abs{A})^{1/64809}}.\]
\end{corollary}

We use, as an input in our proof of the translational $3$-AP portion of \cref{thm:local-3-AP}, the aforementioned result by Kelley and Meka \cite{KelleyMeka} in the abelian setting. In the averaging $3$-AP portion, we use instead the recent substantial improvement of Jaber, Liu, Lovett, Ostuni, and Sawhney \cite{JLLOS} on the corners problem. 

\subsection{Approximate subgroups of matrix groups}\label{sec:intro-approx}

In proving \cref{thm:gen-abelian-struct}, much of our work will be done within the composition series of the ambient group $G$. We thus need to use the classification of finite simple groups, which allows us to reduce many of our questions to questions about subsets of \emph{linear groups} with small doubling. The structure of such subsets is quite well-studied, particularly in the context of finding \emph{diameter bounds} for Cayley graphs on groups. Our work on linear groups, while drawing on the methods used to study such problems, is entirely self-contained.

Our main result on linear groups is the following, which parallels \cref{thm:solv-abelian-struct}.

\begin{theorem}\label{thm:GL} There exists an absolute constant $C>0$ for which the following holds. Let $\FF$ be any field, let $d\geq1$ be any integer, and let $A\subset\GL_d(\FF)$ be a $K$-approximate group for some $K\geq 1$. Let
\[M(d):=\exp\pr*{C\log^2d\cdot\log(1+\log d)}=e^{\log^{2+o(1)}(d)}\]
be a quasi-polynomial function of $d$. Then there exists some abelian subgroup $H\leq\GL_d(\FF)$ for which
\[\abs{A^2\cap H}\geq\frac1{(2K)^{Cd^9}}\abs{A^2}^{1/M(d)}.\]
\end{theorem}

The source of the $\log(1+\log d)$ term in the statement of \cref{thm:GL} is the following proposition, which may be of independent interest.

\begin{proposition}\label{prop:reg-el} Let $d\geq 8000$, let $\FF$ be an algebraically closed field, and let $A\subset\GL_d(\FF)$ be a $K$-approximate subgroup. Then one of the following holds:
\begin{itemize}
    \item (Not-too-irregular element) There exists some element $a$ of $A^2$ whose largest-dimension eigenspace has codimension at least $\frac d{600\log\log d}$.

    \item (Control by well-behaved subgroup) There exist two subspaces $V_1\subset V_2\subset\FF^d$ with $\dim V_2>\frac{9d}{10}$ and $\dim V_1<\frac d{10}$ for which at least $(2K)^{-2d^6}\abs{A^2}$ elements $a$ of $A^2$ satisfy the following property: there exists some $\lambda\in\FF$ for which $(a-\lambda)V_2\subset V_1$.
\end{itemize}
\end{proposition}

Up to an $O(1/\log\log d)$ factor in the first condition, \cref{prop:reg-el} is optimal: the invertible operators $a$ for which $(a-\lambda)V_2\subset V_1$ for some $\lambda$ form a subgroup of $\GL_d(\FF)$, and every element of this subgroup has $\lambda$-eigenspace of dimension at least $4d/5$. So, the ``control by a well-behaved subgroup'' condition is a genuine obstruction to $A$ having some element with no large eigenspace.

We now comment on the relationship between our work and previous work on approximate subgroups of linear groups. Our comments here are restricted to the relationship between the results; for a deeper discussion of the relevant methods, see \cref{sec:GL-history}.

A conjecture of Babai \cite{BabaiSeress} states that every connected Cayley graph on a simple group of order $N$ has diameter $(\log N)^{O(1)}$; in other words, if $G$ is a simple group of order $N$ and $A\subset G$ generates $G$, then $A^D=G$ for some $D=(\log N)^{O(1)}$. Much of the work towards this conjecture, such as the resolution \cite{HelfgottSL2,HelfgottSL3,BGTLinear,PyberSzabo} of the bounded rank case, has focused on  ``product theorems'' for various groups: statements roughly of the form ``if $A\subset G$ generates $A$, then either $A^3=G$ or $A^3$ is much larger than $A$.'' This is a statement much stronger than what we need (assuming the additional assumption that $A$ generates some well-behaved linear group), but the quantitative dependence on the rank in all such works tends to be quite poor. Recent quantitative improvements \cite{BDH}, as well as work in the ``high rank'' case \cite{HelfgottSeress,HMPQ}, prove expansion results for very large iterated product sets, rather than simply of the form ``$A^3$ is much larger than $A$.'' These results, while quite useful for diameter bounds, are not of much help in our setting. 

There has also been interesting work \cite{EMPS} towards a classification approximate subgroups of linear groups (in the line of the structure theorem of Breuillard, Green, and Tao \cite{BGTGeneral}) without assumptions on the group they generate. It seems plausible to us that the results of \cite{EMPS} imply a weakening of \cref{thm:GL} where one allows the function $M(d)$ to grow arbitrarily quickly. Relatedly, if one is willing to look for not-too-irregular elements in $A^m$ instead of $A^2$ for some $m$ roughly exponential in $d$, an argument similar to that of \cite[Proposition~6.2]{BDH} gives a stronger version of \cref{prop:reg-el}.

In this work, we strive to optimize the quantitative dependence on the rank of the ambient linear group. This allows the generality in our application to Alon's conjecture (see \cref{thm:ramsey-cayley}) to be comparable to that of \cite{CFPY} in the abelian setting.

\subsection{Organization} In \cref{sec:prelim} we present some definitions, notation, and preliminary results which will be useful to us. 

\cref{sec:approach} contains an outline of our proofs of \cref{thm:solv-abelian-struct,thm:gen-abelian-struct}, with many of the details deferred to \cref{sec:commutable-proofs}. 

\crefrange{sec:GL-setup}{sec:subspace-to-subgroup} are concerned with approximate subgroups of matrix groups; this is where we prove \cref{thm:GL} and \cref{prop:reg-el}. 

In \cref{sec:roth} we apply our results to Roth-type problems, proving \cref{thm:local-3-AP}, while in \cref{sec:ramsey-cayley} we apply our results to the existence of Ramsey Cayley graphs, proving \cref{thm:ramsey-cayley,thm:ramsey-cayley-general}. Both \cref{thm:ramsey-cayley-general,thm:local-3-AP} require some results in abelian additive combinatorics, which we state when needed and prove in \cref{sec:abelian}. For a reader only interested in the application of \cref{thm:solv-abelian-struct,thm:gen-abelian-struct} to three-term arithmetic progressions or to Alon's conjecture, only \crefrange{sec:roth}{sec:abelian} are relevant.

\cref{sec:conclusion} contains some concluding remarks and conjectures. Finally, proofs of some purely group-theoretic results are deferred to the appendices.

\section{Preliminaries}\label{sec:prelim}

We will need a few (mostly standard) properties of and definitions relating to groups.

\subsection{Groups and subsets of groups}\label{sec:prelim-gp} For a group $G$ and a subset $S\subset G$, we write $S\leq G$ to denote that $S$ is a subgroup of $G$ and $S\unlhd G$ to denote that $S$ is a normal subgroup of $G$. We write $\langle S\rangle$ for the subgroup of $G$ generated by $S$. A quotient of a subgroup of $G$ is termed a \emph{subquotient} of $G$. 

For a group $G$, we write $1_G$ for the identity element of $G$ and $Z(G)$ for the center of $G$. For subgroups $H_1,H_2\leq G$, we write $[H_1,H_2]$ for the group generated by the set of commutators $\{h_1h_2h_1^{-1}h_2^{-1}:h_1\in H_1,h_2\in H_2\}$. For an element $g\in G$, we write
\begin{align*}
C_G(g)&:=\set{h\in G:gh=hg}\\
\Conj_G(g)&:=\set{hgh^{-1}:h\in G}
\end{align*}
for the centralizer and conjugacy class of $g$ in $G$, respectively. For $S\subset G$, we also write $C_G(S)$ for the intersection $\bigcap_{g\in S}C_G(g)$. In all notation, we frequently omit the subscript specifying the group when it is clear from context.

\begin{definition}[Subnormal and normal series] Given a group $G$, a sequence of subgroups $G=G_0\unrhd G_1\unrhd \cdots\unrhd G_s=\{1\}$, each of which is normal in the previous, is termed a \emph{subnormal series}. If each $G_i$ is in fact normal in $G$, then the series is a \emph{normal series}. In either case, we call the quotient groups
\[G_0/G_1,G_1/G_2,\ldots,G_{s-1}/G_s\]
the \emph{successive quotients} of the series. 
\end{definition}

There are two particularly important series one can form from a group.

\begin{definition}[Derived series] Given a group $G$, its \emph{derived series} $G^{(0)}\unrhd G^{(1)}\unrhd\cdots$ is defined by $G^{(0)}=G$ and $G^{(i+1)}=[G^{(i)},G^{(i)}]$. Suppose that the derived series of $G$ is eventually constant, and let $s$ be minimal for which $G^{(s)}=G^{(s+1)}$. If $G^{(s)}$ is the trivial group, we say $G$ is \emph{solvable} of \emph{derived length} $s$. In this case the derived series is a normal series of $G$.
\end{definition}

\begin{definition}[Composition series] A group is \emph{simple} if it is nontrivial and has no proper normal subgroup. Given a group $G$, a \emph{composition series} of $G$ is a subnormal series $G=G_0\unrhd G_1\unrhd \cdots\unrhd G_s=\set{1}$ each successive quotient of which is simple. The simple groups $G_i/G_{i+1}$ are termed \emph{composition factors} of $G$.
\end{definition}

\noindent Every finite group has a composition series, as do some infinite groups. The term ``composition factor'' is justified by the Jordan--H\"older theorem, which states that (when $G$ has a composition series) every composition series of a group $G$ has the same composition factors with the same multiplicities. We also record here that a finite group $G$ is solvable if and only if all of its composition factors are abelian. (Note that the only abelian finite simple groups are $\ZZ/p\ZZ$ for primes $p$.)

If we impose in the definition of a composition series the condition that each $G_i$ is normal in $G$, we obtain instead a \emph{chief series}.

\begin{definition}[Chief series] Given a group $G$, a \emph{chief series} of $G$ is a normal series $G=G_0\unrhd G_1\unrhd\cdots\unrhd G_s=\set{1}$ where, for each $i$, the only normal subgroups $N$ of $G$ satisfying $G_i\geq N\geq G_{i+1}$ are $G_i$ and $G_{i+1}$.
\end{definition}

\noindent We will use the following relationship between chief series and composition series.

\begin{lemma}[{see \cite[Theorems~I.9.12~and~I.11.7]{Huppert}}]\label{lem:chief-series} Every finite group $G$ possesses a chief series. Moreover, in any chief series $G_0\unrhd\cdots\unrhd G_s$ of a finite group $G$, the quotient groups $G_i/G_{i+1}$ are each direct products of isomorphic simple groups, each of which are composition factors of $G$.
\end{lemma}

We also need a few facts following from the classification of finite simple groups. We begin with the following nonstandard definition.

\begin{definition} Let $G$ be a group. The \emph{composition rank} of $G$ is the smallest positive integer $d$ such that every composition factor $S$ of $G$ embeds into $\GL_d(\FF)$ for some field $\FF$ (which may depend on $S$). If $G$ does not have a composition series, or some composition factor of $G$ fails to embed in any $\GL_d(\FF)$, then the composition rank of $G$ is defined to be $\infty$.
\end{definition}

For example, solvable groups have composition rank one, since $\ZZ/p\ZZ$ embeds into $\GL_1(\CC)=\CC^\times$ for every prime $p$. The following lemma is partially a corollary of the classification of finite simple groups. Qualitatively, it says that having large composition rank is equivalent to having a large alternating group as a subquotient.

\begin{lemma}\label{lem:comp-rank-alternating} There exists an absolute constant $C>0$ for which the following holds. Let $d$ be a positive integer and $G$ be a finite group.
\begin{enumerate}[(1)]
    \item If $G$ possesses a subquotient isomorphic to $\A_d$, then $G$ has composition rank at least $d/7-1$.

    \item If $G$ possesses no subquotient isomorphic to $\A_d$, then $G$ has composition rank at most $Cd^2$.
\end{enumerate}
\end{lemma}

We defer the proof of \cref{lem:comp-rank-alternating} to \cref{sec:appendix-cfsg}. For some of our applications, we will also need another purely group-theoretic lemma. This lemma essentially follows from results of Thompson (see \cite[Theorem~III.12.3]{Huppert}) and Mann \cite[Theorem~B]{Mann}, along with some simple lemmas reducing to the case of $p$-groups (which are also useful in the proof of \cref{lem:comp-rank-alternating}(1)). We give a derivation in \cref{sec:appendix-cfsg}.

\begin{lemma}\label{lem:sectional-p-rank} Let $p$ be a prime and $d$ be a positive integer. Every finite group with a subquotient isomorphic to $(\ZZ/p\ZZ)^{2d^2}$ has a subgroup isomorphic to $(\ZZ/p\ZZ)^d$.
\end{lemma}

\subsection{Arithmetic combinatorics}\label{sec:prelim-arith}

We now give some notation and terminology which is more specialized to arithmetic combinatorics. Given a group $G$ and subsets $A,B\subset G$, we write $AB:=\{ab:a\in A,b\in B\}$ for the set of products of an element of $A$ and an element of $B$. For a positive integer $k$, we write $A^k:=\{a_1a_2\cdots a_k:a_1,\ldots,a_k\in A\}$ for the $k$th iterated product set of $A$, and we write $A^{-1}:=\{a^{-1}:a\in A\}$ for the set of inverses of elements of $A$. We sometimes use the corresponding additive notation ($A+B$, $kA$, and $-A$) when the ambient group is abelian; such use is restricted to \cref{sec:abelian}.

A handful of tools from additive combinatorics apply without modification to the non-abelian setting. Two such tools are Ruzsa's triangle inequality and covering lemma (see for example \cite[Lemmas~2.3.4~and~2.4.4]{Tointon}).

\begin{lemma}[Ruzsa's triangle inequality]\label{lem:tri-ineq} For any three finite sets $A,B,C\subset G$, we have
\[\abs{A}\abs{B^{-1}C}\leq\abs{AB}\abs{AC}.\]
\end{lemma}

\begin{lemma}[Ruzsa's covering lemma]\label{lem:covering} Let $A,B\subset G$ be finite sets. There is a set $Y\subset G$ of size at most $\abs{AB}/\abs{B}$ for which $A\subset YBB^{-1}$.
\end{lemma}

Finally, we need the notion of a \emph{dissociated} set, which will be essential in building the sets described in the conclusion of \cref{thm:solv-abelian-struct,thm:gen-abelian-struct}. Dissociated sets are usually only defined in abelian groups; all of our dissociated sets in arbitrary groups will be contained in abelian subgroups.

\begin{definition} Let $G$ be a group and $D=\set{h_1,\ldots,h_m}\subset G$ be a commuting subset (that is, every pair of elements of $D$ commute; this is equivalent to the statement that $\langle D\rangle\leq G$ is an abelian subgroup of $G$). We say that $D$ is \emph{dissociated} if, for any $\eps_1,\ldots,\eps_m\in\set{-1,0,1}$ not all zero, we have
\[h_1^{\eps_1}h_2^{\eps_2}\cdots h_m^{\eps_m}\neq 1.\]
This is equivalent to the distinctness of the $2^m$ subset products $h_1^{\delta_1}\cdots h_m^{\delta_m}$ for $(\delta_1,\ldots,\delta_m)\in\set{0,1}^m$. Given a dissociated set $D=\set{h_1,\ldots,h_m}$, we write
\[\dspan{D}:=\set[\big]{h_1^{\eps_1}h_2^{\eps_2}\cdots h_m^{\eps_m}:(\eps_1,\ldots,\eps_m)\in\set{-1,0,1}^m}.\]
\end{definition}

\subsection{Approximate groups}\label{sec:prelim-approx} We have stated \cref{thm:solv-abelian-struct,thm:gen-abelian-struct} in terms of sets with small tripling. These are closely related to the more structured (and well-studied) \emph{approximate groups}. 

\begin{definition} Let $K\geq 1$ be a parameter. A subset $A$ of a group $G$ is \emph{symmetric} if $A=A^{-1}$ and \emph{centered} if additionally $1\in A$. We call $A$ a \emph{$K$-approximate group} if $A$ is centered and there exists a subset $X\subset G$ of size at most $K$ for which $A^2\subset XA$.
\end{definition}

Iterated product sets of $K$-approximate groups are extremely well-behaved. In particular, if $A$ is an approximate group and $k\geq j$ are positive integers, then $A^k\subset X^{k-j}A^j$, and so $\abs{A^k}\leq K^{k-j}\abs{A^j}$. We will need some other properties of (not necessarily abelian) approximate groups. The following fact, which can be found in Tointon's book \cite{Tointon} and of which we will make repeated use, says that intersections of approximate groups with subgroups are well-behaved.

\begin{lemma}[{Special case of \cite[Proposition~2.6.5]{Tointon}}]\label{lem:intersect-subgroup} Let $G$ be a group, and let $H$ be any (not necessarily normal) subgroup of $G$. If $A\subset G$ is any $K$-approximate subgroup, then $A^2\cap H$ is a $K^3$-approximate group. Moreover, for any $k\geq 2$, we have $\abs{A^k\cap H}\leq K^{k-1}\abs{A^2\cap H}$.
\end{lemma}

We also need the following result of Tao, which says that sets $A\subset G$ with $AA^{-1}$ small are ``controlled'' by a well-behaved centered set.

\begin{lemma}[{\cite[Proposition~4.5]{Tao08}}]\label{lem:nice-set-finder}
If $\abs{AA^{-1}}\leq K\abs{A}$, then there exists a centered set $B\subset A^{-1}A$ satisfying $\abs{B}\geq\abs{A}/(2K)$ and $\abs{AB^nA^{-1}}\leq 2^nK^{2n+1}\abs{B}$ for every positive integer $n$.\footnote{This result, as we state it here, contains two more conditions than the version stated in \cite{Tao08}, namely that $B\subset A^{-1}A$ and that $1\in B$. Both conditions can be easily read from the proof.}
\end{lemma}

We will need one lemma relying on the following definition. Given two groups $G$ and $G'$, a centered subset $A\subset G$, and an integer $\ell\geq 2$, a \emph{centered Freiman $\ell$-homomorphism} $\phi\colon A\to G'$ is a map for which $\phi(1_G)=1_{G'}$ and, if $x_1,\ldots,x_\ell,y_1,\ldots,y_\ell\in A$ satisfy $x_1\cdots x_\ell=y_1\cdots y_\ell$, then $\phi(x_1)\cdots\phi(x_\ell)=\phi(y_1)\cdots\phi(y_\ell)$. We record two properties of this definition: first, any centered Freiman $\ell$-homomorphism satisfies $\phi(x^{-1})=\phi(x)^{-1}$ for each $x\in A$; second, any group homomorphism is (when restricted to any centered set) a centered Freiman $\ell$-homomorphism for any $\ell$. 

The following lemma says that, when viewed in the right way, the fibers of a Freiman homomorphism between approximate groups are not too different. In nearly all of our applications of this lemma, $\phi$ will be the projection onto a quotient of $G$.

\begin{lemma}\label{lem:intersect-subgroup-projection} Let $G$ be a group and let $A$ be a $K$-approximate subgroup of $G$. Let $G'$ be a group, and let $\phi\colon A^6\to G'$ be a map which restricts to a centered Freiman $3$-homomorphism on $A^2$. Then, for every subset $S\subset G'$, 
\[K^{-4}\frac{\abs{A^2\cap\phi^{-1}(S)}}{\abs{A^2}}\leq\frac{\abs{\phi(A^2)\cap S}}{\abs{\phi(A^2)}}\leq \frac{\abs{A^6\cap\phi^{-1}(S)}}{\abs{A^2}}.\]
\end{lemma}

\begin{proof} Select $t\in\phi(A^2)$ to maximize $\abs{A^2\cap \phi^{-1}(t)}$. For each $s\in \phi(A^2)$, pick elements $a_s,b_s\in A$ satisfying $s=\phi(a_sb_s)$. For any such $s$, and any $a'\in A^2$ for which $\phi(a')=t$, we have, using our assumptions on $\phi$, that
\[\phi((a_tb_t)^{-1}a'(a_sb_s))=\phi(a_tb_t)^{-1}\phi(a')\phi(a_sb_s)=t^{-1}ts=s.\]
We conclude that, for each $s\in\phi(A^2)$,
\begin{equation}\label{eq:cap-phi}
\abs{A^6\cap\phi^{-1}(s)}\geq\abs{(a_tb_t)^{-1}A^2a_sb_s\cap\phi^{-1}(s)}\geq\abs{A^2\cap\phi^{-1}(t)}\geq\frac{\abs{A^2}}{\abs{\phi(A^2)}},
\end{equation}
where the last inequality is by our choice of $t$. Summing \eqref{eq:cap-phi} over all $s\in\phi(A^2)\cap S$ gives the second inequality in our statement. To prove the first, we sum \eqref{eq:cap-phi} over all $s\in\phi(A^2)$ to obtain
\[\abs{A^6}\geq\sum_{s\in\phi(A^2)}\abs{A^6\cap\phi^{-1}(s)}\geq\abs{\phi(A^2)}\cdot\abs{A^2\cap\phi^{-1}(t)}.\]
As a result,
\[\frac{\abs{A^2\cap\phi^{-1}(S)}}{\abs{A^6}}=\frac1{\abs{A^6}}\sum_{s\in\phi(A^2)\cap S}\abs{A^2\cap\phi^{-1}(s)}\leq\frac{\abs{A^2\cap\phi^{-1}(t)}}{\abs{A^6}}\cdot\abs{\phi(A^2)\cap S}\leq\frac{\abs{\phi(A^2)\cap S}}{\abs{\phi(A^2)}}.\]
The inequality $\abs{A^6}\leq K^4\abs{A^2}$ finishes the proof.
\end{proof}

The last result we will need is the following proposition of Sanders, related to an ``almost-periodicity'' lemma of Croot and Sisask. If $A\subset G$ is a centered set with small doubling, then this result allows us to find a large iterated product set contained in $A^4$.

\begin{proposition}[{\cite[Proposition~4.1]{Sanders}}]\label{prop:sanders-quasipoly} There exists an absolute constant $C\geq1$ for which the following holds. Let $G$ be any group, let $A\subset G$ satisfy $\abs{AA^{-1}}\leq K\abs{A}$, and let $t\geq 1$ be a parameter. Then there exists some centered set $S\subset G$ such that
\begin{enumerate}[(1)]
    \item $\abs{S}\geq\exp(-Ct^2(\log 2K)^2)\abs{A}$ and
    \item $S^t\subset A^{-1}AA^{-1}A$.
\end{enumerate}
\end{proposition}

For our proofs of \cref{thm:GL,prop:reg-el}, we will need some notation related to matrix groups and linear algebra. We defer this to just before it is needed, in \cref{sec:prelim-matrix}.

\subsection{Asymptotic notation} We write $f=O(g)$, or equivalently $f\lesssim g$ or $g=\Omega(f)$, if there exists some absolute constant $C$ for which $\abs{f}\leq Cg$ when all relevant parameters are sufficiently large. We write $f=\Theta(g)$ when both $f=O(g)$ and $g=O(f)$. Whenever asymptotic notation appears within an argument, any hidden constant is absolute, independent of any parameters. In the discussion surrounding an argument, we occasionally write an expression like $O_d(1)$ to mean ``bounded by a function of $d$, independently of any other parameters.''

\section{Our approach, and a toy problem}\label{sec:approach}

In this section, we describe the approach we use to prove \cref{thm:solv-abelian-struct,thm:gen-abelian-struct}. We begin with a toy problem: given a finite solvable group $G$, how large of an abelian subgroup must it contain? 

In the setting of general finite groups, this problem goes back to Erd\H{o}s and Straus \cite{ErdosStraus}, who first gave a lower bound of $(1-o(1))\log\abs{G}$. (Much earlier, Miller had shown that a finite \emph{$p$-group} $G$ contains an abelian subgroup of order $\exp(\Omega(\sqrt{\log\abs{G}}))$ \cite{Miller}.) The problem was essentially solved by Pyber \cite{Pyber} (see \cref{thm:Pyber}), who gave a lower bound of $\exp(\Omega(\sqrt{\log\abs{G}}))$ for general finite $G$. This cannot be improved, even in the special case where $G$ is $2$-step nilpotent \cite{Olshanskii,BGH}.

We present a different proof of a weak version of \cref{thm:Pyber}, based partially on the argument in \cite{Miller}, which obtains a quasi-polynomial lower bound on the size of the largest abelian subgroup of a finite solvable group $G$. We will then generalize this argument to approximate groups, beginning in \cref{sec:approach-outline}.

\subsection{Abelian subgroups in solvable groups}\label{sec:approach-genuine-gp}

\begin{proposition}[Weak version of {\cite[Theorem~1]{Pyber}}]\label{prop:solvable-gpvsn} Any finite solvable group $G$ of order $n$ contains an abelian subgroup of order at least
\[\exp\pr*{\frac12\log^{1/3}n}.\]
\end{proposition}

The proof we present of \cref{prop:solvable-gpvsn} is the starting point for our proof of \cref{thm:solv-abelian-struct}, and thus of the more general \cref{thm:gen-abelian-struct}. Every ingredient in this proof will transfer over to the approximate subgroup case, except for the following fact.

\begin{fact}\label{fact:expand-with-new-generator}
Let $G$ be a group. If $H$ is a subgroup of $G$ and $h\in G\setminus H$, then $\frac{\abs{\langle H,h\rangle}}{\abs H}\geq 2$.
\end{fact}

\noindent We will wish to apply \cref{fact:expand-with-new-generator} inside approximate groups instead of genuine groups, where it does not necessarily hold. This is where we use the notion of dissociated sets; see \cref{sec:approach-outline} for an elaboration of this fix.

Our proof of \cref{prop:solvable-gpvsn} proceeds via the following lemma (cf.~\cref{lem:abelian-commutable}).

\begin{lemma}\label{lem:solvable-gpvsn-step} Let $G$ be a group, $N$ a normal subgroup of $G$, and $I$ a subgroup of $Z(G)$. Suppose that $I\lneq G$ and that $G/N$ is abelian. Then there exists a positive integer $m$ and an abelian subgroup $I\leq I'\leq G$, generated by adjoining $m$ elements to $I$, satisfying $\abs{I'}/\abs{I}\geq 2^m$ and
\[\abs{G}\leq\abs{N}^m\abs{I'}.\]
\end{lemma}
\begin{proof} Let $h_1,\ldots,h_m\in G$ be a maximal list of elements of $G$ such that $h_1,\ldots,h_m$ commute and
\[I\lneq \langle I,h_1\rangle\lneq \cdots\lneq \langle I,h_1,h_2,\ldots,h_m\rangle.\]
(Such an element $h_1$ exists since $I\leq Z(G)$ but $I\lneq G$, so we have $m\geq 1$.) Define a function $\phi\colon G\to N^m$ by
\[\phi\colon x\mapsto (xh_1x^{-1}h_1^{-1},\ldots,xh_mx^{-1}h_m^{-1}).\]
(This maps $G$ into $N^m$ since $G/N$ is abelian.) If $x,y\in G$ are such that $\phi(x)=\phi(y)$, then $y^{-1}x$ commutes with $h_1,\ldots,h_m$. Thus, by the maximality of $m$, we must have $y^{-1}x\in\langle I,h_1,\ldots,h_m\rangle$. We conclude that the fibers of $\phi$ are exactly the left cosets of $\langle I,h_1,\ldots,h_m\rangle$. This implies
\[\abs{G}\leq\abs{N}^m\abs{\langle I,h_1,\ldots,h_m\rangle}.\]
Setting $I'=\langle I,h_1,\ldots,h_m\rangle$ and using \cref{fact:expand-with-new-generator} gives the result.
\end{proof}

We now give the rest of our proof of \cref{prop:solvable-gpvsn}. The ideas presented here will reappear in the proof of \cref{lem:series-commutable}.

\begin{proof}[Proof of \cref{prop:solvable-gpvsn}] Let $s$ be the derived length of $G$ and let
\[G=G^{(0)}\rhd G^{(1)}\rhd\cdots\rhd G^{(s)}=\set{1}\]
be the derived series of $G$. We construct our abelian subgroup $I$ of $G$ iteratively. We will define a sequence of abelian groups $\set{1}=:I_0\subset I_1\subset\cdots\subset I_\ell$ and an auxiliary sequence $m_0,m_1,\ldots,m_{\ell-1}$ of positive integers ($\ell$, the index at which our procedure terminates, is not chosen beforehand). For each $j$, when we define $I_j$, set $H_j:=C_G(I_j)$. We repeat the following process for each $j\geq 0$ in turn:
\begin{enumerate}[(i)]
    \item If $I_j=H_j$, i.e.~if $I_j$ is its own centralizer, terminate the procedure.

    \item Let $i$ be maximal such that $H_j\cap G^{(i)}\not\subset I_j$; such an $i$ exists since $I_j\lneq H_j$.

    \item Apply \cref{lem:solvable-gpvsn-step} to the group $H_j\cap G^{(i)}$ with the normal subgroup $H_j\cap G^{(i+1)}$ and the abelian subgroup $I_j\leq Z(H_j\cap G^{(i)})$ to obtain a positive integer $m$ and an abelian subgroup $I'\leq H_j\cap G^{(i)}$. Let $m_j:=m$ and $I_{j+1}:=I'$.
\end{enumerate}

We claim that
\begin{equation}\label{eq:H-ratio}
\frac{\abs{H_j}}{\abs{H_{j+1}}}\leq\abs{H_j\cap G^{(i)}}^{m_j}\leq\abs{I_j}^{m_j^2}\abs{I_{j+1}}^{m_j}
\end{equation}
Indeed, let $h_1,\ldots,h_{m_j}\in H_j\cap G^{(i)}$ be such that $I'=\langle I_j,h_1,\ldots,h_{m_j}\rangle$, and define a map $\sigma\colon H_j\to (H_j\cap G^{(i)})^{m_j}$ by
\[\sigma\colon x\mapsto (xh_1x^{-1},\ldots,xh_{m_j}x^{-1}).\]
Since $G^{(i)}$ is normal in $G$, we have that $H_j\cap G^{(i)}$ is normal in $H_j$, and so $xh_kx^{-1}\in H_j\cap G^{(i)}$ for each $x\in H_j$ and each $1\leq k\leq m_j$. If $\sigma(x)=\sigma(y)$, then $y^{-1}x$ commutes with $h_1,\ldots,h_{m_j}$. This implies that $y^{-1}x\in C_G(I_{j+1})=H_{j+1}$. We conclude the first inequality of \eqref{eq:H-ratio}. For the second, we also have from \cref{lem:solvable-gpvsn-step} and our choice of $i$ in step (ii) that
\begin{equation}\label{eq:H-grow}
\abs{H_j\cap G^{(i)}}\leq\abs{H_j\cap G^{(i+1)}}^{m_j}\abs{I_{j+1}}\leq\abs{I_j}^{m_j}\abs{I_{j+1}}.
\end{equation}
The second inequality of \eqref{eq:H-ratio} now follows from \eqref{eq:H-grow}.

What remains is simply computation. Since $I_\ell=H_\ell$, we have
\[n=\abs{G}=\abs{H_0}=\frac{\abs{H_0}}{\abs{H_\ell}}\cdot\frac{\abs{I_\ell}}{\abs{I_0}}=\prod_{j=0}^{\ell-1}\pr*{\frac{\abs{H_j}}{\abs{H_{j+1}}}\cdot\frac{\abs{I_{j+1}}}{\abs{I_j}}}.\]
Using \eqref{eq:H-ratio}, we get
\[n\leq\prod_{j=0}^{\ell-1}\abs{I_j}^{m_j^2-1}\abs{I_{j+1}}^{m_j+1}\leq\prod_{j=0}^\ell\abs{I_\ell}^{m_j^2+m_j}\leq\abs{I_\ell}^{2(m_0+\cdots+m_{\ell-1})^2}.\]
Taking the logarithm and using
\[\abs{I_\ell}=\prod_{j=0}^{\ell-1}\frac{\abs{I_{j+1}}}{\abs{I_j}}\geq 2^{m_0+\cdots+m_{\ell-1}}\]
gives
\[\log_2n\leq 2(m_0+\cdots+m_{\ell-1})^2\log_2\abs{I_\ell}\leq 2(\log_2\abs{I_\ell})^3.\]
The conclusion follows.
\end{proof}

\begin{remark} The bounds towards the end of the proof of \cref{prop:solvable-gpvsn} are rather crude. However, it is not possible to use \eqref{eq:H-ratio} and the bound $\abs{I_{j+1}}/\abs{I_j}\leq2^{m_j}$ alone to obtain a version of \cref{prop:solvable-gpvsn} with a bound exceeding $\exp(O(\log^{1/3}n))$. This may be seen by considering the possibility that $\ell=2$ and $m_0=m_1$.
\end{remark}

We now make some comments on the above proof with the hope that it may elucidate the structure of our argument. Our proof of \cref{prop:solvable-gpvsn} proceeds by first writing down a normal series of the ambient group $G$, then building up an abelian subgroup iteratively throughout this series. In each iteration, we carefully choose a step within this series, find as many elements as we can to add to our abelian group, and then restrict focus to the centralizer of the group we have constructed so far. This process can be done within an arbitrary subset $A$ of a group, as long as we know how to choose appropriately which step of the series to choose in each iteration of our process. The resulting asymptotics (particularly the analogues of \eqref{eq:H-ratio} and the lower bound on $\abs{I'}/\abs{I}$ in \cref{lem:solvable-gpvsn-step}) depend on how well-behaved the subset $A$ is. Elaborating on this observation will lead us to the proof of \cref{thm:solv-abelian-struct}.

A key observation in extending this method to non-solvable groups is that this iteration can be repeated. Suppose we have a group with a normal series whose successive quotients are not necessarily abelian. If we know how to pick commuting elements from the successive quotients, we can run the same argument as in our proof of \cref{prop:solvable-gpvsn}, replacing \cref{lem:solvable-gpvsn-step} with an arbitrary (potentially more complicated) procedure to find commuting elements.\footnote{This does not help us if we only know how to treat solvable groups, since any group with a normal series whose successive quotients are all solvable is in fact solvable. However, we will be able to pass to appropriately chosen abelian subgroups of our successive quotients with \cref{thm:GL} and use \cref{prop:solvable-gpvsn} in those subgroups.} We formalize the concept of such a procedure in \cref{def:commutable} below.

\subsection{Adapting to approximate groups}\label{sec:approach-outline}

We now use the proof strategy of \cref{prop:solvable-gpvsn} to prove corresponding results for approximate subgroups. Just as the proof of \cref{prop:solvable-gpvsn} involved first proving a statement about abelian quotient groups and then bootstrapping it to solvable groups, our proofs of \cref{thm:solv-abelian-struct,thm:gen-abelian-struct} will involve repeatedly bootstrapping statements about quotient groups which lie in certain classes. To this end, we introduce the following definition.

\begin{definition}\label{def:commutable} Let $R\geq 1$ and $\beta\in(0,1)$ be parameters. We say that a subgroup-closed class $\mathscr G$ of groups is \emph{$(R,\beta)$-commutable} if the following holds.

Let $\Gamma$ be a group, let $G\leq\Gamma$, and let $N\unlhd G$. Let $A$ be a $K$-approximate subgroup of $\Gamma$, and let $S$ be a centered subset of $\Gamma$. Suppose that $G/N\in\mathscr G$ and $A^4\cap N\subset S\subset C_\Gamma(G)$. Then the set $A^4\cap G$ contains a commuting dissociated subset $D$ satisfying $\dspan{D}\cap S=\set{1}$ and
\[\abs{D}\geq(\beta\log_3\abs{A^2\cap G})^{1/R}-\log_3(3K\abs{S}).\]
\end{definition}

Before stating some properties that \cref{def:commutable} satisfies, we make one comment (the details of which are borne out in the statement and proof of \cref{lem:commutable-to-abelian-substruct}). \cref{def:commutable} offers ways to find commuting dissociated sets with size polylogarithmic in $\abs{A}$, rather than commuting sets with size quasi-polynomial in $\abs{A}$ (which \cref{thm:solv-abelian-struct,thm:gen-abelian-struct} require). This is due to the aforemented failure of \cref{fact:expand-with-new-generator} in the approximate group setting. We will construct commuting sets by adding elements a few at a time. Without imposing additional restrictions on the elements we add, we have no guarantee that these elements allow us to expand our commuting subset of a bounded power of $A$ more quickly than linearly. (In the genuine subgroup setting, \cref{fact:expand-with-new-generator} guaranteed that the size of our commuting set at least doubled with each new element.) We circumvent this by imposing the restriction that our commuting set remain dissociated at each step. Such a condition does not cost us very much to implement, and it will enable us (later on) to boost a set of polylogarithmic size to a set of quasi-polynomial size: we can take the set of products of subsets of the dissociated set. There is one more wrinkle: we need enough of these products to lie in $A^4$. We ensure this by first applying Sanders' result, \cref{prop:sanders-quasipoly}, to the set $A$. At the cost of increasing the ``tripling constant'' $K$ to some larger $L$, this gives us an $L$-approximate group $B$ satisfying $B^t\subset A^4$ for some large $t$. Applying the condition of \cref{def:commutable} to $B$ rather than $A$, and then considering the set of products of elements of the resulting commuting dissociated set, will give us our results.

\vspace{2mm}

We will show that the notion of $(R,\beta)$-commutability enjoys the following four properties.

\begin{lemma}[Abelian groups are commutable]\label{lem:abelian-commutable} The class $\mathscr{A}$ of abelian groups is $(2,1/3)$-commutable.
\end{lemma}

\begin{lemma}[Groups with commutable successive quotients are commutable] \label{lem:series-commutable} Let $\mathscr H$ be any subgroup-closed class of groups, and let $\mathscr G$ denote the class of groups which possess finite normal series all of whose successive quotients lie in $\mathscr H$. Suppose that $\mathscr H$ is $(R,\beta)$-commutable. Then $\mathscr G$ is $(R+1,\beta')$-commutable with $\beta'=\min(\beta,\frac1{10})$.
\end{lemma}

\begin{lemma}[Groups with enough commutable subgroups are commutable]\label{lem:subgroup-commutable} Let $\mathscr H$ be any subgroup-closed class of groups. Let $C,\eps>0$ be parameters, and let $\mathscr G$ be a subgroup-closed class of groups for which, if $G\in\mathscr G$ and $A\subset G$ is a $K$-approximate group, then there exists a subgroup $H\leq G$ satisfying $H\in\mathscr H$ and
\[\abs{A^2\cap H}\geq\frac1{(2K)^C}\abs{A^2}^\eps.\]
Suppose that $\mathscr H$ is $(R,\beta)$-commutable. Then $\mathscr G$ is $(R,\beta')$-commutable with $\beta'=\beta\eps/(3C+12)$.
\end{lemma}

\begin{lemma}[Approximate subgroups of commutable groups have large abelian substructures]\label{lem:commutable-to-abelian-substruct} There exists an absolute constant $c>0$ for which the following holds. Let $\mathscr G$ be a subgroup-closed class of groups which is $(R,\beta)$-commutable, let $G\in\mathscr G$, and let $A\subset G$ be a finite symmetric subset satisfying $\abs{A^3}\leq K\abs{A}$. Then there exists some abelian subgroup $H$ of $G$ satisfying
\[\abs{A^4\cap H}\geq\exp\pr*{c\frac{(\beta\log\abs{A})^{1/(2R)}}{\log 2K}}\]
\end{lemma}

We will prove \crefrange{lem:abelian-commutable}{lem:commutable-to-abelian-substruct} in \cref{sec:commutable-proofs}. For now, we give some quick comments about their proofs. \cref{lem:abelian-commutable,lem:series-commutable} can be viewed as generalizations of the two main ingredients we used in the proof of \cref{prop:solvable-gpvsn}, and the proofs we will present of these lemmas will look similar to our proof of \cref{prop:solvable-gpvsn}. The proof of \cref{lem:subgroup-commutable} mostly involves unraveling definitions and using \cref{lem:intersect-subgroup-projection}. As previously discussed, the proof of \cref{lem:commutable-to-abelian-substruct} involves applying the commutability condition to an approximate subgroup arising from an application of \cref{prop:sanders-quasipoly}.

\subsection{Proofs of \texorpdfstring{\cref{thm:solv-abelian-struct,thm:gen-abelian-struct}}{Theorems 1.2 and 1.3}}\label{sec:approach-pfs}

To conclude this section, we present proofs of \cref{thm:solv-abelian-struct,thm:gen-abelian-struct} using \crefrange{lem:abelian-commutable}{lem:commutable-to-abelian-substruct}. To prove \cref{thm:gen-abelian-struct}, we also need to assume \cref{thm:GL}.

\begin{lemma}\label{lem:commutable-for-nice-families}
\begin{enumerate}[(1)]
    \item The class $\mathscr S$ of solvable groups is $(3,1/10)$-commutable.

    \item For every $r\geq 1$, the class $\mathscr G_r$ of groups of composition rank at most $r$ is $(4,\beta_r)$-commutable for some
    \[\beta_r\geq e^{-O(\log^2r\log\log r)}.\]
\end{enumerate}
\end{lemma}

\begin{proof} Recall from \cref{lem:abelian-commutable} that the class $\mathscr A$ of abelian groups is $(2,1/3)$-commutable. The derived series of any solvable group $G$ is a normal series all of whose successive quotients are abelian. So, part (1) follows from \cref{lem:abelian-commutable} and \cref{lem:series-commutable}.

For (2), we let $\mathscr L_r$ be the class of groups which can be found as a subgroup of $\GL_r(\FF)$ for some field $\FF$ and $\mathscr L_r'$ be the class of groups which can be found as a subgroup of some direct product $\GL_r(\FF)\times\cdots\times\GL_r(\FF)$ for some field $\FF$. In the notation of \cref{thm:GL}, write
\[\eps_r:=\frac1{M(r)}=e^{-O(\log^2r\log\log r)}\quad\text{and}\quad C_r:=Cr^9\]
so that, for every $K$-approximate subgroup of $\GL_r(\FF)$ for any field $\FF$, there exists some abelian subgroup $H\subset\GL_d(\FF)$ for which
\[\abs{A^2\cap H}\geq\frac1{(2K)^{C_r}}\abs{A^2}^{\eps_r}.\]
By \cref{lem:subgroup-commutable} and \cref{lem:abelian-commutable}, i.e.\ the fact that abelian groups are $(2,1/3)$-commutable, we conclude that $\mathscr L_r$ is $(2,\beta_r)$-commutable for
\[\beta_r:=\frac{(1/3)\eps_r}{3C_r+12}=\frac{\eps_r}{9C_r+36}.\]

Next, we consider the class $\mathscr L_r'$. Take any $G\in\mathscr L_r'$, and write $G\subset H_1\times\cdots\times H_m$ for some positive integer $m$, where $H_1\cong\cdots\cong H_m\cong\GL_r(\FF)$ for some field $\FF$. For each $0\leq i\leq m$, define
\[G_i=G\cap(\set{1}\times\cdots\times\set{1}\times H_{i+1}\times\cdots\times H_m).\]
Each $G_i$ is a normal subgroup of $G$, and so $G=G_0\unrhd G_1\unrhd\cdots\unrhd G_m=\set{1}$ is a normal series of $G$ of finite length $m$. Moreover, each successive quotient $G_{i-1}/G_i$ in this series is embedded naturally into $H_i\cong\GL_r(\FF)$ by the projection map onto the $H_i$ factor. We conclude that any $G\in\mathscr L_r'$ has a normal series of finite length wherein all successive quotients lie in $\mathscr L_r$. Therefore, \cref{lem:series-commutable} implies that $\mathscr L_r'$ is $(3,\beta_r)$-commutable (since $\beta_r<1/10$). 

Finally, by \cref{lem:chief-series}, every group in $\mathscr G_r$ possesses a normal series wherein every successive quotient is contained in a product $\GL_r(\FF)\times\cdots\times\GL_r(\FF)$ for some field $\FF$ --- in particular, wherein every successive quotient lies in $\mathscr L_r'$. We conclude by one final application of \cref{lem:series-commutable} that $\mathscr G_r$ is $(4,\beta_r)$-commutable. The asymptotics on $\eps_r$ and $C_r$ given to us by \cref{thm:GL} imply that $\beta_r=\exp(-O(\log^2r\log\log r))$ as desired.
\end{proof}

Having built up the relevant concepts and definitions, \cref{thm:solv-abelian-struct} follows from \cref{lem:commutable-to-abelian-substruct} and \cref{lem:commutable-for-nice-families}(1).

\begin{proof}[Proof of \cref{thm:solv-abelian-struct} assuming \crefrange{lem:abelian-commutable}{lem:commutable-to-abelian-substruct}] Let $G$ be a solvable group with $A\subset G$ a symmetric subset satisfying $\abs{A^3}\leq K\abs{A}$. By \cref{lem:commutable-for-nice-families}(1), $G$ lies in a $(3,1/10)$-commutable family. As a result, \cref{lem:commutable-to-abelian-substruct} implies that there exists some abelian subgroup $H\leq G$ which satisfies
\[\abs{A^4\cap H}\geq\exp\pr*{\Omega\pr*{\frac{\log^{1/6}\abs{A}}{\log 2K}}}.\qedhere\]
\end{proof}

Correspondingly, \cref{thm:gen-abelian-struct} follows from \cref{lem:commutable-to-abelian-substruct} and \cref{lem:commutable-for-nice-families}(2).

\begin{proof}[Proof of \cref{thm:gen-abelian-struct} assuming \crefrange{lem:abelian-commutable}{lem:commutable-to-abelian-substruct} and \cref{thm:GL}] Let $G$ be a finite group, and suppose that $G$ has no subquotient isomorphic to the alternating group $\A_d$. Let $A\subset G$ be a finite symmetric subset satisfying $\abs{A^3}\leq K\abs{A}$. By \cref{lem:comp-rank-alternating}, $G$ has composition rank at most $r=O(d^2)$. By \cref{lem:commutable-for-nice-families}(2) (in the proof of which we have used \cref{thm:GL}), the group $G$ lies in a $(4,\beta)$-commutable family, where $\beta=\beta_r=\exp(-O(\log^2d\log\log d))$ is as in \cref{lem:commutable-for-nice-families}(2). As a result, \cref{lem:commutable-to-abelian-substruct} implies that there exists some abelian subgroup $H\leq G$ which satisfies
\[\abs{A^4\cap H}\geq\exp\pr*{\Omega\pr*{\beta^{1/8}\frac{\log^{1/8}\abs{A}}{\log 2K}}}\geq\exp\pr*{\Omega\pr*{\frac{\log^{1/8}\abs{A}}{\exp(O(\log^2d\log\log d))\cdot\log 2K}}}.\qedhere\]
\end{proof}

\section{Proofs of \texorpdfstring{$(R,\beta)$}{(k,R,beta)}-commutability properties}\label{sec:commutable-proofs}

In this section, we prove \crefrange{lem:abelian-commutable}{lem:commutable-to-abelian-substruct}.

We begin by proving \cref{lem:abelian-commutable}, the statement that abelian groups are $(2,1/3)$-commutable. The proof is very similar to that of \cref{lem:solvable-gpvsn-step}, in that we take a maximal set of ``new'' commuting elements and use conjugation by said elements to derive our bounds.

\begin{proof}[Proof of \cref{lem:abelian-commutable}] Let $\Gamma$ be a group, let $G\leq\Gamma$, and let $N\unlhd G$ be such that $G/N$ is abelian. Let $A$ be a $K$-approximate subgroup of $\Gamma$ with $n:=\abs{A^2\cap G}$, and let $S\supset A^4\cap N$. We will show that $A^4\cap G$ contains a commuting dissociated subset $D$ satisfying $\dspan{D}\cap S=\{1\}$ and
\[\abs{D}\geq\sqrt{\frac{\log_3n}3}-\log_3(3K\abs{S}).\]

Let $D=\set{h_1,\ldots,h_m}$ be a maximal subset of $A^4\cap G$ which is commuting, dissociated, and satisfies $\dspan{D}\cap S=\{1\}$. Define a map $\phi$ on $A^2\cap G$ by
\[\phi(a)=(ah_1a^{-1}h_1^{-1},\ldots,ah_ma^{-1}h_m^{-1}).\]
Since $G/N$ is abelian, $[G,G]\subset N$. As a result, the image of $\phi$ is contained in $(A^{12}\cap N)\times\cdots\times(A^{12}\cap N)$. On the other hand, if $\phi(a)=\phi(b)$, then $b^{-1}a$ commutes with each $h_i$. As a result, for any fiber $Y$ of $\phi$, every element of the set $Y^{-1}Y\subset A^4\cap G$ commutes with every element of $D$. The maximality of $D$ thus implies that, for each $x\in Y^{-1}Y$, either (a) $x\in D$, (b) $D\cup\set{x}$ fails to be dissociated, or (c) $\dspan{D\cup\set{x}}\cap S$ contains some non-identity element of $G$. Any of these cases imply that there exists some $y\in\dspan{D}$ and some $s\in S$ with $xy=s$. This means that
\[\abs{Y}\leq \abs{Y^{-1}Y}\leq \abs{S\dspan{D}^{-1}}\leq\abs{S}\cdot\abs{\dspan{D}}=3^m\abs{S}.\]
We conclude that
\[n=\abs{A^2\cap G}=\sum_{t\in\im\phi}\abs{\phi^{-1}(t)}\leq 3^m\abs{S}\cdot\abs{\im\phi}\leq (3\cdot\abs{A^{12}\cap N})^m\abs{S}.\]
If $m=0$, then the above implies that $\abs{S}\geq n$ and the result holds trivially. Otherwise, since $A$ is a $K$-approximate group, \cref{lem:intersect-subgroup} furnishes the inequality $\abs{A^{12}\cap N}\leq K^{11}\abs{A^2\cap N}\leq K^{11}\abs{S}$. This means that $n\leq (3K^{11}\abs{S})^m\abs{S}\leq (3K\abs{S})^{11m}$; equivalently,
\[\abs{D}=m\geq\frac{\log_3n}{11(\log_3(3K\abs{S}))}.\]
Applying the inequality $\frac x{11y}\geq \sqrt{\frac x3}-y$ transforms this lower bound on $m$ into the form required.
\end{proof}

Next, we prove \cref{lem:series-commutable}. Let $\mathscr H$ be a $(R,\beta)$-commutable class of groups. Let $\Gamma$ be a group with a subgroup $G$ and let $G=G_0\unrhd G_1\unrhd\cdots\unrhd G_s=N$ be a normal series wherein $G_i/G_{i+1}\in\mathscr H$ for each $0\leq i<s$, arising from lifting such a normal series of $G/N\in\mathscr G$ to $G$. Let $A\subset\Gamma$ be a $K$-approximate subgroup and let $S$ be a centered subset of $\Gamma$ satisfying $A^4\cap N\subset S\subset C_\Gamma(G)$. Recall that we must find some commuting dissociated subset $D$ of $A^4\cap G$ satisfying $\dspan{D}\cap S=\{1\}$ and
\begin{equation}\label{eq:D-total-bound}
\abs{D}\geq\pr*{\min\pr*{\beta,\frac1{10}}\log_3\abs{A^2\cap G}}^{1/(R+1)}-\log_3(3K\abs{S}).
\end{equation}

The proof is quite similar to that of \cref{prop:solvable-gpvsn} in the genuine group setting. We will inductively construct our set $D$ by, at each step, adding elements to $D$ from as far into the series $G=G_0\unrhd G_1\unrhd\cdots\unrhd G_s=N$ as possible. An application of the $(R,\beta)$-commutability of $\mathscr H$, as well as consideration of conjugation by our new elements, will give the necessary bounds at each step.

\begin{proof}[Proof of \cref{lem:series-commutable}]
We construct such a set $D$ iteratively. Initialize $D_0=\emptyset$. Whenever we define $D_j$, we will also define $S_j:=\dspan{D_j}S$ and $H_j:=\bigcap_{x\in D_j}C_G(x)\leq G$. We maintain the property that, for each $j$, the set $D_j$ is commuting, dissociated, and satisfies $\dspan{D_j}\cap S=\{1\}$. Since $S\subset C_\Gamma(G)$ and $S_j$ is generated by $D_j\cup S$, we have for each $j$ that $S_j$ is centered and contained in $C_\Gamma(H_j)$.

Suppose we have constructed some $D_j$ (and thus $H_j$ and $S_j$). If $A^4\cap H_j\subset S_j$, we terminate the procedure; let $\ell$ be the index at which this occurs. (We shall show later in the argument that our procedure eventually terminates.) Otherwise, since
\begin{align*}
A^4\cap H_j\cap G_0&=A^4\cap H_j\not\subset S_j,\\
A^4\cap H_j\cap G_s&\subset A^4\cap N\subset S\subset S_j,
\end{align*}
there exists some index $i\in[0,s)$ for which $A^k\cap H_j\cap G_i\not\subset S_j$ and $A^k\cap H_j\cap G_{i+1}\subset S_j$. Since $G_i/G_{i+1}\in\mathscr H$ and $\mathscr H$ is subgroup-closed, we have $(H_j\cap G_i)/(H_j\cap G_{i+1})\in\mathscr H$. Let $D'$ be a maximum-size commuting dissociated subset of $A^4\cap H_j\cap G_i$ subject to the condition that $\dspan{D'}\cap S_j=\{1\}$. 

On one hand, we have $\abs{D'}\geq 1$; since $S_j$ is centered, we can take $D'$ to be a singleton consisting of any element of the nonempty set $(A^4\cap H_j\cap G_i)\setminus S_j$. On the other hand, we can apply the fact that $\mathscr H$ is $(R,\beta)$-commutable to the subquotient $(H_j\cap G_i)/(H_j\cap G_{i+1})$ of $\Gamma$, with the approximate group $A$ and the centered set $S_j\subset C_\Gamma(H_j)\subset C_\Gamma(H_j\cap G_i)$ containing $A^4\cap H_j\cap G_{i+1}$ to find some such set $D'$ satisfying
\begin{equation}\label{eq:D'-size}
\abs{D'}\geq(\beta\log_3\abs{A^2\cap H_j\cap G_i})^{1/R}-\log_3(3K\abs{S_j}).
\end{equation}
Set $D_{j+1}:=D_j\sqcup D'$. Since $D'\subset A^4\cap H_j$, every element of $D'$ commutes with every element of $D_j$. Since $D_j$ and $D'$ are commuting sets, this implies that $D_{j+1}$ is as well. Now, since $D'$ is commuting and dissociated and $\dspan{D'}\cap\dspan{D_j}S=\{1\}$, we have that $D_{j+1}=D_j\sqcup D'$ is dissociated and $\dspan{D_{j+1}}\cap S=\{1\}$. The fact that $\abs{D'}\geq 1$ ensures that $\abs{D_j}$ is a strictly increasing function of $j$; since each $D_j$ lies inside the finite set $A^4\cap G$, this in turn ensures that our process eventually terminates.

We next understand how $\abs{A^4\cap H_j}$ changes with $j$. First, if $\ell=0$ then $A^4\cap G\subset S$, and in particular
\[(\beta\log_3\abs{A^2\cap G})^{1/R}\leq \beta\log_3\abs{A^2\cap G}\leq\log_3\abs{A^4\cap G}\leq \log_3(3K\abs{S}).\]
This means that \eqref{eq:D-total-bound} holds trivially, and there is nothing to prove. We henceforth assume that $\ell>0$, which implies that $\abs{D_\ell}>0$. On one hand, by our termination condition we have
\begin{equation}\label{eq:H-bound-termination}
\frac{\abs{A^4\cap H_0}}{\abs{A^4\cap H_\ell}}\geq\frac{\abs{A^4\cap G}}{\abs{S_\ell}}\geq\frac{\abs{A^2\cap G}}{3^{\abs{D_\ell}}\abs{S_0}}.
\end{equation}
On the other hand, we claim that, if $m=\abs{D'}=\abs{D_{j+1}\setminus D_j}$, then
\begin{equation}\label{eq:H-bound-approx}
\frac{\abs{A^4\cap H_j}}{\abs{A^4\cap H_{j+1}}}\leq\pr*{K^{10}\abs{A^2\cap H_j\cap G_i}}^m.
\end{equation}
Indeed, define a map $\sigma$ on $A^2\cap H_j$ by
\[\sigma(a)=(ah_1a^{-1},\ldots,ah_ma^{-1}),\]
where $h_1,\ldots,h_m$ is an enumeration of the elements of $D'=D_{j+1}\setminus D_j$. Recall that $D'\subset A^4\cap H_j\cap G_i$. So, since $G_i\unlhd G$, we have $axa^{-1}\in A^8\cap H_j\cap G_i$ for every $a\in A^2\cap H_j$ and $x\in D'$. Thus the image of $\sigma$ is contained in $(A^8\cap H_j\cap G_i)\times\cdots\times(A^8\cap H_j\cap G_i)$. Now, let $Y$ be a maximum-size fiber of $\sigma$. By the definition of $\sigma$, every element of $Y^{-1}Y$ commutes with each element of $D'$ and thus lies in $(A^2\cap H_j)^2\cap H_{j+1}$. We conclude that
\[\abs{A^4\cap H_{j+1}}\geq\abs{Y^{-1}Y}\geq\abs{Y}\geq\frac{\abs{A^2\cap H_j}}{\abs{A^8\cap H_j\cap G_i}^m}.\]
Applying \cref{lem:intersect-subgroup} twice, we get
\[\frac{\abs{A^4\cap H_j}}{\abs{A^4\cap H_{j+1}}}\leq\frac{K^3\abs{A^2\cap H_j}}{\abs{A^4\cap H_{j+1}}}\leq K^3\abs{A^8\cap H_j\cap G_i}^m\leq (K^{10}\abs{A^2\cap H_j\cap G_i})^m,\]
showing \eqref{eq:H-bound-approx}. Now, \eqref{eq:D'-size} implies that
\begin{align*}
\log_3\abs{A^2\cap H_j\cap G_i}
&\leq \frac1\beta\pr*{\abs{D'}+\log_3(3K\abs{S_j})}^R\\
&=\frac1\beta\pr*{\abs{D_{j+1}}+\log_3(3K\abs{S_0})}^R.
\end{align*}
Combining this with \eqref{eq:H-bound-termination} and \eqref{eq:H-bound-approx} gives, letting $\kappa:=\log_3K$ for notational simplicity,
\begin{align*}
\log_3\abs{A^2\cap G}-\abs{D_\ell}-\log_3\abs{S_0}
&\leq\log_3\abs{A^4\cap H_0}-\log_3\abs{A^4\cap H_\ell}\\
&\leq\sum_{j=0}^{\ell-1}\abs{D_{j+1}\setminus D_j}\pr*{10\kappa+\log_3\abs{A^2\cap H_j\cap G_i}}\\
&\leq10\kappa\abs{D_\ell}+\frac1\beta\sum_{j=0}^{\ell-1}\abs{D_{j+1}\setminus D_j}\pr[\big]{\abs{D_{j+1}}+1+\kappa+\log_3\abs{S_0}}^R\\
&\leq10\kappa\abs{D_\ell}+\frac1\beta\abs{D_\ell}\pr[\big]{\abs{D_\ell}+1+\kappa+\log_3\abs{S_0}}^R.
\end{align*}
Let $\beta':=\min(\beta,1/10)$. Since $\abs{D_\ell},R\geq 1$, we have
\begin{align*}
\log_3\abs{A^2\cap G}
&\leq\abs{D_\ell}(1+10\kappa+\log_3\abs{S_0})+\frac1\beta\abs{D_\ell}\pr[\big]{\abs{D_\ell}+1+\kappa+\log_3\abs{S_0}}^R\\
&\leq \frac1{\beta'}\pr*{\abs{D_\ell}(1+\kappa+\log_3\abs{S_0})+\abs{D_\ell}\pr[\big]{\abs{D_\ell}+1+\kappa+\log_3\abs{S_0}}^R}\\
&\leq \frac1{\beta'}\pr[\big]{\abs{D_\ell}+\log_3(3K\abs{S_0})}^{R+1}.
\end{align*}
Rearranging gives \eqref{eq:D-total-bound}, as desired.
\end{proof}

Next, we prove \cref{lem:subgroup-commutable}, which allows us to pass from a $(R,\beta)$-commutable class of groups $\mathscr H$ to a larger class $\mathscr G$ in which every approximate group has large intersection with some group in $\mathscr H$. The lemma essentially follows from \cref{lem:intersect-subgroup-projection}. 

\begin{proof}[Proof of \cref{lem:subgroup-commutable}] Let $\mathscr H$ be a $(R,\beta)$-commutable class of groups, and let $G'$ be a group for which, if $\tilde A\subset G'$ is a $K$-approximate group, then there exists a subgroup $H\leq G$ satisfying $H\in\mathscr H$ and
\begin{equation}\label{eq:tilde-A-subgroup}
\abs[\big]{\tilde A^2\cap H}\geq\frac1{(2K)^C}\abs[\big]{\tilde A^2}^\eps.
\end{equation}
We must show that, for any group $\Gamma$ with subquotient $G/N\cong G'$, any $K$-approximate subgroup $A$ of $\Gamma$, and any centered subset $S\subset C_\Gamma(G)$ with $S\supset A^4\cap N$, the set $A^4\cap G$ contains a commuting dissociated subset $D$ satisfying $\dspan{D}\cap S=\{1\}$ and
\begin{equation}\label{eq:D-total-subgp}
\abs{D}\geq(\beta'\log_3\abs{A^2\cap G})^{1/R}-\log_3(3K\abs{S})
\end{equation}
where $\beta'=\frac{\beta\eps}{3C+12}$. First, let $\pi\colon G\to G/N\cong G'$ be the natural quotient map. The set $\tilde A:=\pi(A^2\cap G)$ is a $K^3$-approximate subgroup of $G'$. So, we can find some subgroup $H\leq G'$ with $H\in\mathscr H$ which satisfies \eqref{eq:tilde-A-subgroup} with $K$ replaced by $K^3$. By applying \cref{lem:intersect-subgroup-projection} to the $K^3$-approximate subgroup $A^2\cap G$ and the projection $\pi\colon G\to G'$, we have
\begin{align*}
\frac{K^{11}\abs{A^2\cap\pi^{-1}(H)}}{\abs{(A^2\cap G)^2}}&\geq\frac{\abs{(A^2\cap G)^6\cap\pi^{-1}(H)}}{\abs{(A^2\cap G)^2}}\\
&\geq\frac{\abs{\pi((A^2\cap G)^2)\cap H}}{\abs{\pi((A^2\cap G)^2)}}\\
&=\frac{\abs{\tilde A^2\cap H}}{\abs{\tilde A^2}}\geq\frac1{(2K^3)^C}\abs{\tilde A^2}^{-(1-\eps)}\geq \frac{\abs{A^2\cap G}^\eps}{(2K^3)^C\abs{(A^2\cap G)^2}}
\end{align*}
(where we have used \cref{lem:intersect-subgroup} in the first inequality). We conclude that
\begin{equation}\label{eq:pi-H-size}
\abs{A^2\cap\pi^{-1}(H)}\geq\frac1{(2K)^{3C+11}}\abs{A^2\cap G}^\eps.
\end{equation}
Now, $\pi^{-1}(H)/N\cong H\in\mathscr H$, and so we may apply the $(R,\beta)$-commutability of $\mathscr H$ to the subquotient $\pi^{-1}(H)/N$ of $\Gamma$, with the approximate group $A$ and the centered set $S\subset C_{\Gamma}(G)\subset C_\Gamma(\pi^{-1}(H))$ containing $A^4\cap N$ to find some commuting dissociated set $D\subset A^4\cap \pi^{-1}(H)\subset A^4\cap G$ which satisfies $\dspan{D}\cap S=\{1\}$ as well as
\begin{equation}\label{eq:D-subgp-1}
\abs{D}\geq(\beta\log_3\abs{A^2\cap\pi^{-1}(H)})^{1/R}-\log_3(3K\abs{S}).
\end{equation}
If $\abs{A^2\cap G}^\eps<(2K)^{3C+12}$, then (since $R>1$)
\[\pr*{\frac{\beta\eps}{3C+12}\log_3\abs{A^2\cap G}}^{1/R}\leq \frac{\beta\eps}{3C+12}\log_3\abs{A^2\cap G}<\log_3(2K)<\log_3(3K\abs{S}),\]
and \eqref{eq:D-total-subgp} holds trivially. Otherwise, \eqref{eq:pi-H-size} implies that $\abs{A^2\cap\pi^{-1}(H)}\geq\abs{A^2\cap G}^{\eps/(3C+12)}$, and so \eqref{eq:D-total-subgp} follows from \eqref{eq:D-subgp-1}.
\end{proof}

Finally, we prove \cref{lem:commutable-to-abelian-substruct}. This lemma enables us to upgrade $(R,\beta)$-commutability to the existence of a (quasi-polynomially) large abelian approximate group inside $A^4$ for any approximate group $A$. Instead of applying $(R,\beta)$-commutability to the set $A$ which we are given, we first apply a result of Sanders (\cref{prop:sanders-quasipoly}) to find a set $S_0$, not much smaller than $A$, a large iterated product set of which lies in $A^4$. We then apply $(R,\beta)$-commutability to $S_0^2$ to find a large commuting, dissociated subset $D$ of $S_0^8$. Since a large iterated product set of $S_0$ lies in $A^4$, we can use the dissociativity of $D$ to find a commuting subset of $A^4$ which has size exponential in $\abs{D}$.

\begin{proof}[Proof of \cref{lem:commutable-to-abelian-substruct}] Let $\mathscr G$ be a $(R,\beta)$-commutable family, and let $G\in\mathscr G$. Let $A\subset G$ be a finite symmetric subset satisfying $\abs{A^3}\leq K\abs{A}$. We must find some abelian subgroup $H$ of $G$ satisfying
\begin{equation}\label{eq:commutable-to-substruct-goal}
\abs{A^4\cap H}\geq\exp\pr*{c\frac{(\beta\log\abs{A})^{1/(2R)}}{\log 2K}}
\end{equation}
for some absolute constant $c>0$.

By Ruzsa's triangle inequality, we have $\abs{A^4}\leq K^2\abs{A}$. Let $t\geq 4$ be a parameter which we will choose later. We begin by applying \cref{prop:sanders-quasipoly} to $A$. For some absolute constant $C\geq1$, we can find some centered set $S_0$ for which $S_0^{2t}\subset A^4$ and
\[\abs{S_0}\geq\exp\pr*{-4Ct^2(\log 2K)^2}\abs{A}.\]
Let $L:=K^2\exp(4Ct^2(\log 2K)^2)$, so that $L\abs{S_0}\geq\abs{A^4}\geq\abs{S_0^{2t}}$. Let $B:=S_0^2$. By Ruzsa's covering lemma applied to $S_0^2$ and $S_0$, there is some subset $Y\subset G$ with $\abs{Y}\leq\abs{S_0^5}/\abs{S_0}\leq L$ and $B^2\subset YB$. The set $B$ is thus an $L$-approximate group satisfying $B^t\subset A^4$.

So, we may apply the $(R,\beta)$-commutability of $\mathscr G$ to the $L$-approximate group $B\subset G$, with $\Gamma=G$ and $N=\set{1}$ and $S=\set{1}$. We conclude that $B^4$ contains a commuting dissociated set $D$ satisfying
\begin{equation}\label{eq:D-bound-in-substruct}
\abs{D}\geq(\beta\log_3\abs{B})^{1/R}-\log_3(3L).
\end{equation}
Let $\ell:=\min(\abs{D},\floor{t/4})$, and let $D'$ be an arbitrary subset of $D$ of size $\ell$, so that $(\set{1}\cup D')^\ell\subset B^t\subset A^4$. Let $H:=\langle D'\rangle$ be the abelian subgroup of $G$ generated by $D'$. The set $(\set{1}\cup D')^\ell$ contains the product of every subset of $D'$; since $D'$ is dissociated, this set thus has size at least $2^\ell$. We conclude that $\abs{A^4\cap H}\geq 2^\ell$.

What remains is to choose $t$ to maximize $\ell$. We select
\[t=\floor*{\frac{(\beta\log_3\abs{A})^{1/(2R)}}{5C\log(2K)}}.\]
We can absorb the case where $t$ as chosen above is less than $4$ into \eqref{eq:commutable-to-substruct-goal} by choosing $c$ sufficiently small. We may thus assume that $t\geq 4$. We claim that $\abs{D}\geq t/4$. This will finish the proof: if $\abs{D}\geq t/4$, then the inequality $\abs{A^4\cap H}\geq 2^\ell\geq\exp(t/8)$ is of the form \eqref{eq:commutable-to-substruct-goal}. Indeed, we have
\begin{align*}
\abs{D}
&\geq(\beta\log_3\abs{B})^{1/R}-\log_3(3L)\\
&\geq(\beta\log_3(\abs{A}/L))^{1/R}-\log_3(3L)\\
&\geq(\beta\log_3\abs{A})^{1/R}-2\log_3(3L)\\
&\geq25C^2t^2(\log 2K)^2-2\log(3K^2)-8Ct^2(\log 2K)^2>15C^2t^2(\log 2K)^2>t.\qedhere
\end{align*}
\end{proof}

\section{Non-solvable groups: strategy and outline}\label{sec:GL-setup}

The next four sections are devoted to the proof of \cref{thm:GL}; we will prove \cref{prop:reg-el} along the way. Throughout these sections, $\FF$ will always be an algebraically closed field (of arbitrary characteristic). This property of $\FF$ can be imposed on the statement of \cref{thm:GL} with no loss of generality. 

In this section, we outline the strategy we will use to prove \cref{thm:GL}, and reduce the proof of \cref{thm:GL} to finding a Borel subgroup\footnote{In $\GL_d(\FF)$, the \emph{Borel subgroups} are those subgroups consisting, in some basis, of the upper triangular matrices.} of $\GL_d(\FF)$ with large intersection with $A^2$.

All of the terminology and definitions we use are presented in \cref{sec:prelim-matrix}.

\subsection{Strategy}\label{sec:GL-strategy}
We begin with $A\subset\GL_d(\FF)$ and wish to find some abelian group $H\leq\GL_d(\FF)$ for which the intersection $A^2\cap H$ is large. We first find such a group $H$ which is contained in a Borel subgroup of $\GL_d(\FF)$, and then use a more tailored version of \cref{thm:solv-abelian-struct} (see \cref{sec:borel-to-abelian}) to conclude \cref{thm:GL}. For the remainder of this sketch, our goal is to find some Borel subgroup $\cB$ for which $\abs{A^2\cap\cB}$ is large.

We will do this by finding successively smaller block-triangular (in some basis) subgroups which have a large intersection with $A^2$. Once we have reduced the sizes of all our blocks each to one, many elements of $A^2$ will be contained in the group of upper-triangular matrices in some basis, and we will be done. 

The starting point of our argument is the following simple inequality, of which we have used some variations implicitly in the previous two sections.

\begin{lemma}\label{lem:pivot} For any subset $A$ of any group $G$, and any element $a\in A$, one has
\[\abs{A^2\cap C(a)}\cdot\abs{A^3\cap\Conj(a)}\geq\abs{A}.\]
\end{lemma}
\begin{proof} Consider the map $\sigma\colon A\to A^3$ defined by $\sigma(g)=gag^{-1}$. On one hand, the image of $\sigma$ is contained within $A^3\cap\Conj(a)$. On the other hand, for each fiber $S$ of $\sigma$, the set $S^{-1}S$ is contained within $A^2\cap C(a)$. We conclude that
\[\abs{A}=\sum_{b\in\im\sigma}\abs{\sigma^{-1}(b)}\leq\abs{\im\sigma}\max_{b\in\im\sigma}\abs{\sigma^{-1}(b)}\leq\abs{A^3\cap\Conj(a)}\cdot\abs{A^2\cap C(a)}.\qedhere\]
\end{proof}

So, we begin by selecting some $a\in A$. By \cref{lem:pivot}, either $\abs{A^2\cap C(a)}$ or $\abs{A^3\cap\Conj(a)}$ is large. We examine the first case first. The set $C(a)$ is exceedingly nice: it is a subgroup of $\GL_d(\FF)$, and it is also exactly the set of invertible matrices in some subspace of $\Mat_d(\FF)$. If $a\in\GL_d(\FF)$ is a regular semisimple element (i.e., it possesses $d$ distinct eigenvalues), then $C(a)$ is a diagonal subgroup, and so we are immediately done if $A^2\cap C(a)$ is large. In practice, we will not necessarily be able to pick $a$ to be regular semisimple. However, if $a$ is ``not too irregular,'' i.e.\ $a$ does not possess an eigenspace of too large dimension, then $C(a)$ will be contained in some block-triangular subgroup which is suitably smaller than $\GL_d(\FF)$ itself (see \cref{lem:centralizer}), and we will gain enough by passing down to $A^2\cap C(a)$. We will use \cref{prop:reg-el} to find such an element $a$, or else to obtain a reduction to a much simpler group directly.

We now examine the case where $\abs{A^3\cap\Conj(a)}$ is large. The set $\Conj(a)$ is not as nice as $C(a)$. However, it is contained within the codimension-one affine subspace $\bX:=\set{x\in\Mat_d(\FF):\tr x = \tr a}$. This set is very unlike a subgroup --- it has poor multiplicative structure. We will be able to utilize this. As a starting point, if $A^3\cap \bX$ is large, then a double-counting argument gives the existence of many $a\in A$ for which $A^4\cap(\bX\cap a\bX)$ is also large. Unless $a$ is a scalar multiple of the identity, $\bX\cap a\bX$ is a codimension-\emph{two} affine subspace of $\Mat_d(\FF)$, and so we can replace $\bX$ by the smaller $\bX\cap a\bX$ without losing too much. If, on the other hand, $a$ is a scalar multiple of the identity for many $a\in A$, then we immediately obtain a reduction to a (much) smaller subgroup. Using a more involved version of this approach (the basis of which is \cref{lem:dimension-decrement}), we will thus be able to upgrade the statement ``many elements of $A^3$ have the same trace'' to a statement (roughly) of the form ``$A^2$ has large intersection with some block-triangular subgroup of $\GL_d(\FF)$ of dimension at most $\frac1{\sqrt2}d^2$.''

The main complication in carrying through this sketch is that, at step $t$ of our iterative procedure, we will not have a subset of $\GL_{d(t)}(\FF)$ for $d(t)$ a decreasing function of $t$, but rather a subset of some block-triangular subgroup of $\GL_d(\FF)$ (with $d$ fixed) whose blocks decrease in size as $t$ increases. Our notation for carrying this data through the proof is based on \emph{flags} in vector spaces; see \cref{sec:prelim-matrix}.

\subsection{Connections to prior work}\label{sec:GL-history}

The proof of \cref{thm:GL} is inspired by the methods used by various authors to prove Babai's conjecture on diameter bounds for Cayley graphs of simple groups in the bounded rank case \cite{HelfgottSL2,HelfgottSL3,BGTLinear,PyberSzabo}. Two central tools in these results are \emph{escape from subvarieties} and the \emph{Larsen--Pink inequality}. 

These results are most naturally stated in $\SL_d(\FF)$. Both rely an a notation of \emph{complexity} for varieties whose definition we omit, except to say that varieties cut out by linear polynomials have low complexity. Let $A\subset\SL_d(\FF)$ be a $K$-approximate group. The Larsen--Pink inequality for centered sets (see \cite[Theorem~4.1]{BGTLinear}) states roughly that either
\begin{enumerate}[(a)]
    \item for every small-complexity subvariety $\bV$ of $\SL_d(\FF)$, we have
    \begin{equation}\label{eq:small-inter-w-variety}
    \abs{A^2\cap\bV}\leq K^{O_{d,\bV}(1)}\abs{A}^{\dim\bV/\dim\SL_d(\FF)},
    \end{equation}
    or else
    
    \item $A$ lives entirely in some small-complexity subvariety of $\SL_d(\FF)$.
\end{enumerate}
Escape from subvarieties (see \cite[Lemma~4.4]{HelfgottSL2}, \cite[Lemma~3.11]{BGTLinear}, or \cite[Lemma~51]{PyberSzabo}) allows one to upgrade condition (b) from a subvariety to a subgroup: if $A^m$ lives in some small-complexity subvariety $\bV$ (for some $m$ sufficiently large in terms of the complexity of $\bV$), then in fact $A$ is contained in some small-complexity subgroup of $\SL_d(\FF)$.

We are unable to use these results as they appear in the literature for two reasons: first, the Larsen--Pink inequality does not directly apply when the ambient group $\SL_d(\FF)$ is replaced by a non-simple group; second, the quantitative dependence on $d$ is generally weaker than we are willing to accept. However, the quantitative dependence on $\abs{A}$ of \eqref{eq:small-inter-w-variety} is much \emph{stronger} than we need: we are perfectly happy to accept a relatively modest upper bound on $\abs{A^2\cap\bV}$ for $\bV$ of interest, while \eqref{eq:small-inter-w-variety} is tight up to the dependence on $K$, $d$, and the complexity of $\bV$.

Let us consider the case of applying \cref{lem:pivot} as discussed in the previous subsection, with the aim of finding a large block-triangular subgroup containing many elements of $A^2$. We obtain for some $a$ that
\[\abs{A^2\cap C(a)}\cdot\abs{A^3\cap\Conj(a)}\geq\abs{A}.\]
If \eqref{eq:small-inter-w-variety} held, then both $\abs{A^2\cap C(a)}$ and $\abs{A^3\cap\Conj(a)}$ would be bounded within a multiplicative factor of $K^{O_d(1)}$. However, even a weak bound like $\abs{A^2\cap C(a)}\geq\abs{A^2}^{d^{-10}}$ is enough for our purposes. This allows us to assume that $\abs{A^3\cap\Conj(a)}$ is exceedingly large. In particular, letting $\bX$ be the (codimension one) subvariety of $\Mat_d(\FF)$ cut out by the equation $\tr x=\tr a$, we obtain a contradiction to \eqref{eq:small-inter-w-variety} for $\bV=\bX$. In fact, we obtain a statement qualitatively close to condition (b) in the Larsen--Pink inequality: such a large proportion of $A^3$ is contained in $\bX$, relative to its complexity, that ``escape from subvarieties'' machinery should give us a subgroup containing a large portion of $A^3$ in this case. In other words, by making quantitative what precisely we require from our iteration, we eliminate the need for a general Larsen--Pink-inspired argument entirely. What we do need is a quantitative version of ``escape from affine subspaces.'' This is given by \cref{lem:subspace-to-subgroup-flag}, stated in the following section, and its proof is motivated by the proofs of aforementioned escape results in the literature.

There is one other application of escape from subvarieties in proofs of diameter bounds which will also be necessary in our work. This is to find, for some $m=O_d(1)$, regular semisimple elements of $A^m$. This is indeed an ``escape from subvarieties'' problem: the set of elements which are not regular and semisimple forms a subvariety $\bZ$ of $\Mat_d(\FF)$ (cut out by the discriminant of the characteristic polynomial). However, the situation for us is poor here in two respects. First, $\bZ$ cannot be placed inside an affine subspace of $\Mat_d(\FF)$, so we cannot treat it the same way we treat $\Conj(a)$. Second, we are given no guarantee that $A$ generates a large subgroup of $\GL_d(\FF)$. As such, there is no reason to expect that $A^m$ contains regular semisimple elements for any $m$, only that we might be able to reduce to a smaller ambient group if this is not the case. This dichotomy is captured by \cref{prop:reg-el} which, although far from guaranteeing the existence of regular semisimple elements, is enough for our purposes. Our proof of \cref{prop:reg-el} is not related to proofs of general escape results but instead relies crucially on the particular subvariety at hand, and is motivated by arguments of Haramaty--Shpilka and Lampert in the additive group $\Mat_d(\FF)$ rather than arguments in the multiplicative setting. The connection and approach are explained at the beginning of \cref{sec:not-too-irreg}.

\subsection{The last step: from solvable to abelian}\label{sec:borel-to-abelian}

To conclude this section, we present an argument, similar to the proofs given in \cref{sec:commutable-proofs}, which forms the conclusion of the proof of \cref{thm:GL}. The result is a stronger version of \cref{thm:solv-abelian-struct} in the case where the solvable group is assumed to lie in a Borel subgroup $\cB$ of $\GL_d(\FF)$ for some field $\FF$ and some $d$. We allow our bounds to depend on $d$ but not on the field $\FF$.

\begin{lemma}\label{lem:borel-abelian-struct} Let $A$ be a $K$-approximate subgroup of a Borel subgroup of $\GL_d(\FF)$. There exists an abelian subgroup $H$ of $\GL_d(\FF)$ for which
\[\abs{A^2\cap H}\geq K^{-33}\abs{A^2}^{1/d^4}.\]
\end{lemma}

This argument we use to prove \cref{lem:borel-abelian-struct} mirrors that used to prove \cref{thm:solv-abelian-struct}. We are able to obtain polynomial bounds on $\abs{A^2\cap H}$ for $H$ abelian by exploiting two crucial properties of $\GL_d(\FF)$. One is the fact that, if $\cB\leq\GL_d(\FF)$ is a Borel subgroup, then the commutator subgroup $[\cB,\cB]$, consisting of the upper triangular matrices with ones on the diagonal, is nilpotent. The other is the following ``finite-dimensional'' condition. We call a group \emph{$t$-centrally bounded} if the following holds: suppose $h_1,\ldots,h_m\in\GL_d(\FF)$ are such that, for each $1\leq k\leq m$, each $h_k$ is a non-central element of $C(h_1)\cap\cdots\cap C(h_{k-1})$. Then $m\leq t$.

\begin{lemma}\label{lem:GL-central-bounded} The group $\GL_d(\FF)$ is $(d^2-1)$-centrally bounded.
\end{lemma}
\begin{proof} The property of $h_1,\ldots,h_m$ defining central boundedness requires that the $\{h_k\}$ are linearly independent (as elements of $\Mat_d(\FF)$); since they are linearly independent and commute, we conclude that $m<d^2$. Indeed, if $h_k$ is a linear combination of $\{h_1,\ldots,h_{k-1}\}$, then $h_k$ commutes with every element of $\GL_d(\FF)$ that each of $h_1,\ldots,h_{k-1}$ all commute with, and thus lies in the center of $C(h_1)\cap\cdots\cap C(h_{k-1})$.
\end{proof}

We can use this property to show the following two lemmas, the proofs of which are quite similar. The first reduces our problem ``down one level'' in the derived series, which in our case is enough to reduce our setting from solvable to nilpotent (since the commutator subgroup of a Borel subgroup of $\GL_d(\FF)$ is nilpotent). The second uses central boundedness to go from a nilpotent group to an abelian group, with a structure similar to that of the proof of \cref{lem:series-commutable}.

\begin{lemma}\label{lem:intersect-commutator} Let $G$ be a $t$-centrally bounded group, and let $A$ be a $K$-approximate subgroup of $G$. Then there exists some abelian subgroup $H\leq G$ for which
\[\abs{A^2\cap H}\cdot\pr*{K^{10}\abs{A^2\cap[G,G]}}^t\geq\abs{A^2}.\]
\end{lemma}
\begin{proof} Construct a sequence of elements $h_1,h_2,\ldots\in A^2$ iteratively such that, for each $k\geq 0$, each $h_{k+1}$ is a non-central element of $H_k:=C(h_1)\cap\cdots\cap C(h_k)$. Terminate the sequence as $h_1,\ldots,h_m$ when such a choice is no longer possible; at this point, the set $A^2\cap H_m$ is contained in the abelian group $Z(H_m)$.

We apply \cref{lem:pivot} to $A':=A^2\cap H_k$ for each $k$. We conclude, using also \cref{lem:intersect-subgroup}, that
\[\abs{A^4\cap H_{k+1}}\geq\abs{(A')^2\cap C(h_{k+1})}\geq\frac{\abs{A'}}{\abs{(A')^3\cap\Conj(h_{k+1})}}.\]
We also have that
\[\abs{A^8\cap[G,G]}\geq\abs{(A')^4\cap \Conj(h_{k+1})h_{k+1}^{-1}}\geq\abs{(A')^3\cap\Conj(h_{k+1})}.\]
Combining the above two inequalities and using \cref{lem:intersect-subgroup} gives
\[\frac{\abs{A^2\cap H_k}}{\abs{A^2\cap H_{k+1}}}\leq K^{10}\abs{A^2\cap[G,G]}.\]
Multiplying the above over all $0\leq k<m$, we conclude
\[\frac{\abs{A^2}}{\abs{A^2\cap H_m}}\leq\pr*{K^{10}\abs{A^2\cap [G,G]}}^m\leq\pr*{K^{10}\abs{A^2\cap[G,G]}}^t.\]
The result follows.
\end{proof}

\begin{lemma}\label{lem:nilpotent-central-bounded} Let $G$ be a $t$-centrally bounded nilpotent group. Let $A$ is a $K$-approximate subgroup of $G$. Then there exists some abelian subgroup $H\leq G$ satisfying
\[\abs{A^2\cap H}\geq\frac1{K^{10}}\abs{A^2}^{1/(t+1)}.\]
\end{lemma}
\begin{proof} Let
\[G=G_0\rhd G_1\rhd\cdots G_s=\{1\}\]
be the lower central series of $G$ (so $G_{i+1}=[G_i,G]$ for each $i\geq 0$; since $G$ is nilpotent, this series terminates at the trivial subgroup after some finitely many steps). Choose a sequence of elements $h_1,h_2,\ldots\in A^2$ iteratively such that, for each $k\geq 0$, each $h_{k+1}$ is a non-central element of $H_k:=C(h_1)\cap\cdots\cap C(h_k)$. While such an element exists, we impose that $h_{k+1}$ is chosen from $G_i$ for as large an $i$ as possible, so that
\begin{equation}\label{eq:center-contain}
A^2\cap H_k\cap G_{i+1}\subset A^2\cap Z(H_k).
\end{equation}
When such an element no longer exists, terminate the construction of the sequence, and let $m$ be its length. At this point, we have $A^2\cap H_m\subset A^2\cap Z(H_m)$.

We claim that, for each $1\leq k\leq m$,
\begin{equation}\label{eq:central-bounded-intermediate}
\frac{\abs{A^2\cap H_k}}{\abs{A^2\cap H_{k+1}}}\leq K^{10}\abs{A^2\cap Z(H_k)}.
\end{equation}
To prove this, let $A':=A^2\cap H_k$, and let $i$ be such that $h_{k+1}\in G_i\setminus G_{i+1}$. Apply \cref{lem:pivot} to $A'$ and $h_{k+1}$ to find
\[\abs{(A')^2\cap C(h_{k+1})}\cdot\abs{(A')^3\cap\Conj_{H_k}(h_{k+1})}\geq\abs{A'}.\]
Now, we observe that $\Conj_{H_k}(h_{k+1})h_{k+1}^{-1}\subset H_k\cap [G,G_i]=H_k\cap G_{i+1}$. This gives \eqref{eq:central-bounded-intermediate}:
\[\abs{A^2\cap H_k}\leq \abs{A^4\cap H_{k+1}}\cdot\abs{A^8\cap H_k\cap G_{i+1}}\leq K^{10}\abs{A^2\cap H_{k+1}}\cdot\abs{A^2\cap Z(H_k)},\]
where we have used \cref{lem:intersect-subgroup} and \eqref{eq:center-contain}.

Finally, multiplying \eqref{eq:central-bounded-intermediate} for each $0\leq k<m$ and using the termination condition $A^2\cap H_m\subset A^2\cap Z(H_m)$ gives
\[\abs{A^2}\leq K^{10m}\prod_{k=0}^m\abs{A^2\cap Z(H_k)}\leq K^{10(m+1)}\pr*{\max_{\substack{H\leq G\\H\text{ abelian}}}\abs{A^2\cap H}}^{m+1}.\]
The result follows from the $t$-central boundedness of $G$.
\end{proof}

\cref{lem:borel-abelian-struct} now follows from combining \crefrange{lem:GL-central-bounded}{lem:nilpotent-central-bounded}.

\begin{proof}[Proof of \cref{lem:borel-abelian-struct}] Let $\cB$ be a Borel subgroup of $\GL_d(\FF)$ and let $A\subset\cB$ be a $K$-approximate subgroup of $\cB$. By \cref{lem:GL-central-bounded}, we know that $\GL_d(\FF)$, and thus $\cB$, is $(d^2-1)$-centrally bounded. We conclude from \cref{lem:intersect-commutator} that there exists some abelian $H_1\leq \cB$ for which
\[\abs{A^2\cap H_1}\cdot\pr*{K^{10}\abs{A^2\cap[\cB,\cB]}}^{d^2-1}\geq\abs{A^2}.\]
Let $A':=A^2\cap[\cB,\cB]$; by \cref{lem:intersect-subgroup}, $A'$ is a $K^3$-approximate subgroup of $[\cB,\cB]$. We can then apply \cref{lem:nilpotent-central-bounded} to $A'$ to find some abelian subgroup $H_2\leq[\cB,\cB]$ for which
\[\abs{A^2\cap H_2}\geq \frac1{K^3}\abs{(A')^2\cap H_2}\geq\frac1{K^{33}}\abs{A'}^{1/d^2}.\]
Combining the above two inequalities gives
\[\abs{A^2}\leq\abs{A^2\cap H_1}\cdot\pr*{K^{10}\pr*{K^{33}\abs{A^2\cap H_2}}^{d^2}}^{d^2-1}\leq \pr*{K^{33}\max_{H\in\{H_1,H_2\}}\abs{A^2\cap H}}^{d^4}.\]
The result follows.
\end{proof}

\section{Preliminary details of the proof of \texorpdfstring{\cref{thm:GL}}{Theorem 1.12}}\label{sec:GL-prelim}

We begin this section by presenting some of the notation we will carry through \crefrange{sec:GL-prelim}{sec:subspace-to-subgroup} and proving some related simple lemmas. Then, we state \cref{lem:subspace-to-subgroup-flag}, to be proven in \cref{sec:subspace-to-subgroup}, an ``escape from affine subspaces'' result which will enable us to treat the case when $\abs{A^3\cap\Conj(a)}$ is large in the approach outlined in \cref{sec:GL-strategy}. Finally, we derive \cref{thm:GL} assuming \cref{prop:reg-el} and \cref{lem:subspace-to-subgroup-flag}, iterating the selection of ``not too irregular'' $a\in A^2$ described in \cref{sec:GL-strategy}.

\subsection{Matrix groups and linear algebra}\label{sec:prelim-matrix}

To prove \cref{thm:GL} and \cref{prop:reg-el} we will need some facts and notation surrounding matrix groups and algebras.\footnote{All of our algebras are (not necessarily unital) subalgebras of some matrix algebra over a field, with matrix multiplication as the product rather than a Lie bracket.} Given a vector space $V$ over a field $\FF$, we write $\End(V)$ for the endomorphism ($\FF$-)algebra of $V$ and $\GL(V)$ for the automorphism group of $V$. We write $\Mat_d(\FF):=\End(\FF^d)$ and $\GL_d(\FF):=\GL(\FF^d)$. Throughout, for an element $x$ of a matrix group, we write $\eigm(x)$ for the largest dimension of any eigenspace of $x$. 

\begin{definition}[Flags and stabilizers] Let $V$ be a finite-dimensional $\FF$-vector space. A \emph{flag} $\cF$ on $V$ is a nested sequence of subspaces $\set{0}=V_0\subsetneq V_1\subsetneq\cdots\subsetneq V_m=V$.
\begin{itemize}
    \item The \emph{trivial flag} on a vector space $V$ of positive dimension is the sequence $\set{0}\subsetneq V$.

    \item The flag $\cF$ is \emph{complete} if $\dim V_i-\dim V_{i-1}=1$ for each $1\leq i\leq m$.

    \item A flag $\mathcal G$ is a \emph{refinement} of $\cF$ if every subspace in the sequence defining $\cF$ also appears in the sequence defining $\mathcal G$.

    \item The stabilizer of $\cF$ is a block upper-triangular subalgebra of $\End(V)$, given by
    \[\set{x\in\End(V):xV_i\subset V_i\text{ for all }i},\]
    which we write as $\bP_\cF$.

    \item We write $H_\cF:=\bP_\cF\cap\GL(V)$ for the invertible elements of the stabilizer of $\cF$.

    \item The \emph{parabolic dimension} of $\cF$ is
    \[\pdim(\cF):=\sum_{i=1}^m\pr*{(\dim V_i-\dim V_{i-1})^2-1}.\]

    \item Each element of $\bP_\cF$ descends to a linear map on $V_i/V_{i-1}$; this gives a map $\pi_i\colon\bP_\cF\to\End(V_i/V_{i-1})$. We define the \emph{vector-valued trace} on $\bP_\cF$ by
    \[\vtr_\cF(x):=(\tr\pi_1(x),\tr\pi_2(x),\ldots,\tr\pi_m(x)).\]
    Since each $\pi_i$ is an $\FF$-algebra homomorphism, $\vtr_\cF$ is a class function on $H_\cF$.

    \item For any affine subspace $\bV\subset\bP_\cF$, we define the \emph{$\cF$-dimension} of $\bV$ to be
    \[\dim_\cF\bV:=\sum_{i=1}^m(\dim\pi_i(\spn\bV)-1).\]
    We note that $\dim_\cF\bP_\cF=\pdim(\cF)$.
\end{itemize}
\end{definition}

We often form refinements of flags in the following way: if we have a flag $\cF=(V_0,V_1,\ldots,V_m)$ on some vector space $V$ and a flag $\cF_i=(W_0,\ldots,W_k)$ on some quotient $V_i/V_{i-1}$, the \emph{refinement of $\cF$ by $\cF_i$} is the flag given by inserting $V_{i-1}+\cF_i$ between $V_{i-1}$ and $V_i$; i.e.,
\begin{align*}
\{0\}=V_0\subsetneq V_1\subsetneq\cdots
&\subsetneq V_{i-1}=V_{i-1}+W_0\\
&\subsetneq V_{i-1}+W_1\subsetneq \cdots\subsetneq V_{i-1}+W_{k-1}\\
&\subsetneq V_{i-1}+W_k=V_i\subsetneq V_{i+1}\subsetneq\cdots\subsetneq V_k=V.
\end{align*}
If $\cF$ is the trivial flag, then $H_\cF=\GL_d(\FF)$, while if $\cF$ is a complete flag then $H_\cF$ is a Borel subgroup of $\GL_d(\FF)$. Most of our work will be done within subgroups $H_\cF$ for a flag $\cF$ which begins as the trivial flag and is refined throughout the iteration. The reader may wish to keep in mind the model case where $\cF$ is the trivial flag $\{0\}\subset\FF^d$.

The following lemma will be useful to us in the case when $A^3$ contains many elements of a fixed trace.

\begin{lemma}\label{lem:frobenius-inner-product} Let $\bV,\bW\subset\Mat_d(\FF)$ be affine subspaces and let $t\in\FF$. Suppose that $\tr(xy)=t$ for all $x\in\bV$ and $y\in\bW$. Then
\[\dim\spn(\bV)+\dim\spn(\bW)\leq d^2+1.\]
\end{lemma}
\begin{proof} The condition that $\tr(xy)=t$ for all $x\in\bV$ and $y\in\bW$ implies that $\tr(zw)=0$ for every $z\in\spn\bV$ and $w\in\bW-\bW$. Both $\spn\bV$ and $\bW-\bW$ are genuine subspaces of $\Mat_d(\FF)$, and they are orthogonal with respect to the Frobenius inner product. Since the Frobenius inner product is nondegenerate, we conclude that
\[\dim\spn(\bV)+\dim(\bW-\bW)\leq d^2.\]
The result now follows from the fact that $\spn\bW$ is spanned by $\bW-\bW$ and any fixed element $w_0$ of $\bW$, and so $\dim(\bW-\bW)+1\geq\dim\spn\bW$.
\end{proof}

We will need to deal with the multiplicative structure of affine subspaces of $\Mat_d(\FF)$. For affine subspaces $\bV,\bW\subset\Mat_d(\FF)$, define
\begin{align*}
\Lmul(\bW\to \bV)&:=\set{x\in\Mat_d(\FF):x\bW\subset \bV};\\
\Rmul(\bW\to \bV)&:=\set{x\in\Mat_d(\FF):\bW x\subset \bV}.
\end{align*}

For a non-commutative ring $R$, a \emph{bitranslate} of a subset $S\subset R$ is a set of the form $aSb$ for units $a,b\in R^\times$.

We will make use of the following three simple lemmas concerning $\Lmul$ and $\Rmul$.

\begin{lemma}\label{lem:mapping-is-affine} If $\bV$ and $\bW$ are affine subspaces of $\Mat_d(\FF)$, then $\Lmul(\bW\to \bV)$ and $\Rmul(\bW\to \bV)$ are affine subspaces of $\Mat_d(\FF)$.
\end{lemma}
\begin{proof} Suppose $x_1,x_2\in\Mat_d(\FF)$ are such that $x_1\bW,x_2\bW\subset \bV$. Then, for each $\lambda\in\FF$,
\[(\lambda x_1+(1-\lambda)x_2)\bW\subset \lambda x_1\bW+(1-\lambda)x_2\bW\subset \lambda \bV+(1-\lambda)\bV\subset \bV.\]
This implies that $\Lmul(\bW\to \bV)$ is an affine subspace of $\Mat_d(\FF)$; the proof for $\Rmul$ is identical.
\end{proof}

\begin{lemma}\label{lem:bitranslate} Let $\bV$ and $\bW$ be affine subspaces of $\Mat_d(\FF)$. If $\bW$ contains an invertible element, then $\Lmul(\bW\to\bV)$ is contained in a right translate of $\bV$. If $\Lmul(\bW\to\bV)$ contains an invertible element, then $\bW$ is contained in a left translate of $\bV$. The same two properties hold with $\Lmul$ replaced by $\Rmul$, with ``left'' and ``right'' swapped.
\end{lemma}
\begin{proof} Taking $a\in\bW$ (resp.\ $b\in\Lmul(\bW\to\bV)$) invertible, we have $\Lmul(\bW\to\bV)\subset\bV a^{-1}$ (resp.\ $\bW\subset b^{-1}\bV$). The proof for $\Rmul$ is identical.
\end{proof}

\begin{lemma}\label{lem:lmul-dim} Let $\bV$ and $\bW$ be proper affine subspaces of $\Mat_d(\FF)$, and suppose that $\bW$ contains an invertible element. Then
\begin{equation}\label{eq:lmul-dim}
\dim\Lmul(\bW\to\bV)\leq\dim\bV.
\end{equation}
Moreover, suppose equality holds in \eqref{eq:lmul-dim} and $\bV$ contains an invertible element. Then there exists a subalgebra $\bU$ of $\Mat_d(\FF)$ for which (i) $(1+\bU)\cap\GL_d(\FF)$ is a subgroup of $\GL_d(\FF)$,\footnote{This condition in fact follows from the fact that $\bU$ is a subalgebra of $\Mat_d(\FF)$. However, as it arises directly from the way in which we construct $\bU$ here, we prefer to simply place it in the statement of the lemma.} (ii) $\bW$ is contained in a left translate of $1+\bU$, and (iii) $1+\bU$ is contained in a left translate of $\bV$.
\end{lemma}

\begin{proof} Write $\bX:=\Lmul(\bW\to\bV)$ for notational simplicity. Let $a\in\bW$ be some invertible element. The map $\bX\to\bV$ given by $x\mapsto xa$ is an injective linear map, from which we conclude \eqref{eq:lmul-dim}. Now, suppose equality holds. We have
\[\bW\subset\Rmul(\bX\to\bV)=a^{-1}\Rmul(\bX a\to\bV)=a^{-1}\Rmul(\bV\to\bV),\]
where the last equality holds since $\dim\bX=\dim\bV$. Now, let $\bU:=\Rmul(\bV\to\bV)-1$. By \cref{lem:bitranslate}, $1+\bU$ is contained in a left translate of $\bV$. Since $\Rmul(\bV\to\bV)$ contains the identity, and is an affine subspace of $\Mat_d(\FF)$ by \cref{lem:mapping-is-affine}, $\bU$ is a genuine subspace of $\Mat_d(\FF)$. It is also closed under multiplication: if $x,y\in\bU$, then
\[x+y+xy=(1+x)(1+y)-1\in\Rmul(\bV\to\bV)-1=\bU,\]
and so $xy\in\bU$. Therefore $\bU$ is a subalgebra of $\Mat_d(\FF)$, the set $1+\bU$ is contained in a left translate of $\bV$, and $\bW$ is contained in $a^{-1}(1+\bU)$, as desired.
\end{proof}

\subsection{A partial proof of \texorpdfstring{\cref{thm:GL}}{Theorem 1.12}}\label{sec:pf-GL} We now prove \cref{thm:GL}, modulo \cref{prop:reg-el} and the following lemma.

\begin{lemma}[Subspace-to-flag lemma]\label{lem:subspace-to-subgroup-flag} Let $d$ be a positive integer and let $\cF=(V_0,V_1,\ldots,V_m)$ be a flag on $\FF^d$ of some parabolic dimension $D$. Let $A$ be a $K$-approximate subgroup of $H_\cF$ and let $\eta>0$. Suppose that, for some $\mathbf t\in\FF^m$, we have
\[\frac{\abs{A^2\cap\vtr_\cF^{-1}(\mathbf t)}}{\abs{A^2}}\geq\eta.\]
Then there exists some refinement $\cG$ of $\cF$ with parabolic dimension at most $\frac1{\sqrt2}D$ which satisfies
\[\frac{\abs{A^2\cap H_\cG}}{\abs{A^2}}\geq\pr*{\frac{\eta}{2K}}^{5d^6}.\]
\end{lemma}

This lemma will be one of the main results of \cref{sec:subspace-to-subgroup}, where it will follow from a more general ``subspace-to-subgroup'' lemma. It is what we use to treat the case arising from \cref{lem:pivot} in which $\abs{A^3\cap\Conj(a)}$ is large.

To prove \cref{thm:GL}, we will also need the following lemma, which tells us how much we gain in the case where $\abs{A^2\cap C(a)}$ is large.

\begin{lemma}\label{lem:centralizer} For every $a\in\GL_d(\FF)$, there exists some flag $\cF$ of parabolic dimension at most $d\cdot\eigm(a)-1$ for which $C(a)\subset H_\cF$.
\end{lemma}
\begin{proof} For each $\lambda$, let $t_\lambda$ denote the size of the largest Jordan block of $a$ with eigenvalue $\lambda$. For each $0\leq i\leq t_\lambda$, let $V_\lambda^i:=\ker((a-\lambda)^i)$ be the space of generalized eigenvectors of $a$ with eigenvalue $\lambda$ and rank at most $i$. For any $v\in V_\lambda^i$ and any $b\in C(a)$,
\[(a-\lambda)^ibv=b(a-\lambda)^iv=0,\]
so $bv\in V_\lambda^i$. Therefore $C(a)$ preserves each $V_\lambda^i$. 

We now construct the flag $\cF$. First, order the eigenvalues of $a$ arbitrarily as $\lambda_1,\ldots,\lambda_m$, and define the flag
\[\cF':=\pr*{\{0\},V_{\lambda_1}^{t_{\lambda_1}},V_{\lambda_1}^{t_{\lambda_1}}+V_{\lambda_2}^{t_{\lambda_2}},\ldots,V_{\lambda_1}^{t_{\lambda_1}}+V_{\lambda_2}^{t_{\lambda_2}}+\cdots+V_{\lambda_{m-1}}^{t_{\lambda_{m-1}}},\FF^d}\]
so that the quotients between consecutive subspaces in $\cF'$ are naturally identified with the generalized eigenspaces $V_{\lambda_j}^{t_{\lambda_j}}$. Now, successively refine $\cF'$ by the flags
\[\cF_j:=\pr*{\{0\},V_{\lambda_j}^1,V_{\lambda_j}^2,\ldots,V_{\lambda_j}^{t_{\lambda_j}}}\]
for each $1\leq j\leq m$ to form the flag $\cF$. Since $C(a)$ preserves each $V_\lambda^i$, we conclude that $C(a)\subset H_\cF$. Moreover, we have
\begin{equation}\label{eq:pdim-centralizer}
\pdim\cF=\sum_{\lambda\text{ eig.~of }a}\sum_{i=1}^{t_\lambda}\pr*{\dim(V_\lambda^i/V_\lambda^{i-1})^2-1}\leq d\pr*{\max_{\lambda,i}\dim(V_\lambda^i/V_\lambda^{i-1})}-1.
\end{equation}
Finally, the map $V_\lambda^{i+1}/V_\lambda^i\to V_\lambda^i/V_\lambda^{i-1}$ given by $v\mapsto (a-\lambda)v$ is injective for each $i\geq 1$, and so the dimension $\dim(V_\lambda^i/V_\lambda^{i-1})$ is maximized when $i=1$. For each $\lambda$, this dimension is at most $\eigm(a)$. So, the result follows from \eqref{eq:pdim-centralizer}.
\end{proof}

Finally, we will need a corollary of \cref{prop:reg-el}.

\begin{lemma}\label{lem:reg-el-spec} Let $d\geq 8000$ and let $A\subset\GL_d(\FF)$ be a $K$-approximate subgroup. Let $\delta\in(0,1]$ be a parameter. Then one of the following holds:
\begin{itemize}
    \item (Not-too-irregular elements) For more than $(1-\delta)\abs{A^2}$ elements $a\in A^2$, it holds that $\eigm(a)\leq d(1-\frac1{1200\log\log d})$.

    \item (Control by well-behaved subgroup) There exist two subspaces $V_1\subset V_2\subset\FF^d$ with $\dim V_2>\frac{9d}{10}$ and $\dim V_1<\frac d{10}$ for which at least $(2K^8\delta^{-1})^{-3d^6}\abs{A^2}$ elements $a$ of $A^2$ satisfy the following property: there exists some $\lambda\in\FF$ for which $(a-\lambda)V_2\subset V_1$.
\end{itemize}
\end{lemma}
\begin{proof} Let $B$ be the set of elements $a\in A^2$ satisfying $\eigm(a)\leq d(1-\frac1{1200\log\log d})$; the set $B$ is centered. If $\abs{B}\leq\delta\abs{A^2}$, the not-too-irregular elements conclusion holds. Otherwise, Ruzsa's covering lemma implies that $B^4$ can be covered by at most $\abs{B^5}/\abs{B}\leq K^8\delta^{-1}$ translates of $B^2$, and so $B^2$ is a $K^8\delta^{-1}$-approximate group. Moreover, every element of $B^2$ has an eigenspace of codimension at most $\frac d{600\log\log d}$. We can thus apply \cref{prop:reg-el} to $B^2$ to obtain that there exist two subspaces $V_1\subset V_2\subset\FF^d$ with $\dim V_2>\frac{9d}{10}$ and $\dim V_1<\frac d{10}$ such that the following holds: if $H\leq\GL_d(\FF)$ is the subgroup consisting of those elements $g$ for which there exists a scalar $\lambda\in\FF$ with $(g-\lambda)V_2\subset V_1$, then
\[\abs{B^4\cap H}\geq(2K^8\delta^{-1})^{-2d^6}\abs{B^4}.\]
We conclude from \cref{lem:intersect-subgroup} that
\[\abs{A^2\cap H}\geq K^{-7}\abs{A^8\cap H}\geq K^{-7}\abs{B^4\cap H}\geq K^{-7}(2K^8\delta^{-1})^{-2d^6}\abs{B^4}\geq K^{-7}(2K^8\delta^{-1})^{-2d^6}\delta\abs{A^2}.\]
The result follows.
\end{proof}

We now embark on the proof of \cref{thm:GL}. By \cref{lem:borel-abelian-struct}, it is enough to find a Borel subgroup of $\GL_d(\FF)$ which has large intersection with $A^2$. We proceed by an iterative argument through the subgroups $H_\cF\leq\GL_d(\FF)$, refining $\cF$ at each step. The following technical lemma describes the iteration.

\begin{lemma}\label{lem:big-GL-iteration} Let $A$ be a $K$-approximate subgroup of $\GL_d(\FF)$, let $\cF$ be a flag on $\FF^d$ of parabolic dimension $D>0$, and let $\lambda,\eps>0$ be such that
\[\abs{A^2\cap H_\cF}\geq\frac1{(2K)^\lambda}\abs{A^2}^\eps.\] Write $\gamma=\min(\frac1{8000},\frac1{1200\log\log d})$ if $d>2$ and $\gamma=\frac14$ otherwise. For some triple $(D',\lambda',\eps')$ satisfying
\[(D',\lambda',\eps')\in\set*{
\pr*{D-1,\lambda+85d^7,\eps},\ \
\pr*{\floor*{\frac1{\sqrt2}D},\lambda+165d^6,\frac\eps2},\ \
\pr*{\floor*{(1-\gamma)D},7,\frac{\eps}{10d^6}}}\]
there exists a refinement $\cG$ of $\cF$ satisfying $\pdim\cG\leq D'$ and $\abs{A^2\cap H_\cG}\geq\frac1{(2K)^{\lambda'}}\abs{A^2}^{\eps'}$.
\end{lemma}
\begin{proof}[Proof of \cref{lem:big-GL-iteration} assuming \cref{lem:subspace-to-subgroup-flag,prop:reg-el}] Write $\cF=(V_0,V_1,\ldots,V_m)$, and let $\pi_i\colon\bP_\cF\to\End(V_i/V_{i-1})$ be the natural ring homomorphism. Write $d_i:=\dim V_i-\dim V_{i-1}$, so that $d_1+\cdots+d_m=d$. For each $i$, write $Z_i$ for the center of $\GL(V_i/V_{i-1})$. Finally, select $\delta:=(dK^{12})^{-1}$.

For each $1\leq i\leq m$, the set $B_i:=\pi_i(A^2\cap H_\cF)\subset\GL(V_i/V_{i-1})$ is a $K^3$-approximate group, by \cref{lem:intersect-subgroup}. So, for each $i$, one of the following holds (in the large $d_i$ case, we apply \cref{lem:reg-el-spec}):

\begin{enumerate}[(1)]
    \item $d_i=1$;
    \item $1<d_i<8000$, and
    \begin{enumerate}[(a)]
        \item $\abs{B_i^2\cap Z_i}\geq\delta\abs{B_i^2}$, or
        \item for strictly more than $(1-\delta)\abs{B_i^2}$ elements $b\in B_i^2$, we have $\eigm(b)\leq d_i-1$; or
    \end{enumerate}
    \item $d_i\geq 8000$, and
    \begin{enumerate}[(a)]
        \item there exists some flag $\cF_i$ on $V_i/V_{i-1}$ which satisfies\footnote{In the notation of \cref{lem:reg-el-spec}, this flag is obtained by refining the flag $\{0\}\subsetneq V_1\subsetneq V_2\subsetneq \FF^{d_i}$ by any complete flag on $V_2/V_1$. The bound $\pdim\cF_i\leq\frac1{50}(d_i^2-1)$ follows from the fact that $\dim V_1<\frac{d_i}{10}<\frac{9d_i}{10}<\dim V_2$.} $\pdim\cF_i\leq\frac1{50}(d_i^2-1)$ and
        \[\abs{B_i^2\cap H_{\cF_i}}\geq(2K^{24}\delta^{-1})^{-2d_i^6}\abs{B_i^2},\]
        or

        \item for strictly more than $(1-\delta)\abs{B_i^2}$ elements $b\in B_i^2$, we have $\eigm(b)\leq d_i\pr*{1-\frac1{1200\log\log d_i}}$.
    \end{enumerate}
\end{enumerate}

\noindent \textbf{If (2a) or (3a) holds for some $i$:}
In either case, we have some $1\leq i\leq m$ and some flag $\cF_i$ on $V_i/V_{i-1}$ (in case (2a), this is a complete flag) which satisfies $\pdim\cF_i\leq d_i^2-2$ and
\[\frac{\abs{B_i^2\cap H_{\cF_i}}}{\abs{B_i^2}}\geq\frac{\delta^{2d^6}}{(2K)^{48d^6}}.\]
Apply \cref{lem:intersect-subgroup-projection} to the $K^3$-approximate subgroup $A^2\cap H_\cF$ and the projection map $\pi_i$ to get that, in fact,
\[\frac{\delta^{2d^6}}{(2K)^{48d^6}}\leq \frac{\abs{(A^2\cap H_\cF)^6\cap\pi_i^{-1}(H_{\cF_i})}}{\abs{(A^2\cap H_\cF)^2}}\leq\frac{\abs{A^{12}\cap\pi_i^{-1}(H_{\cF_i})}}{\abs{A^2\cap H_\cF}}.\]
Let $\cG$ be the flag given by refining $\cF$ by $\cF_i$, so that $H_\cG=\pi_i^{-1}(H_{\cF_i})$. Using \cref{lem:intersect-subgroup}, we obtain
\[\frac{\delta^{2d^6}}{(2K)^{48d^6}}\leq K^{11}\frac{\abs{A^2\cap H_\cG}}{\abs{A^2\cap H_\cF}}.\]
Using our choice of $\delta$ and the inequality $2^d>d$, we have
\[\abs{A^2\cap H_\cG}\geq\frac{\delta^{2d^6}}{(2K)^{48d^6+11}}\abs{A^2\cap H_\cF}\geq\frac1{(2K)^{72d^6+2d^7+11}}\abs{A^2\cap H_\cF}\geq\frac1{(2K)^{85d^7+\lambda}}\abs{A}^\eps.\]
This gives the result with $D'=D-1$, $\lambda'=\lambda+85d^7$, and $\eps'=\eps$.

\vspace{2mm}

\noindent \textbf{If (1), (2b), or (3b) holds for each $i$:}
For each $i$, define
\[e_i:=\begin{cases}1&\text{if }d_i=1\\d_i-1&\text{if }d_i<8000\\d_i-\frac{d_i}{1200\log\log d_i}&\text{if }d_i\geq 8000,\end{cases}\]
and let $S_i:=\set{b\in\GL(V_i/V_{i-1}):\eigm(b)>e_i}$. For every $i$, we have $\abs{B_i^2\cap S_i}/\abs{B_i^2}<\delta$. Applying \cref{lem:intersect-subgroup-projection} to the $K^3$-approximate subgroup $A^2\cap H_\cF$ and the projection map $\pi_i$, we get that
\[\frac{\abs{(A^2\cap H_\cF)^2\cap\pi^{-1}(S_i)}}{\abs{(A^2\cap H_\cF)^2}}\leq K^{12}\frac{\abs{B_i^2\cap S_i}}{\abs{B_i^2}}<K^{12}\delta=\frac1d.\]
Since $m\leq d$, a union bound gives that there exists some $a_0\in(A^2\cap H_\cF)^2$ for which $\pi(a_0)\not\in S_i$ for any $1\leq i\leq m$, i.e.~that $\eigm(\pi(a_0))\leq e_i$ for each $i$. Let $B:=(A^2\cap H_\cF)^2$. Since $A^2\cap H_\cF$ is a $K^3$-approximate group, $B$ is a $K^6$-approximate group.

\vspace{2mm}

We now apply \cref{lem:pivot} to find that
\[\abs{B^2\cap C(a_0)}\cdot\abs{B^3\cap\Conj(a_0)}\geq\abs{B}.\]
Let $\eta:=\abs{A}^{-\eps/(10d^6)}$ be a parameter. Either
\begin{enumerate}[(i)]
    \item $\abs{B^3\cap\Conj(a_0)}\geq\eta\abs{B}$ or 
    \item $\abs{B^2\cap C(a_0)}\geq1/\eta$.
\end{enumerate}

\noindent \textbf{If case (i) holds:}
Let $\mathbf t:=\vtr_\cF(a_0)$ so that $\Conj(a_0)\subset\vtr_\cF^{-1}(\mathbf t)$. Let $\tilde B:=B^2$ be a $K^{12}$-approximate group. We have
\[\frac{\abs{\tilde B^2\cap\vtr_{\cF}^{-1}(\mathbf t)}}{\abs{\tilde B^2}}\geq\frac{\abs{B^3\cap\Conj(a_0)}}{\abs{B^4}}\geq\frac\eta{K^{18}}.\]
So, we may apply \cref{lem:subspace-to-subgroup-flag} to the $K^{12}$-approximate subgroup $\tilde B$ of $H_\cF$ to find some refinement $\cG$ of $\cF$ with parabolic dimension at most $\frac1{\sqrt2}D$ which satisfies
\[\frac{\abs{\tilde B^2\cap H_\cG}}{\abs{\tilde B^2}}\geq\pr*{\frac{\eta}{2K^{30}}}^{5d^6}.\]
Unwinding this inequality into a statement about $A$ using \cref{lem:intersect-subgroup} gives
\[\pr*{\frac{\eta}{2K^{30}}}^{5d^6}\leq \frac{\abs{\tilde B^2\cap H_\cG}}{\abs{\tilde B^2}}\leq\frac{\abs{A^{16}\cap H_\cG}}{\abs{A^2\cap H_\cF}}\leq\frac{K^{15}\abs{A^2\cap H_\cG}}{\abs{A^2\cap H_\cF}}.\]
We conclude that
\[\abs{A^2\cap H_\cG}\geq\frac{\eta^{5d^6}}{(2K)^{150d^6+15}}\abs{A^2\cap H_\cF}\geq\frac1{(2K)^{165d^6+\lambda}}\abs{A}^{\eps/2}.\]
This gives the result with $D'=\frac1{\sqrt2}D$, $\lambda'=\lambda+165d^6$, and $\eps'=\eps/2$.

\vspace{2mm}

\noindent \textbf{If case (ii) holds:}
For each $i$, by \cref{lem:centralizer}, $\pi_i(C_{H_\cF}(a_0))=C_{\GL(V_i/V_{i-1})}(\pi_i(a_0))$ is contained within $H_{\cF_i}$ for some flag $\cF_i$ of parabolic dimension at most $d_ie_i-1$. Our choice of $\gamma=1/4$ if $d\leq 2$ and $\gamma=\min(\frac1{8000},\frac1{1200\log\log d})$ otherwise is enough to guarantee that $e_i\leq (1-\gamma)d_i$ whenever $d_i>1$. We thus have $d_ie_i-1\leq(1-\gamma)(d_i^2-1)$ for every $i$ (including those $i$ satisfying $d_i=1$). As a result, the refinement $\cG$ of $\cF$ by all of the $\cF_i$ satisfies
\[\pdim\cG=\sum_{i=1}^m\pdim\cF_i\leq (1-\gamma)\pdim\cF\]
and $H_\cG=\bigcap_{i=1}^m\pi_i^{-1}(H_{\cF_i})$. This means that $C(a_0)\subset H_\cG$. We conclude, using \cref{lem:intersect-subgroup}, that
\[\abs{A^2\cap H_\cG}\geq\abs{A^2\cap C(a_0)}\geq K^{-7}\abs{A^8\cap C(a_0)}\geq K^{-7}\abs{B^2\cap C(a_0)}\geq\frac{1/\eta}{K^7}=K^{-7}\abs{A}^{\eps/(10d^6)}.\]
This gives the result with $D'=(1-\gamma)D$, $\lambda'=7$, and $\eps'=\eps/(10d^6)$.
\end{proof}

Finally, we can use the iteration described in \cref{lem:big-GL-iteration} to prove \cref{thm:GL}.

\begin{proof}[Proof of \cref{thm:GL} assuming \cref{lem:subspace-to-subgroup-flag,prop:reg-el}] If $d=1$ the result follows from taking $H=\GL_1(\FF)$, so we henceforth assume $d>1$. We begin with a $K$-approximate group $A\subset\GL_d(\FF)$. Initialize $\cF_0:=(\set{0},\FF^d)$ to be the trivial flag. Also initialize $D_0:=d^2-1$ and $\lambda_0:=0$ and $\eps_0:=1$. From $(\cF_i,D_i,\lambda_i,\eps_i)$, we will find $(\cF_{i+1},D_{i+1},\lambda_{i+1},\eps_{i+1})$, maintaining the properties that $\pdim\cF_i\leq D_i$ and
\begin{equation}\label{eq:intersect-w-flag}
\abs{A^2\cap H_{\cF_i}}\geq\frac1{(2K)^{\lambda_i}}\abs{A^2}^{\eps_i}.
\end{equation}
The construction of $(\cF_{i+1},D_{i+1},\lambda_{i+1},\eps_{i+1})$ is given by applying \cref{lem:big-GL-iteration} to $(\cF_i,D_i,\lambda_i,\eps_i)$. The sequence $\set{D_i}$ is strictly decreasing in $i$, and so we have $D_\ell=0$ for some $\ell<d^2$; we terminate the process the first time this occurs. The condition $D_\ell=0$ implies that $\cF_\ell$ is a complete flag, and thus that $H_{\cF_\ell}$ is Borel. By \cref{lem:intersect-subgroup}, $A^2\cap H_{\cF_\ell}$ is a $K^3$-approximate subgroup of $\GL_d(\FF)$. We can therefore use \cref{lem:borel-abelian-struct} to find some abelian subgroup $H$ of $\GL_d(\FF)$ for which, by \eqref{eq:intersect-w-flag},
\[\abs{A^2\cap H}\geq K^{-3}\abs{(A^2\cap H_{\cF_\ell})^2\cap H}\geq K^{-102}\abs{A^2\cap H_{\cF_\ell}}^{1/d^4}\geq\frac1{(2K)^{\lambda_\ell+102}}\abs{A^2}^{d^{-4}\eps_\ell}.\]

What remains is thus to upper-bound $\lambda_\ell$ and lower-bound $\eps_\ell$.

To bound $\lambda_\ell$, we see that $\lambda_{i+1}-\lambda_i\leq 165d^7$ for each $i$. Therefore $\lambda_\ell=\lambda_\ell-\lambda_0\leq 165d^9$. To bound $\eps_\ell$, we see that $\eps_{i+1}=\eps_i$ whenever $D_{i+1}=D_i-1$. Therefore, whenever $\eps_{i+1}<\eps_i$, we have $D_{i+1}/D_i\leq 1-\gamma$, where $\gamma=\min(\frac1{8000},\frac1{1200\log\log d})$ if $d>2$ and $\gamma=1/4$ otherwise. Moreover, $\eps_i$ may decrease only by a factor of $10d^6$ at each such step. Since $D_{\ell-1}\geq 1$, we conclude that
\[\eps_\ell\geq \frac{\eps_{\ell-1}}{10d^6\eps_0}\geq\pr*{\frac1{10d^6}}^{1+\log_{1-\gamma}(d^2)}\geq\exp\pr*{-\frac{80\log^2 d}\gamma}.\]
This is enough to give the result.
\end{proof}

\section{Not-too-irregular elements in matrix groups}\label{sec:not-too-irreg}

The object of this section is to prove the following weak version of \cref{prop:reg-el}. We write $\rk m$ for the rank of a matrix $m$.\footnote{Throughout this section, we endow every vector space with an inner product structure so that we may freely discuss orthogonal complements. This is only for notational convenience; none of our results have any material dependence on the inner product structure chosen.}

\begin{lemma}\label{lem:reg-el} Let $d,s\geq 6$ and let $A\subset\GL_d(\FF)$ be a $K$-approximate subgroup. Suppose that every element $a$ of $A^2$ satisfies $\rk(a-1)\leq s$. Then there exist subspaces $U$ and $W$ of $\FF^d$ satisfying $\dim U,\dim W\leq30s\log\log s$ for which at least $(2K)^{-60s}\abs{A^2}$ elements $a\in A^2$ satisfy $(a-1)W^\perp\subset U$.
\end{lemma}

We will upgrade this statement to \cref{prop:reg-el} in the next section.

Our proof of \cref{lem:reg-el} roughly follows an iterative procedure developed by Haramaty and Shpilka \cite{HaramatyShpilka} in an additive setting. The clearest analogue of our approach in the literature is the following lemma of Lampert, which is quite related to \cite[Lemma~3.5]{HaramatyShpilka}.

\begin{lemma}[{Weak version of \cite[Lemma~2.3]{Lampert}}]\label{lem:haramaty-shpilka-lampert} Let $\bX\subset\Mat_d(\FF)$ be a subspace, and suppose that $\rk x\leq r$ for every $x\in\bX$. Then there exist subspaces $U$ and $W$ of $\FF^d$ satisfying $\dim U,\dim W\leq 2r$ such that $\bX W^\perp\subset U$.
\end{lemma}

Our situation in \cref{lem:reg-el} is similar to that of \cref{lem:haramaty-shpilka-lampert}. The main difference is that, instead of having a set of low-rank operators with rich additive structure, our set $A$ has rich \emph{multiplicative} structure. This complicates our argument significantly, and is ultimately the source of the $\log\log s$ terms in \cref{lem:reg-el}.

We begin with the following simple fact from linear algebra, which is essentially a ``coordinate-free'' version of a computation appearing in the proof of \cite[Lemma~3.5]{HaramatyShpilka}.

\begin{lemma}\label{lem:addition-rank} Let $x,y\colon U\to W$ be linear operators. We have
\[\dim(\im x+\im y)+\dim\pr*{(\ker x)^\perp+(\ker y)^\perp}\leq \rk x+\rk y+\rk (x+y).\]
\end{lemma}
\begin{proof} The map $x$ factors through a vector space of dimension $\rk x$. We can thus find some vector space $V_x$ of dimension $\rk x$, some surjection $m_x\colon U\twoheadrightarrow V_x$, and some injection $n_x\colon V_x\hookrightarrow W$ for which $x=n_xm_x$. Define $V_y$, $m_y$, and $n_y$ similarly. Let $V:=V_x\oplus V_y$, and consider the maps $m\colon U\to V$ and $n\colon V\to W$ defined by $mu=(m_xu,m_yu)$ and $n(v_x,v_y)=n_xv_x+n_yv_y$. We have $nm=x+y$. As a result, we may write
\begin{align*}
\dim\ker(x+y)
&=\dim\ker(nm)\\
&\leq\dim\ker n+\dim\ker m\\
&=\dim V-\dim\im n+\dim(\ker m_x\cap \ker m_y)\\
&=\dim V_x+\dim V_y-\dim(\im n_x+\im n_y)+\dim(\ker x\cap\ker y)\\
&=\rk x+\rk y-\dim(\im x+\im y)+\dim U-\dim\pr*{(\ker x\cap\ker y)^\perp}.
\end{align*}
The left side of the above is $\dim U-\rk(x+y)$. So, using the fact that $(U_1\cap U_2)^\perp=U_1^\perp+U_2^\perp$ for subspaces $U_1$ and $U_2$ of $U$ and rearranging gives the result.
\end{proof}

\begin{corollary}\label{cor:addition-rank} Let $x,y\colon U\to W$ be linear operators. Let $\iota$ be the inclusion $\ker y\hookrightarrow U$ and $\pi$ be the projection $U\twoheadrightarrow U/\im y$. Then
\[2\rk(\pi x\iota)\leq \rk(\pi x)+\rk(x\iota)\leq \rk(x+y)+\rk x-\rk y.\]
\end{corollary}
\begin{proof} Note that $\rk y=\dim\im y=\dim(\ker y)^\perp$. We compute
\begin{align*}
\rk(\pi x)
&=\dim(\im\pi x)=\dim(\im x/\im y)=\dim(\im x+\im y)-\dim(\im y);\\
\rk(x\iota)&=\dim(\ker y)-\dim((\ker x)\cap(\ker y))=\dim\pr[\big]{(\ker x)^\perp+(\ker y)^\perp}-\dim(\ker y)^\perp.
\end{align*}
Combining these with \cref{lem:addition-rank} gives the result.
\end{proof}

We now outline a proof of \cref{lem:haramaty-shpilka-lampert}, as it will be helpful in understanding our approach. Suppose $\bX\subset\End(V)$ is a subspace. One may show the following statement by induction on $k$: there exist vector spaces $U_k$ and $W_k$ of dimension at least $d-2r(1-2^{-k})$, an injection $\iota_k\colon U_k\hookrightarrow V$, and a surjection $\pi_k\colon V\twoheadrightarrow W_k$ such that $\rk \pi_kx\iota_k\leq 2^{-k}r$ for each $x\in\bX$. \cref{lem:haramaty-shpilka-lampert} then follows from taking $k=\ceil{\log_2r}$. The base case $k=0$ is exactly the hypothesis of \cref{lem:haramaty-shpilka-lampert}, and the inductive step comes from applying \cref{cor:addition-rank}: choose $y\in\bX$ to maximize $\rk\pi_ky\iota_k$, define
\[U_{k+1}:=(\ker\pi_ky\iota_k)^\perp,\ \ W_{k+1}:=W_k/\im(\pi_ky\iota_k),\]
and let $\iota_{k+1}$ and $\pi_{k+1}$ be the natural compositions $U_{k+1}\hookrightarrow U_k\hookrightarrow V$ and $V\twoheadrightarrow W_k\twoheadrightarrow W_{k+1}$, respectively. Since $\rk\pi_ky\iota_k\leq 2^{-k}r$, the spaces $U_{k+1}$ and $W_{k+1}$ have codimension at most $2r(1-2^{-(k+1)})$. Finally, \cref{cor:addition-rank}, along with our choice of $y$, is enough to conclude the necessary rank bound on $\pi_{k+1}x\iota_{k+1}$.

Our proof of \cref{lem:reg-el} follows a similar strategy, with an iteration step centered around \cref{cor:addition-rank}.

\begin{lemma}\label{lem:reg-iteration} Let $U$, $V$, and $W$ be vector spaces, let $K,r,s,\eps,\eta>0$, let $A$ be a $K$-approximate subgroup of $\GL(V)$, and let $\iota\colon U\hookrightarrow V$ and $\pi\colon V\twoheadrightarrow W$ be injective and surjective linear maps, respectively. Suppose $S\subset A^2$ satisfies the following properties:
\begin{enumerate}[(i)]
    \item $\abs{S}\geq\eps\abs{A^2}$,
    \item for each $a\in S$, we have $\rk(\pi(a-1))\leq s$, and
    \item for each $a\in S$, we have $\rk(\pi(a-1)\iota)\leq r$.
\end{enumerate}
Then there exist vector spaces $U'$ and $W'$, maps $\iota'\colon U'\hookrightarrow V$ and $\pi'\colon V\twoheadrightarrow W'$, a subset $S'\subset A^2$, and a triple $(\rho,\ell,m)\in\set{(1-\eta,r,1),(3\eta,2s,2)}$ such that
\begin{enumerate}[(1)]
    \item $\iota'$ is injective and $\pi'$ is surjective,
    \item both $\dim U-\dim U'$ and $\dim W-\dim W'$ are at most $\ell$,
    \item $\abs{S'}\geq\frac{\eps^m}{4K}\abs{A^2}$, and
    \item for every $a\in S'$, we have $\rk(\pi'(a-1)\iota')\leq\rho r$.
\end{enumerate}
\end{lemma}

\begin{proof} We begin by defining the following sets.
\begin{align*}
S_0&:=\set[\big]{a\in S:\rk(\pi(a-1)\iota)\leq r(1-\eta)};\ \ S_1:=S\setminus S_0.\\
T&:=\set[\big]{(b,c)\in S_1\times S_1:bc^{-1}\in A^2}.\\
T_0&:=\set[\big]{(b,c)\in T:\rk(\pi(b-c)\iota)\leq r(2-3\eta)};\\
T_1&:=\set[\big]{(b,c)\in T:\rk(\pi(b-c)\iota)>r(2-3\eta)}.
\end{align*}

\noindent \textbf{If $\abs{S_0}\geq\frac14\abs{S}$:}
We can take $(\rho,\ell,m)=(1-\eta,r,1)$. Select $S'=S_0$; we can simply choose $(\pi',\iota',U',W')=(\pi,\iota,U,W)$.

\vspace{2mm}

\noindent \textbf{If $\abs{S_0}<\frac14\abs{S}$:} We have $\abs{S_1}>\frac34\abs{S}$. We can lower-bound the size of $T$ using the Cauchy--Schwarz inequality. Let $X$ be such that $\abs{X}\leq K$ and $A^2\subset AX$, and let $f\colon A^2\to X$ be defined such that $a\in Af(a)$ for each $a\in A^2$. We have
\begin{align*}
\abs{T}
&\geq\#\set{(b,c)\in S_1\times S_1:f(b)=f(c)}\\
&=\sum_{x\in X}\abs{f^{-1}(x)\cap S_1}^2\\
&\geq\frac1K\pr*{\sum_{x\in X}\abs{f^{-1}(x)\cap S_1}}^2=\frac{\abs{S_1}^2}K\geq \frac3{4K}\abs{S}\abs{S_1}.
\end{align*}
As a result, we have
\begin{equation}\label{eq:T-dichotomy}
\max\pr[\big]{\abs{T_0},\abs{T_1}}\geq\frac3{8K}\abs{S}\abs{S_1}.
\end{equation}

\noindent\textbf{If $\abs{T_0}\geq\abs{T_1}$:}
We take $(\rho,\ell,m)=(1-\eta,r,1)$. By \eqref{eq:T-dichotomy}, there exists some $c\in S_1$ for which $S':=\set{bc^{-1}:(b,c)\in T_0}\subset A^2$ has size at least $\frac3{8K}\abs{S}\geq\frac{\eps}{3K}\abs{A^2}$. Write $y:=\pi(c-1)\iota$. Let $\iota_c$ be the inclusion $\ker y\hookrightarrow U$, let $\pi_c$ be the projection $W\twoheadrightarrow W/\im y$, and define
\[U':=\ker y,\ W':=W/\im y,\ \iota':=c\iota\iota_c,\ \pi':=\pi_c\pi.\]
Since $c\in S_1$, we have $r(1-\eta)\leq\rk y\leq r$. This implies that $\dim U-\dim U'\leq r$ and $\dim W-\dim W'\leq r$. So, it is enough to show that
\[\rk\pr*{\pi'(a-1)\iota'}\leq r(1-\eta)\]
for every $a\in S'$. Indeed, let $b=ac$, so that $(b,c)\in T_0$. Write $x:=\pi(a-1)c\iota=\pi(b-c)\iota$, and apply \cref{cor:addition-rank} to $x$ and $y$ to get
\begin{equation}\label{eq:case-2-rank-bound}
2\rk(\pi'(a-1)\iota')=2\rk(\pi_cx\iota_c)\leq\rk(x+y)+\rk x-\rk y.
\end{equation}
Since $(b,c)\in T_0$, we have $\rk x=\rk(\pi(b-c)\iota)\leq r(2-3\eta)$, and since $x+y=\pi(b-1)\iota$ with $b\in S$, we have $\rk(x+y)\leq r$. We conclude from \eqref{eq:case-2-rank-bound} that
\[2\rk\pr*{\pi'(a-1)\iota'}\leq r(2-3\eta)+r-r(1-\eta)=2r(1-\eta),\]
which is enough.

\vspace{2mm}

\noindent\textbf{If $\abs{T_1}>\abs{T_0}$:}
We take $(\rho,\ell,m)=(3\eta,2s,2)$. Since $\abs{T_1}=\max(\abs{T_0},\abs{T_1})\geq\frac3{8K}\abs{S}\abs{S_1}$, there exists some $a\in SS^{-1}\cap A^2$ for which at least
\[\frac{\frac3{8K}\abs{S}\abs{S_1}}{\abs{A^2}}\geq\frac9{32K}\frac{\abs{S}^2}{\abs{A^2}}>\frac{\eps^2}{4K}\abs{A^2}\]
pairs $(b,ab)$ lie in $T_1$. Let $S':=\set{b\in S:(b,ab)\in T_1}$. Write $z:=\pi(a-1)$. Let $\pi_a$ be the projection $W\to W/\im z$ and set
\[U':=U,\ W':=W/\im z,\ \iota':=\iota,\ \pi':=\pi_a\pi.\]
Since $a\in SS^{-1}$ can be written as $bc^{-1}$ for some $b,c\in S$, property (ii) implies that
\[\rk z=\rk\pi(b-c)c^{-1}=\rk\pi(b-c)\leq\rk\pi(b-1)+\rk\pi(c-1)\leq 2s,\]
and so $\dim W-\dim W'\leq 2s$. Thus what remains to prove is that, for every $b\in S'$, we have $\rk(\pi'(b-1)\iota')\leq 3\eta r$. Set $x:=\pi(b-1)\iota$ and $y:=\pi(a-1)b\iota$. Applying \cref{cor:addition-rank} to $x$ and $y$ gives
\begin{equation}\label{eq:case-3-rank-bound}
\dim((\im x)/(\im y))\leq \rk(x+y)+\rk x-\rk y.
\end{equation}
Since $b\in S$, we have $\rk x\leq r$. Since $(b,ab)\in T_1$, we have $\rk y=\rk \pi(ab-b)\iota>(2-3\eta)r$. Also, $(b,ab)\in T_1$ implies $ab\in S$, so we have $\rk(x+y)=\rk\pi(ab-1)\iota\leq r$. Lastly, since $y=zb\iota$, we have $\im y\subset\im z$, so
\[\dim(\im x/\im y)\geq\dim(\im x/\im z)=\rk(\pi_ax)=\rk\pr*{\pi_a\pi(b-1)\iota}=\rk\pr*{\pi'(b-1)\iota'}.\]
Combining this with \eqref{eq:case-3-rank-bound} gives
\[\rk\pr*{\pi'(b-1)\iota'}\leq \rk(x+y)+\rk x-\rk y\leq r+r-(2-3\eta)r=3\eta r\]
for every $b\in S'$, as desired.
\end{proof}

We now have most of what we need to prove \cref{lem:reg-el}. What remains is to describe exactly how we iterate \cref{lem:reg-iteration}. For this, we need the following technical lemma, which describes an asymptotically optimal iteration procedure.

\begin{lemma}\label{lem:reg-iteration-technical} Let $V$ be a vector space and let $A$ be a $K$-approximate subgroup of $\GL(V)$. Suppose that $\rk(a-1)\leq s$ for every element $a\in A^2$. 
For every positive integer $k$, the following holds. For $\alpha_k:=2^{1-2^{k/2}}$, there exist spaces $U_k$ and $W_k$ of dimension at least $(\dim V)-4ks$, an injection $\iota_k\colon U_k\hookrightarrow V$, and a surjection $\pi_k\colon V\twoheadrightarrow W_k$ such that, for at least $(4K)^{-30/\alpha_k^{1/2}}\abs{A^2}$ elements $a$ of $A^2$, one has $\rk(\pi_k(a-1)\iota_k)\leq \alpha_ks$.
\end{lemma}
\begin{proof} Write $\alpha_k:=2^{1-2^{k/2}}$, define a sequence $\set{t_k}$ by $t_0=0$ and $t_{k+1}=2t_k+\frac{6\alpha_k}{\alpha_{k+1}}-5$. It may be manually verified that $t_k\leq30\alpha_k^{-1/2}$ for each $k\geq 0$. (For large $k$, one has $t_k\sim \alpha_k^{-0.292\ldots}=\alpha_k^{-(1-2^{-1/2})}$.) So, it suffices to prove the result of the lemma with $(4K)^{-30/\alpha_k^{1/2}}$ replaced by $(4K)^{-t_k}$. We do this by induction on $k$. For $k=0$ we may take $U_k=W_k=V$ and $\iota_k$ and $\pi_k$ to be the identity maps.

Now, suppose we have $(U_k,W_k,\iota_k,\pi_k)$ and we wish to construct $(U_{k+1},W_{k+1},\iota_{k+1},\pi_{k+1})$. To do this, we repeatedly apply \cref{lem:reg-iteration}. Initialize
\[\pr*{U^{(0)},W^{(0)},\iota^{(0)},\pi^{(0)},r^{(0)},t^{(0)}}=\pr*{U_k,W_k,\iota_k,\pi_k,\alpha_ks,t_k}.\]
Given $(U^{(i)},W^{(i)},\iota^{(i)},\pi^{(i)},r^{(i)},t^{(i)})$, we apply \cref{lem:reg-iteration} with $\eps=(4K)^{-t^{(i)}}$ and $\eta=\frac{\alpha_{k+1}s}{3r^{(i)}}$. Unrolling the result, we obtain some state
\[\pr*{U^{(i+1)},W^{(i+1)},\iota^{(i+1)},\pi^{(i+1)},r^{(i+1)},t^{(i+1)}}\]
for which at least $(4K)^{-t^{(i+1)}}\abs{A^2}$ elements $a$ of $A^2$ satisfy
\[\rk\pr*{\pi^{(i+1)}(a-1)\iota^{(i+1)}}\leq r^{(i+1)},\]
and our parameters are related either by
\begin{enumerate}[(a)]
    \item
    \begin{itemize}
        \item $t^{(i+1)}=t^{(i)}+1$,
        \item $r^{(i+1)}=(1-\eta)r^{(i)}=r^{(i)}-\frac13\alpha_{k+1}s$, and
        \item $\dim U^{(i)}-\dim U^{(i+1)},\dim W^{(i)}-\dim W^{(i+1)}\leq r^{(i)}$; or
    \end{itemize}
    \item
    \begin{itemize}
        \item $t^{(i+1)}=2t^{(i)}+1$,
        \item $r^{(i+1)}=3\eta r^{(i)}=\alpha_{k+1}s$, and
        \item $\dim U^{(i)}-\dim U^{(i+1)},\dim W^{(i)}-\dim W^{(i+1)}\leq 2s$.
    \end{itemize}
\end{enumerate}
(Moreover, all such elements satisfy $\rk(\pi^{(i+1)}(a-1))\leq\rk(a-1)\leq s$.) We terminate the procedure at $(U^{(j)},W^{(j)},\iota^{(j)},\pi^{(j)},r^{(j)},t^{(j)})$ when $r^{(j)}\leq\alpha_{k+1}s$, and set
\[(U_{k+1},W_{k+1},\iota_{k+1},\pi_{k+1},t_{k+1})=(U^{(j)},W^{(j)},\iota^{(j)},\pi^{(j)},t^{(j)}).\]
Note that the inequality $r^{(j)}\leq\alpha_{k+1}s$ automatically holds after any step falling into case (b). Therefore, since $r^{(i)}\leq r^{(0)}-\frac i3\alpha_{k+1}s$ for each $0\leq i\leq j$, we must have
\begin{equation}\label{eq:j-bound}
j+2\leq\frac{3\alpha_k}{\alpha_{k+1}}.
\end{equation}
Therefore, before the termination of the procedure, we may have done at most $\frac{3\alpha_k}{\alpha_{k+1}}-1$ case (a) steps, followed by at most one case (b) step. We conclude that
\[\dim U_k-\dim U_{k+1},\dim W_k-\dim W_{k+1}\leq jr^{(0)}+2s\leq\pr*{2+\frac{3\alpha_k^2}{\alpha_{k+1}}}s\leq 4s\]
and, using \eqref{eq:j-bound},
\[t^{(j)}\leq 2t^{(j-1)}+1=2(t^{(0)}+j-1)+1\leq 2t_k+\frac{6\alpha_k}{\alpha_{k+1}}-5=t_{k+1}.\]
Setting $t_{k+1}=t^{(j)}$, we have that at least $(4K)^{-t_{k+1}}\abs{A^2}$ elements $a$ of $A^2$ have
\[\rk(\pi_{k+1}(a-1)\iota_{k+1})\leq r^{(j)}\leq\alpha_{k+1}s,\]
as desired.
\end{proof}

\begin{remark} The iterative procedure described in \cref{lem:reg-iteration-technical} has some trade-off between the different parameters. For example, the $(4K)^{-30/\alpha^{-1/2}}$ term may be replaced by $(4K)^{-30/\alpha^{-\eps}}$ for any fixed $\eps>0$, at the cost of increasing the dimension decrement from $4ks$ to $Cks$ for some larger $C=C(\eps)$. (This may be done by replacing the double-exponential decay rate $\alpha_k\sim 2^{-2^{k/2}}$ with $\alpha_k\sim 2^{-\beta^k}$ for some $\beta\in(1,\sqrt2)$.) Regardless of the loss one is willing to incur in the exponent on $4K$, however, it is not possible to obtain a density decrement asymptotically less than $s\log\log(1/\alpha)$ using our methods.
\end{remark}

We are now in a position to prove \cref{lem:reg-el}.

\begin{proof}[Proof of \cref{lem:reg-el}] Set $V:=\FF^d$. Select $k=\floor{2\log_2\log_2(2s)}+1$, so that $\alpha:=2^{1-2^{k/2}}$ lies in the interval $[2^{1-\sqrt2}s^{-\sqrt2},s^{-1})$. Apply \cref{lem:reg-iteration-technical} to $A\subset\GL(V)$. We find vector spaces $U_k$ and $W_k$ of dimension at least $(\dim V)-4ks$, an injection $\iota_k\colon U_k\hookrightarrow V$, and a surjection $\pi_k\colon V\twoheadrightarrow W_k$ such that, for at least $(4K)^{-30/\alpha^{1/2}}\abs{A^2}$ elements $a$ of $A^2$, one has $\rk(\pi_k(a-1)\iota_k)\leq \alpha s$. Let $S$ be the set of such $a\in A^2$. Since $\alpha s<1$, we have $\pi_k(a-1)\iota_k=0$ for all $a\in S$. We also have
\[\abs{S}\geq (4K)^{-30\alpha^{-1/2}}\abs{A^2}\geq(2K)^{-60s}\abs{A^2},\]
where the latter inequality holds since $s\geq 6$. Let $U:=\im\iota_k$ and $W:=(\ker\pi_k)^\perp$ be of dimension at most
\[4ks\leq 8s\log_2\log_2(2s)+4s\leq 30s\log\log s,\]
where the latter inequality holds since $s\geq 6$. Since $\pi_k(a-1)\iota_k=0$ for all $a\in S$, we have $(a-1)W^\perp\subset U$ for all such $a$. This concludes the proof.
\end{proof}

\section{The subspace-to-subgroup lemma}\label{sec:subspace-to-subgroup}

The main object of this section is to prove \cref{lem:subspace-to-subgroup-flag}, as well as an analogue which will enable us to deduce \cref{prop:reg-el} from \cref{lem:reg-el}. This will complete the proofs of \cref{thm:GL,prop:reg-el}.

We begin by setting some notation which we carry through this section. We call the following ``setting $(\star)$.'' Let $d$ be a positive integer and let $\cF=(V_0,V_1,V_2,\ldots,V_m)$ be a flag on $\FF^d$ of parabolic dimension $D$. Let $A$ be a $K$-approximate subgroup of $H_\cF$, and let $\bV\subset\bP_\cF$ be a proper affine subspace which intersects $A^2$. Let $\eta:=\abs{A^2\cap\bV}/\abs{A^2}$. Write $\pi_i$ for the ring homomorphism $\bP_\cF\to\End(V_i/V_{i-1})$.

The backbone of this section is the following single step, which we will iterate.

\begin{lemma}\label{lem:dimension-decrement} Work in setting $(\star)$, and let $k$ be any positive integer. There exists some positive integer $j\leq k$ and some affine subspace $\bW$ of $\bP_\cF$ satisfying $\dim\bV-\dim\bW\geq j-1$ and
\[\pr*{\frac{\abs{A^2\cap\bW}}{\abs{A^2}}}\pr*{\frac{\abs{A^2\cap \Lmul(\bW\to \bV)}}{\abs{A^2}}}^{k-j}\geq\pr*{\frac{\eta}{2K}}^k.\]
\end{lemma}
In the conclusion of this lemma, we have $\dim\bV-\dim\bW\geq j-1$ rather than $\dim_\cF\bV-\dim_\cF\bW\geq j-1$. This is simply because $\dim$ is slightly easier to control than $\dim_\cF$ (and in the model case where $\cF$ is the trivial flag, they are nearly the same). In our iterative procedure, this condition is only needed to ensure termination of the iteration, and the precise quantitative properties of this termination are not so important.

\begin{proof}
We begin by showing that
\begin{equation}\label{eq:translate-int}
\sum_{a_1,\ldots,a_k\in A^2}\abs*{A^2\cap\bigcap_{i=1}^ka_i^{-1}\bV}\geq\pr*{\frac\eta K}^k\abs{A^2}^{k+1}.
\end{equation}
Indeed, we first note that $\abs{A^2x\cap A^2}\geq \abs{A}$ for each $x\in A^2$: if $x\in Aa$ for some $a\in A$, then $Aa\subset A^2x\cap A^2$ since $A$ is centered. This enables us to compute
\begin{align*}
\sum_{a_1,\ldots,a_k\in A^2}\abs*{A^2\cap\bigcap_{i=1}^ka_i^{-1}\bV}
&=\sum_{b\in A^2}\sum_{a_1,\ldots,a_k\in A^2}\prod_{i=1}^k1_{a_i^{-1}\bV}(b)\\
&=\sum_{b\in A^2}\pr*{\sum_{a\in A^2}1_{a^{-1}\bV}(b)}^k\\
&=\sum_{b\in A^2}\abs{A^2b\cap \bV}^k\\
&\geq\abs{A^2}\pr*{\frac1{\abs{A^2}}\sum_{b\in A^2}\abs{A^2b\cap \bV}}^k\\
&\geq\abs{A^2}\pr*{\frac1{\abs{A^2}}\sum_{b\in A^2}\sum_{x\in A^2\cap \bV}1_{A^2b}(x)}^k\\
&=\abs{A^2}\pr*{\frac1{\abs{A^2}}\sum_{x\in A^2\cap \bV}\abs{A^2x\cap A^2}}^k\\
&\geq\abs{A^2}\pr*{\frac1{\abs{A^2}}\sum_{x\in A^2\cap \bV}\abs{A}}^k=\frac{\abs{A^2\cap \bV}^k\abs{A}^k}{\abs{A^2}^{k-1}}\geq\pr*{\frac\eta K}^k\abs{A^2}^{k+1}.
\end{align*}
We now split up the sum in \eqref{eq:translate-int}. For each tuple $\mathbf a=(a_1,\ldots,a_k)$ of elements of $A^2$, let $S_{\mathbf a}\subset\set{1,\ldots,k}$ be a set of minimal size for which
\[\bigcap_{i=1}^{k}a_i^{-1}\bV=\bigcap_{i\in S_{\mathbf a}}a_i^{-1}\bV.\]
By \eqref{eq:translate-int}, there exists some set $S_0\subset\set{1,\ldots,k}$ for which
\[\sum_{\substack{\mathbf a:S_{\mathbf a}=S_0}}\abs*{A^2\cap\bigcap_{i=1}^ka_i^{-1}\bV}\geq\pr*{\frac\eta{2K}}^k\abs{A^2}^{k+1}.\]
Without loss of generality, we may assume $S_0=\set{1,\ldots,j}$ for some nonnegative integer $j\in\set{0,1,\ldots,k}$. We may then write
\[\sum_{a_1,\ldots,a_j\in A^2}\abs*{A^2\cap\bigcap_{i=1}^ja_i^{-1}\bV}\cdot\abs*{\set*{a\in A^2:a^{-1}\bV\supset \bigcap_{i=1}^ja_i^{-1}\bV}}^{k-j}\geq\pr*{\frac{\eta}{2K}}^k\abs{A^2}^{k+1}.\]
This implies that, for some affine subspace $\bW:=a_1^{-1}\bV\cap\cdots\cap a_j^{-1}\bV$ of $\bP_{\cF}$,
\[\pr*{\frac{\abs{A^2\cap \bW}}{\abs{A^2}}}\pr*{\frac{\abs[\big]{\set{a\in A^2:a\bW\subset \bV}}}{\abs{A^2}}}^{k-j}\geq\pr*{\frac{\eta}{2K}}^k.\]
If $j=0$, then $\bW=\bP_\cF$, and so no $a\in A^2$ can satisfy $a\bW\subset \bV$ (since $\bV$ is a proper affine subspace of $\bP_\cF$). We conclude that $j>0$. Moreover, by our choice of $S_0=\set{1,\ldots,j}$, the intersections $\bW^{(i)}:=a_1^{-1}\bV\cap\cdots\cap a_i^{-1}\bV$ for $1\leq i\leq j$ must each be distinct --- if $\bW^{(i-1)}=\bW^{(i)}$, then the term $a_i^{-1}$ can be removed from the intersection defining $\bW$, contradicting the definition of $S_{\mathbf a}$. As a result, since $\dim\bW^{(1)}=\dim\bV$ and all the inclusions in the chain
\[\bW^{(1)}\supset \bW^{(2)}\supset\cdots\supset \bW^{(j)}=\bW\]
are proper, we have $\dim\bV-\dim\bW\geq j-1$, as desired.
\end{proof}

Our first application of \cref{lem:dimension-decrement} is, in the case where $\bV=\vtr^{-1}(\mathbf t)$, to reduce from $\bV$ to a subspace of $\cF$-dimension $D/2$.

\begin{lemma}\label{lem:dimension-halving} Work in setting $(\star)$, and suppose that $\bV\subset\vtr^{-1}(\mathbf t)$ for some $\mathbf t\in\FF^m$. Then there exists some affine subspace $\bV'\subset\bP_\cF$ of $\cF$-dimension at most $D/2$ which is contained in some bitranslate of $\bV$ and satisfies
\[\frac{\abs{A^2\cap\bV'}}{\abs{A^2}}\geq\pr*{\frac{\eta}{2K}}^{d^2}.\]
\end{lemma}
\begin{proof} Apply \cref{lem:dimension-decrement} to $A$ and $\bV$ with $k=d^2$. We obtain some positive integer $j\leq d^2$ and some affine subspace $\bW$ of $\Mat_d(\FF)$ of dimension at most $\dim\bV-(j-1)$ satisfying
\begin{equation}\label{eq:dim-halfway-1}
\pr*{\frac{\abs{A^2\cap\bW}}{\abs{A^2}}}\pr*{\frac{\abs{A^2\cap \Lmul(\bW\to \bV)}}{\abs{A^2}}}^{d^2-j}\geq\pr*{\frac{\eta}{2K}}^{d^2}.
\end{equation}
Write $\bX:=\Lmul(\bW\to\bV)$. We first note that both $\bW$ and $\bX$ are contained in some bitranslate of $\bV$ by \cref{lem:bitranslate}. Now, if $j=d^2$, then $\dim\bW=0$, and we can conclude the result from \eqref{eq:dim-halfway-1} with $\bV'=\bW$. Otherwise, \eqref{eq:dim-halfway-1} leads to the weaker conclusion
\begin{equation}\label{eq:dim-halfway-2}
\min\pr*{\frac{\abs{A^2\cap\bW}}{\abs{A^2}},\frac{\abs{A^2\cap\bX}}{\abs{A^2}}}\geq\pr*{\frac{\eta}{2K}}^{d^2}.
\end{equation}

Now, since $\vtr(xy)=\mathbf t$ for every $x\in\bX$ and $y\in\bW$, we have $\tr(\pi_i(x)\pi_i(y))=t_i$ for each $1\leq i\leq m$. By \cref{lem:frobenius-inner-product}, we thus have $\dim\spn\pi_i(\bX)+\dim\spn\pi_i(\bW)\leq d_i^2+1$. This gives us that
\begin{align*}
\dim_\cF\bW+\dim_\cF\bX
&=\sum_{i=1}^m(\dim(\spn\pi_i(\bW))+\dim(\spn\pi_i(\bX))-2)\\
&\leq\sum_{i=1}^m(d_i^2-1)=D.
\end{align*}
So, one of $\set{\bW,\bX}$ has $\cF$-dimension at most $D/2$. The result follows from \eqref{eq:dim-halfway-2} by taking $\bV'$ to be whichever of $\set{\bW,\bX}$ has smaller $\cF$-dimension.
\end{proof}

Secondly, we iterate \cref{lem:dimension-decrement} to obtain that $A^2$ has large intersection with a set which is more structured than a general affine subspace.

\begin{lemma}\label{lem:subspace-to-subgroup} In setting $(\star)$, there exists some subalgebra $\bU\subset\bP_\cF$ for which $1+\bU$ is contained in some bitranslate of $\bV$ and
\[\frac{\abs{A^2\cap(1+\bU)}}{\abs{A^2}}\geq\frac{\eta^{d^4}}{(2K)^{d^6+3}}.\]
\end{lemma}
\begin{proof} Set $\bV_0:=\bV$ and $k_0:=d^4$. For each $i\geq 0$, given $\bV_i$ and $k_i$, use \cref{lem:dimension-decrement} to find some positive integer $j=:j_i\leq k_i$ and some affine subspace $\bW_i$ of $\Mat_d(\FF)$ satisfying $\dim\bV_i-\dim\bW_i\geq j_i-1$ and
\begin{equation}\label{eq:iter-fancy}
\pr*{\frac{\abs{A^2\cap\bW_i}}{\abs{A^2}}}\pr*{\frac{\abs{A^2\cap \Lmul(\bW_i\to\bV_i)}}{\abs{A^2}}}^{k_i-j_i}\geq\pr*{\frac{\eta}{2K}}^{k_i}.
\end{equation}
Define $k_{i+1}:=k_i-j_i$ and $\bV_{i+1}:=\Lmul(\bW_i\to\bV_i)$. Terminate the procedure if $\dim\bV_{i+1}=\dim\bV_i$, and let $\ell$ be the index at which this occurs, so that $\dim\bV_\ell=\dim\bV_{\ell+1}$. We note the following properties:
\begin{enumerate}
    \item By \cref{lem:lmul-dim}, $\dim\bV_{i+1}\leq\dim\bV_i$ for each $i$. This implies that $\dim\bV_i\leq\dim\bV_0-i$ for each $i\leq\ell$. In particular, $\ell\leq\dim\bV_0<d^2$.

    \item Since $\dim\bV_i-\dim\bW_i\geq j_i-1$, we must have $j_i\leq\dim\bV_i+1\leq\dim\bV_0+1$ for each $i$. Together with (1) and the choice $k_0=d^4$, this ensures that $k_i$ remains positive throughout the iteration.

    \item By \cref{lem:bitranslate}, each $\bV_{i+1}$ is contained in a right translate of $\bV_i$. So $\bV_\ell$ is contained in a right translate of $\bV$.

    \item Again by \cref{lem:lmul-dim}, since
    \[\dim\bV_{\ell+1}=\dim\Lmul(\bW_\ell\to\bV_\ell)=\dim\bV_\ell\]
    and both $\bW_\ell$ and $\bV_\ell$ contain invertible elements, there exists some subalgebra $\bU$ of $\Mat_d(\FF)$ and some $b\in\GL_d(\FF)$ for which (i) $(1+\bU)\cap\GL_d(\FF)$ is a subgroup of $\GL_d(\FF)$, (ii) $\bW_\ell\subset b(1+\bU)$, and (iii) $1+\bU$ is contained in a left translate of $\bV_\ell$.
\end{enumerate}

Let $H:=(1+\bU)\cap\GL_d(\FF)$, which is a subgroup of $\GL_d(\FF)$ by (4). Combining (3) and (4) gives that $1+\bU$ is contained in some bitranslate of $\bV$. Let $\eta_i:=\abs{A^2\cap\bV_i}/\abs{A^2}$. By \eqref{eq:iter-fancy}, we have
\[\eta_{i+1}^{k_{i+1}}\geq\pr*{\frac{\eta_i}{2K}}^{k_i}\]
for each $i$. This implies that, using \eqref{eq:iter-fancy},
\[\frac{\abs{A^2\cap\bW_\ell}}{\abs{A^2}}\geq\pr*{\frac{\eta_\ell}{2K}}^{k_\ell}\geq\eta_0^{k_0}\prod_{i=0}^{\ell}\pr*{\frac1{2K}}^{k_i}\geq\pr*{\frac{\eta_0}{(2K)^{\ell+1}}}^{k_0}\geq\pr*{\frac{\eta}{(2K)^{d^2}}}^{d^4}.\]
In particular, $\abs{A^2\cap bH}$ is nonempty. Let $a_0\in A^2\cap bH$ be an arbitrary element. Since $a_0^{-1}(A^2)\subset A^4$, we find, using \cref{lem:intersect-subgroup}, that
\[K^3\abs{A^2\cap H}\geq\abs{A^4\cap H}\geq\abs{a_0^{-1}A^2\cap a_0^{-1}bH}=\abs{A^2\cap bH}\geq\abs{A^2\cap\bW_\ell}.\]
We conclude
\[\frac{\abs{A^2\cap(1+\bU)}}{\abs{A^2}}=\frac{\abs{A^2\cap H}}{\abs{A^2}}\geq\frac{\abs{A^2\cap\bW_\ell}}{K^3\abs{A^2}}\geq\frac{\eta^{d^4}}{(2K)^{d^6+3}}.\qedhere\]
\end{proof}

The last step necessary in our proof of \cref{lem:subspace-to-subgroup-flag} is to move from the subalgebra $\bU$ to some $\bP_\cG$ with $\cG$ a refinement of $\cF$. To do this, we use Burnside's theorem \cite{Burnside}, which states that every proper subalgebra of a matrix algebra over an algebraically closed field has some nontrivial invariant subspace. (We need to apply Burnside's theorem in positive characteristic; a good reference for this is \cite{LR}.)

\begin{lemma}\label{lem:subalg-to-flag} Let $\rho\in(0,1)$. For every subalgebra $\bU\subset\bP_\cF$ which satisfies $\dim_\cF\bU\leq\rho\pdim(\cF)$, there exists a refinement $\cG$ of $\cF$ for which $\spn\bU\subset\bP_\cG$ and $\pdim(\cG)\leq\rho^{1/2}\pdim(\cF)$.
\end{lemma}
\begin{proof} Write $\overline\bU=\spn\bU$. For each $1\leq i\leq m$, let $\cF_i=(U_{i,0},\ldots,U_{i,m_i})$ be a maximal flag stabilized by $\pi_i(\overline\bU)$. For each $1\leq j\leq m_i$, we have a map $\phi_{i,j}$ from $\pi_i(\overline\bU)$ to $\End(U_{i,j}/U_{i,j-1})$. The image of $\phi_{i,j}$ is a subalgebra of $\End(U_{i,j}/U_{i,j-1})$ which contains the identity. So, Burnside's theorem implies that $\im\phi_{i,j}$ either is all of $\End(U_{i,j}/U_{i,j-1})$ or possesses some invariant subspace. The existence of such an invariant subspace, however, contradicts the maximality of the flag $\cF_i$. So, we conclude that
\[\dim\pi_i(\overline\bU)\geq\max_{1\leq j\leq m_i}\dim(U_{i,j}/U_{i,j-1})^2.\]
Write $d_{i,j}:=\dim U_{i,j}/U_{i,j-1}$, and let $\Delta_i:=\dim\pi_i(\overline\bU)$, so that $d_{i,j}\leq\sqrt{\Delta_i}\leq d_i$ for each $1\leq j\leq m_i$. We have
\begin{align}
\sum_{j=1}^{m_i}(d_{i,j}^2-1)
\notag&\leq \sum_{j=1}^{m_i}\pr*{d_{i,j}\sqrt{\Delta_i}-\frac{d_{i,j}}{\sqrt{\Delta_i}}}\\
\notag&=d_i\cdot\frac{\Delta_i-1}{\sqrt{\Delta_i}}\\
\notag&=\pr*{(\Delta_i-1)d_i^2\pr*{1-\frac1{\Delta_i}}}^{1/2}\\
\label{eq:subflag-bound}&\leq\pr*{(\Delta_i-1)d_i^2\pr*{1-\frac1{d_i^2}}}^{1/2}=\sqrt{d_i^2-1}\cdot\sqrt{\Delta_i-1}.
\end{align}
Let $\cG$ be the refinement of $\cF$ attained by refining each $V_i/V_{i-1}$ by inserting the terms in $\cF_i$, so that $\overline\bU\subset\bP_\cG$ by the definition of $\cF_i$. Using \eqref{eq:subflag-bound} and the Cauchy--Schwarz inequality, we may bound
\begin{align*}
\pdim\cG
&=\sum_{i=1}^m\sum_{j=1}^{m_i}(d_{i,j}^2-1)\\
&\leq\sum_{i=1}^m\sqrt{d_i^2-1}\sqrt{(\dim\pi_i(\overline\bU))-1}\\
&\leq\pr*{\sum_{i=1}^m(d_i^2-1)}^{1/2}\pr*{\sum_{i=1}^m\pr*{\dim\pi_i(\overline\bU)-1}}^{1/2}\\
&=(\pdim\cF)^{1/2}(\dim_\cF\bU)^{1/2}\leq \rho^{1/2}\pdim\cF.\qedhere
\end{align*}
\end{proof}

We now have everything we need to prove \cref{lem:subspace-to-subgroup-flag}.

\begin{proof}[Proof of \cref{lem:subspace-to-subgroup-flag}] Recall the statement of \cref{lem:subspace-to-subgroup-flag}: we are in the setting of $(\star)$ with $\bV=\vtr_\cF^{-1}(\mathbf t)$ for some $\mathbf t\in\FF^m$. By \cref{lem:dimension-halving}, we can find some affine subspace $\bV'\subset\bP_\cF$ of $\cF$-dimension at most $D/2$ which satisfies
\[\frac{\abs{A^2\cap\bV'}}{\abs{A^2}}\geq\pr*{\frac{\eta}{2K}}^{d^2}.\]
We are now again in the setting of $(\star)$, with $\eta$ replaced by $(\eta/(2K))^{d^2}$ and $\bV$ replaced by the smaller $\bV'$. So, we may apply \cref{lem:subspace-to-subgroup} to find some subalgebra $\bU\subset\bP_\cF$ for which $1+\bU$ is contained in some bitranslate of $\bV'$ and satisfies
\[\frac{\abs{A^2\cap(1+\bU)}}{\abs{A^2}}\geq\frac1{(2K)^{d^6+3}}\pr*{\frac{\eta}{(2K)^{d^2}}}^{d^4}=\frac{\eta^{d^6}}{(2K)^{2d^6+3}}.\]
Since $1+\bU$ is contained in a bitranslate of $\bV'$, we have $\dim_\cF\bU\leq\dim_\cF\bV'\leq D/2$. By \cref{lem:subalg-to-flag}, we can find some refinement $\cG$ of $\cF$ for which $\spn\bU\subset\bP_\cG$ and $\pdim(\cG)\leq\frac1{\sqrt2}D$. Since $\bU\subset\bP_\cG$, we also have $1+\bU\subset\bP_\cG$. In particular, every invertible element of $1+\bU$ is contained in $H_\cG$. We conclude
\[\frac{\abs{A^2\cap H_\cG}}{\abs{A^2}}\geq\frac{\abs{A^2\cap(1+\bU)}}{\abs{A^2}}\geq\pr*{\frac{\eta}{2K}}^{5d^6},\]
as desired.
\end{proof}

Lastly, we must deduce \cref{prop:reg-el}. The only missing ingredient is the following lemma, which places restrictions on subalgebras of $\FF^d$ which map every element of some vector space into some other fixed vector space.

\begin{lemma}\label{lem:small-rank-subalg} Let $U$ and $W$ be subspaces of $\FF^d$, and let $\bX$ be a subalgebra of $\Mat_d(\FF)$ satisfying $\bX W^\perp\subset U$. Then there exist subspaces $V_1\subset V_2\subset\FF^d$ satisfying $V_1\subseteq U$ and $W^\perp\subseteq V_2$ for which $\bX V_2\subset V_1$.
\end{lemma}
\begin{proof} Set $V_2:=\set{v:\bX v\subset U}$ and $V_1:=U\cap V_2$. From the fact that $\bX W^\perp\subset U$, we see that $V_1\subseteq U$ and $W^\perp\subseteq V_2$. Moreover, by the definition of $V_2$, we have $\bX V_2\subset U$. Finally, since $\bX$ is closed under multiplication, we have
\[\bX(\bX V_2)\subset \bX V_2\subset U.\]
So, every $v\in\bX V_2$ satisfies $\bX v\subset U$, and thus $v\in V_2$. This means $\bX V_2\subset U\cap V_2=V_1$, as desired.
\end{proof}

\begin{proof}[Proof of \cref{prop:reg-el}] Throughout what follows, we write $r:=\frac d{600\log\log d}$ and $s:=2r$. Note that $r<d/12$.

Let $A\subset\GL_d(\FF)$ be a $K$-approximate subgroup, and suppose that we are not in the ``not-too-irregular element'' case, i.e.~that every element $a\in A^2$ satisfy $\eigm(a)\leq d-r$. Our first step is to construct an approximate group $B$ whose structure controls that of $A$ and which satisfies the preconditions of \cref{lem:reg-el}.

We begin by noting that, for each $k\geq 1$, every element of $A^k$ has an eigenspace of codimension at most $rk$. Indeed, if $a_1,\ldots,a_k\in A$ have eigenspaces $V_1,\ldots,V_k$, then $V_1\cap\cdots\cap V_k$ is an eigenspace of $a_1\cdots a_k$ of codimension at most $\codim V_1+\cdots+\codim V_k$. In particular, since $r<d/12$, every element of $A^6$ has an eigenspace of codimension strictly less than $d/2$. For $a\in A^6$, let $\lambda_a$ be the eigenvalue associated to this (unique) eigenspace.

Define a map $\phi\colon A^6\to \GL_d(\FF)$ by $a\mapsto \lambda_a^{-1}a$. We claim that, if $k\ell\leq 6$, then $\phi|_{A^k}$ is a Freiman $\ell$-homomorphism. Indeed, if $a_1,\ldots,a_\ell,b_1,\ldots,b_\ell\in A^k$ satisfy $a_1\cdots a_\ell=b_1\cdots b_\ell$, then
\[\lambda_{a_1}\cdots\lambda_{a_\ell}=\lambda_{a_1\cdots a_\ell}=\lambda_{b_1\cdots b_\ell}=\lambda_{b_1}\cdots\lambda_{b_\ell},\]
and the result follows. Applied to $k=3$ and $\ell=2$, we obtain from \cite[Lemma~2.7.4]{Tointon} that $B:=\phi(A)$ is a $K$-approximate group. Moreover, for every $b\in B$, the $1$-eigenspace of $B$ has codimension at most $r$. This implies that the $1$-eigenspace of any element of $B^2$ has codimension at most $s=2r$. Equivalently, we have $\rk(b-1)\leq s$ for every $b\in B^2$.

We now apply \cref{lem:reg-el} to $B$, which we may do since $s\geq 6$ (as $d\geq 8000$). We obtain subspaces $U$ and $W$ of $\FF^d$ satisfying
\[\dim U,\dim W\leq 30s\log\log s<30s\log\log d=\frac d{10}\]
which satisfy
\[\abs{\set{b\in B^2:(b-1)W^\perp\subset U}}\geq\frac1{(2K)^{60s}}\abs{B^2}.\]
Let $\bV:=\set{x\in\Mat_d(\FF):(x-1)W^\perp\subset U}$. This is an affine subspace of $\Mat_d(\FF)$, and since $\dim U,\dim W<d/10$ it is properly contained in $\Mat_d(\FF)$. We are thus in setting $(\star)$ with $\cF$ the trivial flag $(\set{0},\FF^d)$, approximate group $B$, and $\eta:=\abs{B^2\cap\bV}/\abs{B^2}\geq (2K)^{-60s}$. So, we may apply \cref{lem:subspace-to-subgroup} to find some subalgebra $\bX\subset\bP_\cF$ for which $1+\bX$ is contained in some bitranslate $c_1\bV c_2$ of $\bV$ and for which
\begin{equation}\label{eq:B-int}
\frac{\abs{B^2\cap(1+\bX)}}{\abs{B^2}}\geq\frac{\eta^{d^4}}{(2K)^{d^6+3}}\geq\frac1{(2K)^{d^6+3+60sd^4}}\geq\frac1{(2K)^{2d^6-5}}.
\end{equation}
Let $\bV':=c_1(\bV-1)c_2$, so that $\bV'c_2^{-1}W^\perp\subset c_1U$; we have that $\bV'$ is a genuine subspace of $\Mat_d(\FF)$. Since $1\in 1+\bX\subset c_1\bV c_2$, we have $1-c_1c_2\in\bV'$. This implies that
\[\bX\subset\bV'-(1-c_1c_2)=\bV'.\]
Thus $\bX(c_2^{-1}W^\perp)\subset c_1U$. Applying \cref{lem:small-rank-subalg} to $\bX$, we conclude that there exist subspaces $V_1\subset V_2\subset\FF^d$ satisfying
\begin{align*}
\dim V_1&\leq\dim(c_1U)=\dim U<\frac d{10},\\
\dim V_2&\geq\dim(c_2^{-1}W^\perp)=\dim W^\perp>\frac{9d}{10},
\end{align*}
and $\bX V_2\subset V_1$. 

Now, let $H\leq\GL_d(\FF)$ be the set of $a\in\GL_d(\FF)$ for which there exists some scalar $\lambda$ for which $(a-\lambda)V_2\subset V_1$. This is a subgroup of $\GL_d(\FF)$ and contains every invertible element of $1+\bX$. By \cref{lem:intersect-subgroup-projection}, noting that $\phi|_{A^2}$ is a Freiman $3$-homomorphism and $\phi(A^2)=B^2$, we have
\[\frac{\abs{A^2\cap H}}{\abs{A^2}}\geq K^{-5}\frac{\abs{A^6\cap H}}{\abs{A^2}}\geq K^{-5}\frac{\abs{A^6\cap\phi^{-1}(1+\bX)}}{\abs{A^2}}\geq K^{-5}\frac{\abs{B^2\cap(1+\bX)}}{\abs{B^2}},\]
where we have used \cref{lem:intersect-subgroup} in the first step. Combining this with \eqref{eq:B-int} gives the result.
\end{proof}

\section{Application: Quasi-polynomial bounds in non-abelian Roth’s theorem}\label{sec:roth}

In this section, we prove \cref{thm:local-3-AP} on the growth of sets with no three-term arithmetic progressions. Some of our intermediate steps will also be useful for the proofs of \cref{thm:ramsey-cayley,thm:ramsey-cayley-general}, which we give in the following section. This section focuses primarily on the portions of the proofs which concern the interface between abelian and non-abelian groups. We will state all results we need from the abelian setting (here \cref{lem:3-AP-abelian,lem:corner-abelian}) and prove them in \cref{sec:abelian}.

We begin this section by providing a corollary of \cref{thm:solv-abelian-struct,thm:gen-abelian-struct} which will be useful for both of our applications. The first part gives us, if $A$ is as in the setting of \cref{thm:solv-abelian-struct,thm:gen-abelian-struct}, a large commuting set with small doubling contained in a translate of $A$. The second part is a similar, ``two-sided'' result which will be useful for understanding averaging $3$-APs.

\begin{lemma}\label{lem:abelian-struct-near}  Let $G$ be a group and $d\geq3$ an integer. Suppose either that (a) $G$ is solvable or (b) $G$ is finite and has no subquotient isomorphic to $\A_d$.
\begin{enumerate}[(1)]
    \item Let $A\subset G$ satisfy $\abs{AA^{-1}}\leq K\abs{A}$. Then there exists some abelian subgroup $H\leq G$, some element $g\in G$, and some set $A_0\subset gA\cap H$ which satisfies $\abs{A_0A_0^{-1}}\leq (2K)^{O(1)}\abs{A_0}$ and 
    \begin{equation}\label{eq:A0-size-nec}
    \abs{A_0}\geq\exp\pr*{\Omega\pr*{\frac{\log^{1/8}\abs{A}}{\exp(O(\log^2d\log\log d))\log 2K}}}.
    \end{equation}

    \item Let $A\subset G$ satisfy $\abs{A^2}\leq K\abs{A}$. Then there exists some abelian subgroup $H\leq G$, some $g\in G$, and some centered set $B_0\subset H$ which satisfies $\abs{B_0B_0^{-1}}\leq (2K)^{O(1)}\abs{B_0}$,
    \[\#\big\{(x,y)\in B_0\times B_0:xgy\in A\big\}\geq(2K)^{-O(1)}\abs{B_0}^2,\]
    and
    \begin{equation}\label{eq:B0-size-nec}
    \abs{B_0}\geq\exp\pr*{\Omega\pr*{\frac{\log^{1/8}\abs{A}}{\exp(O(\log^2d\log\log d))\log 2K}}}.
    \end{equation}
\end{enumerate}
\end{lemma}
\begin{proof} We first prove (1). Since $\abs{AA^{-1}}\leq K\abs{A}$, \cref{lem:nice-set-finder} gives us the existence of a centered set $B\subset A^{-1}A$ with $\abs{B}\geq\abs{A}/2K$ and
\[\abs{B^m}\leq\abs{AB^m}\leq\abs{AB^mA^{-1}}\leq2^mK^{2m+1}\abs{B}\]
for all positive integers $m$. In particular, $\abs{B^3}\leq (2K)^{O(1)}\abs{B}$. By \cref{thm:solv-abelian-struct,thm:gen-abelian-struct}, there exists some commuting subset $T\subset B^4$ satisfying
\begin{equation}\label{eq:T-size-bound-0}
\abs{T}\geq\exp\pr*{\Omega\pr*{\frac{\log^{1/8}\abs{B}}{\exp(O(\log^2d\log\log d))\log2K}}}.
\end{equation}
Let $H\leq G$ be the abelian subgroup of $G$ which $T$ generates; by replacing $T$ with a larger commuting set if necessary, we may assume $T=B^4\cap H$. Now,
\[\sum_{g\in TT^{-1}A^{-1}}\abs{gA\cap TT^{-1}}=\sum_{a\in A}\sum_{b\in TT^{-1}}1_{TT^{-1}A^{-1}}(ba^{-1})=\abs{A}\cdot\abs{TT^{-1}}.\]
Therefore there exists some $g\in TT^{-1}A^{-1}$ for which $A_0:=gA\cap TT^{-1}\subset gA\cap H$ satisfies
\begin{equation}\label{eq:A0-size}
\abs{A_0}\geq\frac{\abs{A}\cdot\abs{TT^{-1}}}{\abs{TT^{-1}A^{-1}}}=\frac{\abs{A}\cdot\abs{TT^{-1}}}{\abs{ATT^{-1}}}\geq\frac{\abs{A}\cdot\abs{TT^{-1}}}{\abs{AB^8}}\geq(2K)^{-O(1)}\abs{TT^{-1}}.
\end{equation}
Combining \eqref{eq:T-size-bound-0} and \eqref{eq:A0-size} gives us that $A_0$ is large enough to satisfy \eqref{eq:A0-size-nec}. We need only upper-bound $\abs{A_0A_0^{-1}}$. We have $\abs{B^{10}}\leq 2^{10}K^{21}\abs{B}\leq 2^{10}K^{21}\abs{B^2}$. Ruzsa's covering lemma thus implies that $B^4$ is a $(2K)^{O(1)}$-approximate group. Therefore, by \cref{lem:intersect-subgroup} and \eqref{eq:A0-size}, we have that
\[\abs{A_0A_0^{-1}}\leq\abs{(TT^{-1})(TT^{-1})^{-1}}\leq\abs{B^{16}\cap H}\leq(2K)^{O(1)}\abs{T}\leq(2K)^{O(1)}\abs{A_0}.\]
This concludes the proof of (1).

The proof of (2) is quite similar, but in place of \cref{lem:nice-set-finder} we use \cite[Theorem~2.5.6]{Tointon}. Given the stronger hypothesis of $\abs{A^2}\leq K\abs{A}$, we find a subset $B\subset A$ with $\abs{B}\geq\abs{A}/K$ for which $\abs{B^m}\leq K^{m-1}\abs{B}$ for all $m$. Moreover, the set $B':=(B\cup B^{-1}\cup\{1\})^2$ is a $(2K)^{O(1)}$-approximate group of size at most $7K^2\abs{A}$, at most $K^2$ left- or right-translates of which cover $A$. By \cref{thm:solv-abelian-struct,thm:gen-abelian-struct}, there exists some commuting subset $T\subset (B\cup B^{-1}\cup\{1\})^4\subset(B')^2$ satisfying the size bound \eqref{eq:T-size-bound-0}. As before, we let $H\leq G$ be the abelian subgroup of $G$ which $T$ generates. Let $B_0:=(B')^2\cap H$. Since $B_0\supset T$, \eqref{eq:T-size-bound-0} is enough to imply that $B_0$ is large enough to satisfy \eqref{eq:B0-size-nec}. Now, \cref{lem:intersect-subgroup} implies that $B_0$ is a $(2K)^{O(1)}$-approximate group, and so we have $\abs{B_0B_0^{-1}}\leq(2K)^{O(1)}\abs{B_0}$. Finally, we have (using that $B_0$ is centered) that
\[\sum_{g\in B_0BB_0}\#\{(x,y)\in B_0\times B_0:xgy\in A\}\geq\sum_{g\in B_0BB_0}\#\{(x,y)\in B_0\times B_0:xgy\in B\}=\abs{B}\cdot\abs{B_0}^2,\]
and so there exists some $g\in G$ for which
\begin{align*}
\#\{(x,y)\in B_0\times B_0:xgy\in A\}
&\geq\frac{\abs{B}\cdot\abs{B_0}^2}{\abs{B_0BB_0}}\\
&\geq\frac{\abs{B}\cdot\abs{B_0}^2}{\abs{(B')^5}}\geq\frac{\abs{B}\cdot\abs{B_0}^2}{(2K)^{O(1)}\abs{B'}}\geq(2K)^{-O(1)}\abs{B_0}^2.\qedhere
\end{align*}
\end{proof}

\subsection{Three-term arithmetic progressions}\label{sec:3-AP}

We now turn our attention to sets without three-term arithmetic progressions. We begin by stating the following two results in the abelian setting, which we will prove in \cref{sec:abelian}. Their proofs result from combining the main result of \cite{Sanders} --- quasi-polynomial bounds in (essentially) the Freiman--Ruzsa theorem --- with the recent advances on Roth's theorem by Kelley and Meka \cite{KelleyMeka} and on the corners theorem by Jaber, Liu, Lovett, Ostuni, and Sawhney \cite{JLLOS}. In the following two statements, we use additive notation for abelian groups. Elsewhere in this section we exclusively use multiplicative notation, even for sets which are contained in abelian groups.

\begin{lemma}[Abelian local Roth theorem]\label{lem:3-AP-abelian}
There exists an absolute constant $c>0$ for which the following holds. Let $G_0$ be an abelian group. Let $A\subset G_0$ be a subset for which, if $x,x+y,x+2y\in A$, then $y=0$.\footnote{As noted in the statement of \cref{thm:local-3-AP}, if $G_0$ has $2$-torsion, we allow $2y=0$.} Then
\[\frac{\abs{A-A}}{\abs{A}}\geq\exp\pr*{c\log^{1/144}\abs{A}}.\]
\end{lemma}

\begin{lemma}[Abelian local corners theorem]\label{lem:corner-abelian}
There exists an absolute constant $c>0$ for which the following holds. Let $G_0$ be an abelian group. Let $B\subset G_0$ be an arbitrary subset, and suppose that $A\subset B\times B$ is such that, if $(x,y),(x+d,y),(x,y+d)\in A$, then $d=0$. Then
\[\frac{\abs{B-B}^2}{\abs{A}}\geq\exp\pr*{c\log^{1/7200}\abs{B}}.\]
\end{lemma}

We now prove \cref{thm:local-3-AP} using \cref{lem:abelian-struct-near,lem:3-AP-abelian,lem:corner-abelian}. The first part of \cref{lem:abelian-struct-near}, along with the Roth-type \cref{lem:3-AP-abelian}, will be used to prove the translational version, while the second part of \cref{lem:abelian-struct-near}, along with the corners-type \cref{lem:corner-abelian}, will be used to prove the averaging version.

\begin{proof}[Proof of \cref{thm:local-3-AP}] Recall our setting: we have some group $G$ and some finite subset $A\subset G$. The group $G$ is either (a) is solvable or (b) is finite and has no subgroup isomorphic to $(\ZZ/2\ZZ)^d$, where $d=\exp(c\sqrt{\log\log\abs{A}/\log\log\log\abs{A}})$ (and $c>0$ is a constant we are allowed to choose). The set $A$ avoids translational (resp.\ averaging) $3$-APs, and we wish to lower-bound $\abs{AA^{-1}}/\abs{A}$ (resp.\ $\abs{A^2}/\abs{A}$).

By \cref{lem:sectional-p-rank}, the condition in (b) implies that $G$ has no subquotient isomorphic to $(\ZZ/2\ZZ)^{2d^2}$. In particular, since the alternating group $\A_{d'}$ has a subgroup isomorphic to $(\ZZ/2\ZZ)^{\floor{d'/2}-1}$ and subgroups of subquotients of $G$ are subquotients of $G$, we conclude that $G$ has no subquotient isomorphic to $\A_{4d^2+3}$. With $c$ chosen small enough, our choice of $d$ is small enough to apply \cref{lem:abelian-struct-near}. 

We now specialize to translational $3$-APs. Suppose that $A$ contains no nontrivial translational $3$-APs and let $K:=\abs{AA^{-1}}/\abs{A}$. We find from \cref{lem:abelian-struct-near}(1) some abelian subgroup $H\leq G$, some $g\in G$, and some set $A_0\subset gA\cap H$ which satisfies $\abs{A_0A_0^{-1}}\leq (2K)^{O(1)}\abs{A_0}$ and
\begin{equation}\label{eq:A0-large-appl}
\abs{A_0}\geq\exp\left(\Omega\left(\frac{\log^{1/8}\abs{A}}{\exp(O(\log^2(4d^2+3)\log\log(4d^2+3)))\log 2K}\right)\right)\geq\exp\left(\Omega\left(\frac{\log^{1/9}\abs{A}}{\log 2K}\right)\right)
\end{equation}
(where the latter inequality follows from a sufficiently small choice of $c$). As $A_0$ is contained in a left translate of $A$, the property that $A_0$ avoids nontrivial translational $3$-APs implies the same of $A_0$. Moreover, $A_0\subset H$ is a subset of an abelian group. We thus apply \cref{lem:3-AP-abelian} to conclude
\begin{equation}\label{eq:A0-small-appl}
(2K)^{O(1)}\geq\frac{\abs{A_0A_0^{-1}}}{\abs{A_0}}\geq\exp\pr*{\Omega\pr*{\log^{1/144}\abs{A_0}}}.
\end{equation}
Combining \eqref{eq:A0-large-appl} and \eqref{eq:A0-small-appl} gives
\[\log 2K\gtrsim \log^{1/144}\abs{A_0}\gtrsim\frac{\log^{1/1296}\abs{A}}{\log^{1/144}2K}\implies\log 2K\gtrsim \log^{1/1305}\abs{A}.\]
This proves (1).

Now we treat averaging $3$-APs. Suppose that $A$ contains no nontrivial translational $3$-APs and let $K:=\abs{A^2}/\abs{A}$. We find from \cref{lem:abelian-struct-near}(2) some abelian subgroup $H\leq G$, some $g\in G$, and some centered set $B_0\subset H$ which satisfies $\abs{B_0B_0^{-1}}\leq (2K)^{O(1)}\abs{B_0}$ and
\begin{align}
(2K)^{-O(1)}\abs{B_0}^2&\leq\#\big\{(x,y)\in B_0\times B_0:xgy\in A\big\};\label{eq:B0-dense}\\
\abs{B_0}&\geq\exp\left(\Omega\left(\frac{\log^{1/9}\abs{A}}{\log 2K}\right)\right).\label{eq:B0-large-appl}
\end{align}
Let $S:=\{(x,y)\in B_0\times B_0:xgy^{-1}\in A\}$. Since $B_0$ is centered, \eqref{eq:B0-dense} implies that $\abs{S}/\abs{B_0}^2\geq(2K)^{-O(1)}$. On the other hand, $S$ avoids corners: if $(x,y),(xd,y),(x,yd)\in S$ then
\[(xgy^{-1})^2=xgy^{-1}d^{-1}dxgy^{-1}=(xg(yd)^{-1})((xd)gy^{-1})\]
gives a nontrivial averaging $3$-AP in $A$ as long as $d\neq 1$. Therefore, since $S\subset B_0\times B_0$ and $B_0$ lies in an abelian group $H$, \cref{lem:corner-abelian} implies
\[\exp\pr*{\Omega\pr*{\log^{1/7200}\abs{B_0}}}\leq\frac{\abs{B_0B_0^{-1}}^2}{\abs{S}}=\frac{\abs{B_0}^2}{\abs{S}}\left(\frac{\abs{B_0B_0^{-1}}}{\abs{B_0}}\right)^2\leq(2K)^{O(1)}.\]
Comparing this bound to \eqref{eq:B0-large-appl} gives
\[\log^{7200}2K\gtrsim\log\abs{B_0}\gtrsim\frac{\log^{1/9}\abs{A}}{\log 2K},\]
which is enough to give the result.
\end{proof}

\section{Application: Ramsey Cayley graphs}\label{sec:ramsey-cayley}

In this section, we prove \cref{thm:ramsey-cayley,thm:ramsey-cayley-general} using \cref{thm:gen-abelian-struct} and input from \cite{CFPY}. Like the previous section, we will need some inputs from (abelian) additive combinatorics whose proofs we defer to \cref{sec:abelian}.

We begin by outlining the proof in \cite{CFPY} of the existence of Ramsey Cayley graphs in the abelian setting. A major component of that work is the following counting result, which is proven using purely combinatorial arguments.

\begin{theorem}[{\cite[Theorem~1.2]{CFPY}}]\label{thm:CFPY-main} With $C=6000$, the following holds: if $G$ is a finite group of order $N$ and $n<N$ is a positive integer, the number of subsets $A\subset G$ of size $n$ with $\abs{AA^{-1}}\leq Kn$ is at most $N^{C(K+\log n)}(2^CK)^n$.\footnote{The value of the constant $C$ is not given in \cite{CFPY}, but a value of $6000$ can easily be extracted from their arguments and we include it here for concreteness. Also, in \cite[Theorem~1.2]{CFPY}, the $2^C$ is written as $C$, but in their result the constant multiplier on $K$ is approximately $2^{1000}$ so we bound it by $2^C$ to keep the value of $C$ manageable.}
\end{theorem}

\cref{thm:CFPY-main} relates to the Ramsey property for the following reason: a set $A\subset G$ is a clique (resp.\ independent set) in $\Cay(G,S)$ if $AA^{-1}\setminus\{1\}\subset S$ (resp.\ $AA^{-1}\setminus\{1\}\subset G\setminus S$). 

Consider taking a uniformly random symmetric subset $S$ of a group $G$ and constructing the Cayley graph $\Gamma:=\Cay(G,S)$. The probability that a set $A$ forms a clique or independent set in $\Gamma$ is at most $2^{-(\abs{AA^{-1}}-3)/2}$. Let $n:=\ceil{C_1\log N}$ for some large constant $C_1$. Taking a union bound over all $A\subset G$ of $\abs{A}=n$ and using \cref{thm:CFPY-main}, one gets that, with high probability, the graph $\Gamma$ has no clique or independent set $A$ of size $n$ satisfying $\abs{AA^{-1}}\geq C_2n\log n$ (for some constant $C_2$ depending on $C_1$). So, we need only concern ourselves with those sets $A$ satisfying $\abs{AA^{-1}}=O(\abs{A}\log\abs{A})$. On the other hand, for $G=(\ZZ/5\ZZ)^n$, say, \cref{thm:CFPY-main} is nearly tight even for $K=1$, and the union bound fails for small $K$: the random Cayley graph $\Gamma$ has clique number on the order of $\log\abs{G}\log\log\abs{G}$ rather than $\log\abs{G}$.

The way such issues are handled in \cite{CFPY} is to choose the random generating set $S$ non-uniformly. In particular, for $G$ abelian with $\gcd(\abs{G},6)=1$, Conlon, Fox, Pham, and Yepremyan construct a distribution (see \cite[Section~4.2]{CFPY}) on subsets $S$ of $G$ which has enough independence for the arguments in the previous paragraph to go through, but which deterministically has the ``subgroup-avoiding'' property that, if $x\neq 1$, then neither $\set{x,x^2,x^4,x^8}\cap S$ nor $\set{x,x^2,x^4,x^8}\cap(G\setminus S)$ is empty. This ensures that no set $A$ for which $AA^{-1}$ contains a nontrivial instance of the pattern $\set{x,x^2,x^4,x^8}$ may be a clique or independent set in the resulting random Cayley graph. An argument using Pl\"unnecke's inequality \cite[Lemma~4.5]{CFPY} implies that such a set must satisfy $\abs{AA^{-1}}\geq \abs{A}^{1+1/144}$, and so the union bound using \cref{thm:CFPY-main} is strong enough to reach the desired conclusion.

This construction can be carried through for non-abelian groups $G$. However, the argument using Pl\"unnecke's inequality no longer applies. We use \cref{thm:gen-abelian-struct} to circumvent this issue. If $A\subset G$ is such that $\abs{AA^{-1}}$ is small, then \cref{thm:gen-abelian-struct} furnishes (see \cref{lem:abelian-struct-near}) a large commuting set contained in a translate of $A$. We will then be able to use the translation-invariance properties of the map $A\mapsto AA^{-1}$ to prove the following.

\begin{lemma}\label{lem:cliques-large-doubling} For each positive integer $r$, there exists some $\eta_r\geq\exp(-O(\log^2r\log\log r))$ for which the following holds. Suppose that neither $2^r$ nor $3^r$ divide $\abs{G}$. Suppose $A\subset G$ satisfies the property that, if $\gcd(\ord(x),6)=1$ and $x,x^2,x^4,x^8\in AA^{-1}$, then $x=1$. Then
\[\abs{AA^{-1}}>\abs{A}\exp\pr*{\eta_r\log^{1/16}\abs{A}}.\]
\end{lemma}

We will prove \cref{lem:cliques-large-doubling} in \cref{sec:avoid-pattern}.

\subsection{Constructing Ramsey Cayley graphs} Now that we have outlined the approach, we state our randomized construction of Cayley graphs. We then prove \cref{thm:ramsey-cayley}, using \cref{thm:CFPY-main} (that is, \cite[Theorem~1.2]{CFPY}) and \cref{lem:cliques-large-doubling}. This construction is very similar to that used in \cite{CFPY}.

\begin{construction}\label{constr:ramsey-cayley} Let $G$ be a group. Pick a subset $S$ of $G$ randomly as follows. First, partition the elements of $G\setminus\set{1}$ into equivalence classes under the relation $g\sim g^{-1}$. Next, form a graph $\Pi$ between these equivalence classes where $x\sim x^{-1}$ is connected to $x^2\sim x^{-2}$ for each $x\in G\setminus\{1\}$ with odd order. The graph $\Pi$ is a union of cycles, self-loops (corresponding to elements of $G$ of order $3$), and disconnected vertices (corresponding to elements of $G$ of even order). Find a maximal matching $\Sigma$ in $\Pi$. Perform the following independently:
\begin{itemize}
    \item For each vertex $g\sim g^{-1}$ not incident to any edge of the matching $\Sigma$, place $\{g,g^{-1}\}$ in $S$ with probability $1/2$.

    \item For each edge $(g\sim g^{-1},g^2\sim g^{-2})$ of the matching $\Sigma$, place one of $\{g,g^{-1}\}$ or $\{g^2,g^{-2}\}$ in $S$, uniformly at random.
\end{itemize}
\end{construction}

We state the essential properties of this construction as the following two simple lemmas, which are very similar to \cite[Claims~4.3~and~4.4]{CFPY}.

\begin{lemma}\label{lem:independence} For any set $A\subset G$, the probability that $A$ is a clique or independent set in $\Gamma=\Cay(G,S)$ generated by \cref{constr:ramsey-cayley} is at most $2^{-(\abs{AA^{-1}}-3)/2}$.
\end{lemma}
\begin{proof} The set $A$ is a clique or independent set in $\Gamma$ if and only if either $AA^{-1}\setminus\set{1}\subset S$ or $AA^{-1}\setminus\set{1}\subset G\setminus S$. Since $S$ and $G\setminus S$ have the same distribution, it suffices to show that $AA^{-1}\setminus\set{1}\subset S$ with probability at most $2^{-(\abs{AA^{-1}}-1)/2}$.

Let $\Sigma$ be the matching constructed in \cref{constr:ramsey-cayley}. If $AA^{-1}$ contains any edges of $\Sigma$, then $AA^{-1}\setminus\{1\}$ can never lie in $S$. Otherwise, the equivalence classes $g\sim g^{-1}$ for $g\in AA^{-1}\setminus\{1\}$ each lie in $S$ independently with probability $1/2$. The result follows.
\end{proof}

\begin{lemma}\label{lem:constr-1248-avoid} Let $A\subset G$, and suppose that there exists some $x\in G\setminus\set{1}$ satisfying $\gcd(\ord(x),6)=1$ and $x,x^2,x^4,x^8\in AA^{-1}$. Then, deterministically, $A$ is neither a clique nor independent set in any Cayley graph $\Gamma=\Cay(G,S)$ generated by \cref{constr:ramsey-cayley}.
\end{lemma}
\begin{proof} The equivalence classes $\set{x\sim x^{-1},x^2\sim x^{-2},x^4\sim x^{-4},x^8\sim x^{-8}}$ (not all of which must be distinct) lie, in that order, in a cycle in the auxiliary graph $\Pi$ constructed in \cref{constr:ramsey-cayley}. Since $\gcd(\ord(x),6)\neq 1$, this cycle is not a self-loop. Now, since $\Sigma$ is a maximal matching in $\Pi$, it must contain at least one of the edges
\[\set[\big]{(x\sim x^{-1},x^2\sim x^{-2}),(x^2\sim x^{-2},x^4\sim x^{-4}),(x^4\sim x^{-4},x^8\sim x^{-8})}.\]
This implies that these four pairs of elements of $G$ can neither all lie in $S$ nor all lie in $G\setminus S$. Therefore there exist some $b_1,b_2\in AA^{-1}\setminus\{1\}$ for which $b_1\in S$ and $b_2\not\in S$. We conclude that $A$ is neither a clique nor an independent set in $\Gamma$.
\end{proof}

\cref{thm:ramsey-cayley} will now follow from the following general result.

\begin{lemma}\label{lem:ramsey-cayley-general} Let $C_0>0$ and $c_0\in(0,1)$ be constants. Let $G$ be a finite group and let $\mathcal D$ be a probability distribution on subsets of $G$ such that, if $S\sim\mathcal D$ and $\Gamma=\Cay(G,S)$, then the following two properties hold:
\begin{enumerate}[(i)]
    \item For any set $A\subset G$, the probability that $A$ is a clique or independent set in $\Gamma$ is at most $C_0\cdot 2^{-c_0\abs{AA^{-1}}}$.

    \item Any $A\subset G$ of size $n=\ceil{12000c_0^{-1}\log\abs{G}}$ with $\abs{AA^{-1}}\leq 20c_0^{-1}n\log n$, is, deterministically, neither a clique nor an independent set in $\Gamma$.
\end{enumerate}
Suppose that $N$ is sufficiently large in terms of $C_0$ and $c_0$. Then, with probability at least $1-1/N$, the graph $\Gamma$ is a $12000c_0^{-1}$-Ramsey Cayley graph.
\end{lemma}
\begin{proof} This lemma is essentially a corollary of \cref{thm:CFPY-main}, using arguments presented in \cite[Proof~of~Theorem~4.2]{CFPY}. We present it here for completeness.

Let $C:=6000$ be the constant given in \cref{thm:CFPY-main}. Write $C_1:=2Cc_0^{-1}$. Let $N:=\abs{G}$ and $n:=\ceil{C_1\log N}$. We can assume $n\geq 2^{4C}$. By condition (ii), any subset $A\subset G$ of size $n$ with $\abs{AA^{-1}}\leq 20c_0^{-1}n\log n$ is neither a clique nor an independent set in $\Gamma$.

It remains to show that those sets $A$ of size $n$ with $\abs{AA^{-1}}$ larger than this threshold are unlikely to be cliques or independent sets. Fix some $m\in [20c_0^{-1}n\log n,n^2]$, and let $K:=m/n$. Consider those sets with $\abs{A}=n$ and $\abs{AA^{-1}}=m$. By \cref{thm:CFPY-main}, there are at most $N^{C(K+\log n)}(CK)^n$ such sets. By (i), each such set is a clique or independent set in $\Gamma$ with probability no more than $C_0\cdot 2^{-c_0m}$. We conclude by the union bound that the probability that \emph{any} such set is a clique or independent set in $\Gamma$ is at most
\begin{align*}
C_02^{-c_0Kn}N^{C(K+\log n)}(2^CK)^n
&\leq C_0\left(2^{-c_0K}e^{C_1^{-1}C(K+\log n)}2^Cn\right)^n\\
&=C_0\exp\pr*{n\pr*{-c_0K\log 2+\frac{c_0}2(K+\log n)+C\log 2+\log n}}\\
&\leq C_0\exp\pr*{n\pr*{\frac52\log n-c_0K\left(\log 2-\frac12\right)}}\leq C_0\exp(-n\log n),
\end{align*}
since $K\geq 20c_0^{-1}\log n$. Taking a union bound over the possible values of $m$, we conclude that, with probability at least $1-C_0n^{-(n-2)}$, no set of size $n$ is a clique or independent set in $\Gamma$. This concludes the proof.    
\end{proof}

We are now ready to prove \cref{thm:ramsey-cayley}.

\begin{proof}[Proof of \cref{thm:ramsey-cayley}] Let $r$ be a positive integer and let $G$ be a group whose order is neither a multiple of $2^r$ nor $3^r$. We must show that, as long as $\abs{G}$ is sufficiently large in terms of $r$, the group $G$ has a $24000$-Ramsey Cayley graph.

We apply \cref{lem:ramsey-cayley-general} with $C_0=2^{3/2}$ and $c_0=1/2$ to the probability distribution on subsets of $G$ given by \cref{constr:ramsey-cayley}. Let $S$ be chosen from that distribution and let $\Gamma:=\Cay(G,S)$. Condition (i), with our choice of $C_0$ and $c_0$, becomes that a set $A$ is a clique or independent set in $\Gamma$ with probability at most $2^{(\abs{AA^{-1}}-3)/2}$; this is precisely \cref{lem:independence}. 

Take $n:=\lceil 24000\log\abs{G}\rceil$. To verify condition (ii) of \cref{lem:ramsey-cayley-general}, we must show that any set $A$ of size $n$ which can potentially (i.e.~with nonzero probability) be a clique or an independent set in $\Gamma$ must satisfy  $\abs{AA^{-1}}>40n\log n$. This follows from \cref{lem:cliques-large-doubling,lem:constr-1248-avoid} and some computation. Indeed, by \cref{lem:constr-1248-avoid}, any such set $A$ must avoid the pattern $\{x,x^2,x^4,x^8\}$ for $x\neq 1$ of order relatively prime to $6$. Then, by \cref{lem:cliques-large-doubling}, any such set must satisfy
\begin{equation}\label{eq:AA-large-in-r}
\abs{AA^{-1}}\geq n\exp\pr*{\eta_r\log^{1/16}n},
\end{equation}
for some $\eta_r\geq\exp(-O(\log^2r\log\log r))$. Recall our condition that the largest product of powers of $2$ and $3$ that divide $N$ is at most $\exp(\exp((\log\log\log N)^{1/3}))$. This implies the bound $r\leq2\exp((\log\log\log N)^{1/3})$, and so (for $N$ larger than some absolute constant)
\[\log^2r\log\log r\leq2(\log\log\log N)^{3/4}.\]
Therefore we have $\eta_r\geq\exp(-(\log\log\log N)^{3/4})>\log^{-1/32}n$ as long as $N$ is large enough. By \eqref{eq:AA-large-in-r}, we then have $\abs{AA^{-1}}\geq n\exp(\log^{1/32}n)>40n\log n$. This proves condition (ii) and concludes the proof.
\end{proof}

\subsection{Extending \texorpdfstring{\cref{thm:ramsey-cayley}}{Theorem 1.6}}

We now modify \cref{constr:ramsey-cayley} to prove \cref{thm:ramsey-cayley-general}, finding Ramsey Cayley graphs in groups without large elementary $2$- or $3$-subgroups and without (particular) large cyclic subgroups. See \cref{sec:conc-ramsey} for some discussion of why our method fails when any of these subgroups appear.

The construction we use to prove \cref{thm:ramsey-cayley-general} generalizes \cref{constr:ramsey-cayley}. That construction is not especially useful for groups of even order, since the ``subgroup-avoiding'' property \cref{lem:constr-1248-avoid} completely neglects $2$- and $3$-torsion. To this end, we must replace the pairing $\{x,x^2\}$ with some other pairing for elements of even order. A natural extension is to pair $x$ and $x^q$ for $q$ the smallest prime not dividing the order of $x$. This is essentially what we do:

\begin{construction}\label{constr:ramsey-cayley-general} Let $G$ be a group. Pick a subset $S$ of $G$ randomly as follows. First, partition the elements of $G\setminus\set{1}$ into equivalence classes under the relation $g\sim g^{-1}$. Next, form a graph $\Pi$ between these equivalence classes where each $x\sim x^{-1}$ is connected to $x^q\sim x^{-q}$ for $q$ the smallest prime not dividing the order of $x$. (This choice of $q$ is also the smallest prime not dividing the order of $x^q$.) The graph $\Pi$ is a union of cycles and self-loops. Find a maximal matching $\Sigma$ in $\Pi$. Perform the following independently:
\begin{itemize}
    \item For each vertex $g\sim g^{-1}$ not incident to any edge of the matching $\Sigma$, place $\{g,g^{-1}\}$ in $S$ with probability $1/2$.

    \item For each edge $(g\sim g^{-1},g^q\sim g^{-q})$ of the matching $\Sigma$, place one of $\{g,g^{-1}\}$ or $\{g^q,g^{-q}\}$ in $S$, uniformly at random.
\end{itemize}
\end{construction}

We now must prove analogues of \cref{lem:independence,lem:constr-1248-avoid} for \cref{constr:ramsey-cayley-general}. \cref{lem:independence} remains true verbatim in this setting; the analogue of \cref{lem:constr-1248-avoid} is rather more complex.

\begin{lemma}\label{lem:constr-general-properties} For any set $A\subset G$, the following properties hold.
\begin{enumerate}[(1)]
    \item The probability that $A$ is a clique or independent set in $\Gamma:=\Cay(G,S)$ generated by \cref{constr:ramsey-cayley-general} is at most $2^{-(\abs{AA^{-1}}-3)/2}$.    

    \item Suppose that there exists $x\in G$, not of order dividing $24$, such that, if $q$ is the smallest prime not dividing the order of $x$, then $x,x^q,x^{q^2},x^{q^3}\in AA^{-1}$. Then, deterministically, $A$ is neither a clique nor an independent set in any Cayley graph $\Gamma:=\Cay(G,S)$ generated by \cref{constr:ramsey-cayley-general}. 
\end{enumerate}
\end{lemma}

Before proving \cref{lem:constr-general-properties}, we need the following simple lemma which justifies the appearance of the constant $24$ in \cref{lem:constr-general-properties}(2). 

\begin{lemma}\label{lem:technical-div} For any $m\nmid 24$, the smallest prime $q$ which does not divide $m$ satisfies $m\nmid q^2-1$.    
\end{lemma}
\begin{proof} The definition of $q$ implies $\prod_{q'<q}q'\mid m$, where the product runs over primes only. The elementary Bonse's inequality states that $\prod_{q'<q}q'>q^2$ if $q>7$; this gives $m>q^2$ if $q\geq 11$. Otherwise, $m\mid q^2-1$ implies that $m\mid 48$. The result is then easy to check by hand.
\end{proof}

\begin{proof}[Proof of \cref{lem:constr-general-properties}] Firstly, the proof of part (1) is the same as that of \cref{lem:independence}. If $AA^{-1}$ contains any edges of $\Sigma$ then $AA^{-1}\setminus\{1\}$ can be contained neither in $S$ nor in $G\setminus S$. Otherwise, the event that a particular element of $AA^{-1}\setminus\{1\}$ lies in $S$ is independent from that of at most one other element (its inverse), and the result follows.

For part (2), suppose $x$ is of some order $m\nmid 24$, $q$ is the smallest prime not dividing $m$, and $A$ is such that $x,x^q,x^{q^2},x^{q^3}\in AA^{-1}$. The equivalence classes $T:=\{x\sim x^{-1},x^q\sim x^{-q},x^{q^2}\sim x^{-q^2},x^{q^3}\sim x^{-q^3}\}$ (not all of which must be distinct) lie, in that order, in a cycle of the auxiliary graph $\Pi$ constructed in \cref{constr:ramsey-cayley-general}. Let $\ell$ be the length of the corresponding cycle. Since $m\nmid 24$, \cref{lem:technical-div} implies that $m\nmid q^2-1$. Therefore $x$ neither equals $x^q$ nor $x^{-q}$, and $\ell>1$.

If $\ell\geq 4$ then $T$ forms four consecutive vertices of a cycle in $\Pi$, and thus contains some edge of the maximal matching $\Sigma$. If $2\leq\ell<4$ then $T$ forms an entire cycle, and thus contains some edge $\Sigma$. In either case, there is an edge between $x^{q^i}\sim x^{-q^i}$ and $x^{q^{i+1}}\sim x^{-q^{i+1}}$ for some $0\leq i<3$ in $\Sigma$, and so if $x^{q^i},x^{q^{i+1}}\in AA^{-1}$ then $A$ is neither a clique nor an independent set in $\Sigma$.
\end{proof}

\cref{thm:ramsey-cayley-general} will now follow, in much the same way as \cref{thm:ramsey-cayley}, from \cref{thm:CFPY-main}, \cref{lem:constr-general-properties}, and the following analogue of \cref{lem:cliques-large-doubling}.

\begin{lemma}\label{lem:cliques-large-doubling-general} There exists an absolute constant $c>0$ for which the following holds.

Let $G$ be a finite group. Let $Q$, $r$, and $n$ be positive integers with $n$ sufficiently large, $r\leq\exp(c\sqrt{\log\log n/\log\log\log n})$, and $Q\leq\exp(\log^{1/27}n)$. Suppose that $G$ has no cyclic subgroup of order $\lcm(1,\ldots,Q)$, no subgroup isomorphic to $(\ZZ/2\ZZ)^r$, and no subgroup isomorphic to $(\ZZ/3\ZZ)^r$. 

Suppose $A\subset G$ of size $n$ satisfies the following property: if $x$ is any element of $G$ whose order does not divide $24$ and $q$ is the smallest prime not dividing the order of $x$, then $\{x,x^q,x^{q^2},x^{q^3}\}\not\subset AA^{-1}$. Then
\[\abs{AA^{-1}}\geq \abs{A}\exp\pr*{\log^{1/118}\abs{A}}.\]
\end{lemma}

We will prove \cref{lem:cliques-large-doubling-general} in the following subsection.

\begin{proof}[Proof of \cref{thm:ramsey-cayley-general}] Let $N:=\abs{G}$ and let $n:=\ceil{24000\log N}$. We may assume that $n$ is sufficiently large. We must show that, as long as $G$ contains no subgroup isomorphic to $(\ZZ/2\ZZ)^r$ or $(\ZZ/3\ZZ)^r$ or $\ZZ/\lcm(1,\ldots,Q)\ZZ$ with $r:=\floor{\exp((\log\log\log N)^{1/3})}$ and $Q:=\floor{\exp((\log\log N)^{1/27})}$, then $G$ has a $24000$-Ramsey Cayley graph. 

We apply \cref{lem:ramsey-cayley-general} with $C_0=2^{3/2}$ and $c_0=1/2$ to the probability distribution on subsets of $G$ given by \cref{constr:ramsey-cayley-general}. Let $S$ be chosen from that distribution and let $\Gamma:=\Cay(G,S)$. Condition (i) of \cref{lem:ramsey-cayley-general}, with our choice of $C_0$ and $c_0$, is precisely \cref{lem:constr-general-properties}(1). It remains to verify condition (ii) of \cref{lem:ramsey-cayley-general}. 

It is enough to show that, for $N$ sufficiently large, any set $A$ of size $n$ which can potentially (i.e., with nonzero probability) be a clique or independent set in $\Gamma$ must satisfy $\abs{AA^{-1}}\geq n\exp(\log^{1/118}n)$. Consider such a set $A$. We first apply \cref{lem:constr-general-properties}(2) to find that, for any $x\in G$ of order not dividing $24$, if $q$ is the smallest prime not dividing the order of $x$, then $x,x^q,x^{q^2},x^{q^3}$ cannot all lie in $AA^{-1}$. Since $n\geq\log N$, the subgroup-avoiding conditions on $G$ given in the statement of \cref{thm:ramsey-cayley-general} are enough to imply the subgroup-avoiding conditions of \cref{lem:cliques-large-doubling-general}. We may thus apply \cref{lem:cliques-large-doubling-general} to $A\subset G$ to conclude the result.
\end{proof}

\subsection{Avoiding patterns in \texorpdfstring{$AA^{-1}$}{AA-1}}\label{sec:avoid-pattern}

To finish the proofs of \cref{thm:ramsey-cayley,thm:ramsey-cayley-general}, we need only prove \cref{lem:cliques-large-doubling,lem:cliques-large-doubling-general}. We will deduce these two lemmas from the following two statements in an abelian setting. The first is sourced from \cite{CFPY}, and follows from the pigeonhole principle and the Pl\"unnecke--Ruzsa inequality; we will prove the second in \cref{sec:abelian}.

\begin{lemma}[{Special case of \cite[Lemma~4.5]{CFPY}}]\label{lem:CFPY-A-A} Let $G_0$ be an abelian group and let $A\subset G_0$ be a subset for which, if $x,2x,4x,8x\in A-A$, then $x=0$. Then
\[\abs{A-A}\geq\abs{A}^{1+1/144}.\]
\end{lemma}

\begin{lemma}\label{lem:find-multiples-general} Let $G_0$ be a finite abelian group. Let $Q$, $r$, and $n$ be positive integers with $n$ sufficiently large, $r\leq\log^{1/3}n$, and $Q\leq\exp(\log^{1/3}n)$. Suppose that $G_0$ has no cyclic subgroup of order $\lcm(1,\ldots,Q)$, no subgroup isomorphic to $(\ZZ/2\ZZ)^r$, and no subgroup isomorphic to $(\ZZ/3\ZZ)^r$.

Let $A$ be a subset of $G_0$, and suppose $A$ satisfies $\abs{A}=n$ and $\abs{A-A}\leq Kn$. Suppose that $K\leq\exp(\log^{1/13}n)$. Then there exists some prime $q$ and some $t\in G_0$ for which
\begin{enumerate}
    \item $t$ has order not dividing $24$,
    \item every prime $q'<q$ divides the order of $t$, but $q$ does not divide the order of $t$, and
    \item $t,qt,q^2t,q^3t\in A-A$.
\end{enumerate}
\end{lemma}

\begin{proof}[Proof of \cref{lem:cliques-large-doubling}] Recall that we have some set $A\subset G$ in some ambient finite group $G$ whose order is neither a multiple of $2^r$ nor $3^r$. Moreover, $A$ satisfies the property that, if $\gcd(\ord(x),6)=1$ and $x,x^2,x^4,x^8\in AA^{-1}$, then $x=1$. We endeavor to lower-bound $\abs{AA^{-1}}$. Let $K:=\abs{AA^{-1}}/\abs{A}$.

We first note that, since $2^r$ cannot divide the order of any subquotient of $G$, the group $G$ has no subquotient isomorphic to $\A_{2r+2}$. This enables us to apply \cref{lem:abelian-struct-near}(1) to $A\subset G$. We find some abelian subgroup $H\leq G$, some element $g\in G$, and some set $A_0\subset gA\cap H$ which satisfies $\abs{A_0A_0^{-1}}\leq(2K)^{O(1)}\abs{A_0}$ and
\begin{equation}\label{eq:A0-size-rc}
\abs{A_0}\geq\exp\pr*{\Omega\pr*{\frac{\log^{1/8}\abs{A}}{\exp(O(\log^2r\log\log r))\log 2K}}}.
\end{equation}
Let $H_1\leq H$ be the maximal subgroup of $H$ whose order is relatively prime to $2$ and $3$ (such a subgroup is unique since $H$ is abelian). Since the order of $H$ is neither a multiple of $2^r$ nor $3^r$, the index of $H_1$ in $H$ is at most $2^{r-1}3^{r-1}$. In particular, we can find some subset $A_1\subset A_0$ which lies in some left-translate of $H_1$ and satisfies $\abs{A_1}\geq 6^{-r}\abs{A_0}$. 

The pattern-avoiding property we assumed of $A$ (that $AA^{-1}$ avoids the pattern $x,x^2,x^4,x^8$ subject to some conditions) is preserved under left-translation and taking subsets, and so $A_1$ also satisfies the same pattern-avoiding property. Moreover, since $A_1A_1^{-1}$ lies in some conjugate of $H_1$, every element of $A_1A_1^{-1}$ has order relatively prime to $6$. We conclude that $A_1$ satisfies the stronger property that $A_1A_1^{-1}$ avoids \emph{all} instances of the pattern $\{x,x^2,x^4,x^8\}$ with $x\neq 1$. Since $A_1$ is contained in an abelian group, we can apply \cref{lem:CFPY-A-A} to show
\[(2K)^{O(1)}\abs{A_0}\geq\abs{A_0A_0^{-1}}\geq\abs{A_1A_1^{-1}}\geq\abs{A_1}^{1+1/144}\geq (6^{-r}\abs{A_0})^{1+1/144}.\]
Rearranging the above gives that $\abs{A_0}\leq 2^{O(r)}(2K)^{O(1)}$. We compare this bound with \eqref{eq:A0-size-rc} to obtain
\[r+\log 2K\gtrsim \log\abs{A_0}\gtrsim\frac{\log^{1/8}\abs{A}}{\exp(O(\log^2r\log\log r))\log 2K}.\]
This is enough to imply the result.    
\end{proof}

\begin{proof}[Proof of \cref{lem:cliques-large-doubling-general}] Recall $n=\abs{A}$. Let $K:=\abs{AA^{-1}}/\abs{A}$, and suppose for the sake of contradiction that $K\leq\exp(\log^{1/118}n)$.

First, we apply \cref{lem:sectional-p-rank}: our assumption that $G$ has no subgroup isomorphic to $(\ZZ/2\ZZ)^r$ implies that $G$ has no subquotient isomorphic to $(\ZZ/2\ZZ)^{2r^2}$. In particular, since the alternating group $\A_d$ has a subgroup isomorphic to $(\ZZ/2\ZZ)^{\floor{d/2}-1}$, we conclude that $G$ has no subquotient isomorphic to $\A_{4r^2+3}$. We apply \cref{lem:abelian-struct-near}(1) to $A\subset G$. We find some abelian subgroup $H\leq G$, some element $g\in G$, and some set $A_0\subset gA\cap H$ which satisfies $\abs{A_0A_0^{-1}}\leq(2K)^{O(1)}\abs{A_0}$ and
\begin{equation}\label{eq:A0-size-rcg}
\abs{A_0}\geq\exp\pr*{\Omega\pr*{\frac{\log^{1/8}\abs{A}}{\exp(O(\log^2r\log\log r))\log 2K}}}.
\end{equation}
Let $n':=\abs{A_0}$. Recall that $K\leq\exp(\log^{1/117}n)$. Choosing $c$ small enough, \eqref{eq:A0-size-rcg} is enough to imply that $n'\geq\exp(\log^{1/9}n)$.

The pattern-avoiding property we assumed of $A$ is preserved under left-translation and taking subsets, and so $A_0\subset H$ also satisfies the same pattern-avoiding property. Moreover, $H$ is a subgroup of $G$, and so it has no cyclic subgroup of order $\lcm(1,\ldots,Q)$ with $Q\leq\exp(\log^{1/27}n)\leq\exp(\log^{1/3}n')$, no subgroup isomorphic to $(\ZZ/2\ZZ)^r$, and no subgroup isomorphic to $(\ZZ/3\ZZ)^r$, where
\[r\leq\exp\pr*{c\sqrt{\frac{\log\log n}{\log\log\log n}}}\leq\log^{1/27}n\leq\log^{1/3}n'.\]
Since $H$ is abelian, we can thus apply \cref{lem:find-multiples-general} to $A_0\subset H$ of size $n'$ to find that $\abs{A_0A_0^{-1}}/\abs{A_0}=(2K)^{O(1)}\geq\exp(\log^{1/13}n')\geq\exp(\log^{1/117}n)$, which contradicts our assumption on $K$.
\end{proof}

\section{Working within abelian groups}\label{sec:abelian}

We have now proven \cref{thm:local-3-AP,thm:ramsey-cayley-general} under the assumption of the purely abelian results \cref{lem:3-AP-abelian,lem:corner-abelian,lem:find-multiples-general}. What remains is to prove these three lemmas, which we do in this section. We make one comment: all of the results in this section have been applied to abelian sets found in \cref{lem:abelian-struct-near}. Since \cref{lem:abelian-struct-near} has quasi-polynomial dependence on the size of the input set, we are always willing to accept quasi-polynomial losses in our parameters throughout the rest of the argument. In particular, Sanders' quasi-polynomial Bogolyubov--Ruzsa lemma \cite{Sanders} plays a prominent role.

\subsection{Finding dense models}

One essential property of sets of small doubling in abelian groups is that their additive structure can be captured by sets which take up a large fraction of some ambient abelian group. In other words, sets with small doubling have \emph{dense models}. Obtaining such dense models is a major step in Ruzsa's proof \cite{RuzsaFre} of Freiman's theorem. We need a version of this step with good quantitative bounds; to this end, we use the main result of \cite{Sanders}. To introduce this result, we first need to define coset progressions. We shall be rather sparse with the material surrounding coset progressions; for more detail and background, we refer the reader to \cite[Section~4.4]{TaoVu}.

\begin{definition}[Proper symmetric coset progression] Let $G$ be an abelian group. A \emph{symmetric coset progression} $M$ in $G$ is a set of the form $H+P$ where $H$ is a subgroup of $G$ and $P$ is a \emph{symmetric generalized arithmetic progression}
\[P=\set[\big]{a_1g_1+a_2g_2+\cdots+a_rg_r:(a_1,\ldots,a_r)\in\br{-N_1,N_1}\times\cdots\times\br{-N_r,N_r}}\]
for some $g_1,\ldots,g_r\in G$ and $N_1,\ldots,N_r\in\ZZ_{>0}$. The number $r$ is called the \emph{rank} of the symmetric coset progression $M$. We say that $M$ is \emph{proper} if the map $H\times[-N_1,N_1]\times\cdots\times[-N_r,N_r]\to M$ defined by
\[(h,a_1,\ldots,a_r)\mapsto h+a_1g_1+a_2g_2+\cdots+a_rg_r\]
is a bijection.
\end{definition}

\begin{theorem}[{\cite[Theorem~1.1]{Sanders}}]\label{thm:sanders-main} Let $A$ be a subset of an abelian group, and suppose $\abs{A-A}\leq K\abs{A}$. Then $2A-2A$ contains a proper symmetric coset progression $M$ of rank $O(\log^62K)$ and size
\[\abs{M}\geq\exp\pr*{-\log^{6+o(1)}2K}\abs{A}.\]
\end{theorem}

As a corollary of \cref{thm:sanders-main}, a simple application of the pigeonhole principle gives that, for some translate $S$ of some proper symmetric coset progression of size comparable to $\abs{A}$, the \emph{density} $\abs{A\cap S}/\abs{S}$ of $A$ within $S$ is large. We need to find sets $S$ with a slightly more precise structure, but we are willing to sacrifice considerably in the size of $S$.

To this end, we define an \emph{interval} in a finite abelian group to be a set of the form $g_1+[m]g_2=\{g_1+g_2,g_1+2g_2,\ldots,g_1+mg_2\}$ for some $g_1,g_2\in G$ and some positive integer $m$ which is at least the order of $g_2$. The following lemma is essentially a corollary of \cref{thm:sanders-main}. It states that any set with small doubling in an abelian group contains a large portion of a (somewhat) large very well-behaved set, either an interval or a translate of an elementary $p$-group (for some prime $p$ which we do not control). The specificity of this lemma will be needed to prove \cref{lem:find-multiples-general}, and the quasi-polynomial loss we incur will not be of substantial consequence in the proofs of \cref{lem:3-AP-abelian,lem:corner-abelian} (in both statements, the exponent in the quasi-polynomial behavior can be improved by using Freiman isomorphisms instead of the following lemma.)

\begin{lemma}\label{lem:dense-on-nice-set} Let $G_0$ be an abelian group. Let $A\subset G_0$ be of size $n$ and suppose that $\abs{A-A}\leq K\abs{A}$ with $K\leq\exp(\log^{1/12}n)$.

There exists a set $I\subset G_0$ such that the following conditions hold:
\begin{enumerate}
    \item $I$ is either (a) an elementary $p$-group or (b) an interval in $G_0$.
    \item $A$ can be covered by at most $\exp(O(\log^62K))\cdot\frac{\abs{A}}{\abs{I}}$ translates of $I$.
    \item $\abs{I}\geq\exp(\Omega(\log^{1/2}n))$.
\end{enumerate}
In particular, one can find a set $S$ which is a translate of such an $I$ for which $\abs{A\cap S}/\abs{S}\geq\exp(-O(\log^62K))$.
\end{lemma}
\begin{proof} By modifying the constant hidden in (3), we can assume $n$ is sufficiently large. Now, by \cref{thm:sanders-main}, we can find a proper coset progression $M=H+P$ contained in $2A-2A$ of size at least $\exp(-O(\log^72K))\abs{A-A}$ and rank $r$ at most $O(\log^62K)$. Let $\eps_0:=\abs{H+P}/\abs{A}$. Using the Pl\"unnecke--Ruzsa inequality to obtain $\abs{3A-2A}\leq K^5\abs{A}$, Ruzsa's covering lemma allows us to cover $A$ by $K^5\eps_0^{-1}$ translates of $(H+P)-(H+P)=H+2P$. Write
\[P=\sum_{j=1}^r[-N_j,N_j]g_j\]
for positive integers $N_1\geq N_2\geq\cdots\geq N_r$ and elements $g_j\in G$. The set $H+2P$ can be covered by $\abs{H+2P}/\abs{I_j}$ translates of the interval $I_j:=[-2N_j,2N_j]g_j$ for any $j$, or alternatively by $\abs{H+2P}/\abs{J}$ translates of any subgroup $J$ of $H$. Since $A$ can be covered by
\[K^5\eps_0^{-1}=\frac{K^5\abs{A}}{\abs{H+P}}\leq\frac{K^52^r\abs{A}}{\abs{H+2P}}\]
translates of $H+2P$, condition (2) holds for any $I_j$ or $J$ as described. We need only find such an $I_j$ or $J$ which is (1) of the desired type and (3) large enough.

Our upper bound on $K$ implies that, for large enough $n$, we have $\abs{H+P}\geq n^{1/2}$. So either $\abs{P}\geq n^{1/4}$ or $\abs{H}\geq n^{1/4}$. In the former case we have, using the properness of $P$, that
\[\abs{I_1}=4N_1+1\geq\left(\prod_{j=1}^r(2N_j+1)\right)^{1/r}=\abs{P}^{1/r}\geq n^{1/(4r)}\geq\exp\pr*{\Omega\pr*{\frac{\log n}{\log^62K}}}\geq\exp\pr*{\Omega(\log^{1/2}n)}\]
by our upper bound on $K$. The result follows in this case.

For the latter case, we consider the structure of $H$. By the fundamental theorem of abelian groups, we can write
\[H\cong\prod_{i=1}^d\ZZ/M_i\ZZ\]
for some positive integer $d$ and some $M_1\mid M_2\mid\cdots\mid M_d$. Let $p$ be any prime factor of $M_1$. Note that $H$ has subgroups isomorphic to each of $\ZZ/M_d\ZZ$ and $(\ZZ/p\ZZ)^d$. Let $L:=\max(M_d,p^d)$, so that $H$ has either a cyclic subgroup or an elementary $p$-subgroup $J$ of order $L$. We have
\[\abs{H}\leq (M_d)^d\leq L^{\log_2L},\]
and so $\log^2L=\Omega(\log\abs{H})=\Omega(\log n)$. The result follows. (Note that a cyclic group is by definition an interval.)
\end{proof}

\subsection{Three-term arithmetic progressions and corners}

\cref{lem:dense-on-nice-set} provides enough structural information on sets with small doubling in abelian groups to enable us to deduce \cref{lem:3-AP-abelian,lem:corner-abelian} from the more standard results on the density of $3$-AP-free and corner-free subsets in abelian groups. The input we use about three-term arithmetic progressions in the abelian setting is the following.

\begin{theorem}[Slight generalization of {\cite[Theorem~21]{BloomSisask}}]\label{thm:general-kelley-meka} Let $G$ be an abelian group of odd order $N$ and suppose $A\subset G$ is a set avoiding nontrivial instances of the pattern $\{x,x+y,x+2y\}$. Then
\[\abs{A}\leq N\exp\pr*{-\Omega\pr*{\log^{1/12}N}}.\]
\end{theorem}

As far as we know, the statement of \cref{thm:general-kelley-meka} has not appeared exactly in the literature surrounding the recent improvement of Kelley and Meka to Roth's theorem. However, the proof of \cite[Theorem~21]{BloomSisask}, which is nominally specific to the case $G\cong\ZZ/N\ZZ$, applies verbatim to our more general setting. (The entire argument in \cite{BloomSisask}, up until the ``extraction'' of the statement of their Theorem 21, is in the language of Bohr sets, and is thus completely agnostic to the structure of the ambient abelian group.)

The input we use about corners is the following, the main result of \cite{JLLOS}.

\begin{theorem}[{\cite[Theorem~1.1]{JLLOS}}]\label{thm:corners} Let $G$ be an abelian group of order $N$ and suppose $A\subset G\times G$ is a set avoiding nontrivial instances of the pattern $\{(x,y),(x+d,y),(x,y+d)\}$. Then
\[\abs{A}\leq N^2\exp\pr*{-\Omega\pr*{\log^{1/600}N}}.\]
\end{theorem}

\begin{proof}[Proof of \cref{lem:3-AP-abelian}] We are given some set $A\subset G_0$ with no nontrivial $3$-term arithmetic progressions and we wish to lower-bound $\abs{A-A}/\abs{A}$. Let $n:=\abs{A}$ and $K:=\abs{A-A}/\abs{A}$, and suppose $K\leq\exp(\log^{1/12}n)$ (otherwise, we are done). We begin by applying \cref{lem:dense-on-nice-set} to find some set $S\subset G_0$ which is either (a) a translate of an elementary $p$-group or (b) an interval in $G_0$ and satisfies
\begin{align}
\label{eq:A-dense-in-S}\frac{\abs{A\cap S}}{\abs{S}}&\geq\exp\pr*{-O\pr*{\log^62K}};\\
\label{eq:S-large}\abs{S}&\geq\exp\pr*{\Omega\pr*{\log^{1/2}n}}.
\end{align}
We claim that
\begin{equation}\label{eq:A-sparse-in-S}
\frac{\abs{A\cap S}}{\abs{S}}\leq\exp\pr*{-\Omega\pr*{\log^{1/12}\abs{S}}}.
\end{equation}
We first show this in case (a). Let $s\in G_0$ and $H$ be an elementary $p$-group such that $S=s+H$, and let $B:=(A-s)\cap H$ so that $\abs{B}=\abs{A\cap S}$. If $p=2$ then \eqref{eq:A-dense-in-S} and \eqref{eq:S-large} are enough to imply that $A$ contains two elements in the same translate of some elementary $2$-group. Therefore $A$ contains two distinct elements $x$ and $z$ whose difference $y:=z-x$ has order $2$, and thus $A$ contains each of $x$, $x+y=z$, and $x+2y=x$, and we obtain a contradiction. So $p\neq 2$, and we apply \cref{thm:general-kelley-meka} to $B\subset H$ to obtain \eqref{eq:A-sparse-in-S}. 

If we are in case (b), let $m:=\abs{S}$ and let $g_1,g_2\in G$ be such that $S=\{g_1+ig_2:1\leq i\leq m\}$. Let $B\subset[m]\subset\ZZ/(2m+1)\ZZ$ be the set of $i\in[m]$ for which $g_1+ig_2\in A\cap S$. Since $A\cap S$ contains no three-term arithmetic progressions, and the ambient group of $B$ has been chosen large enough, the set $B$ contains no three-term arithmetic progressions either, and we conclude \eqref{eq:A-sparse-in-S} from \cref{thm:general-kelley-meka} applied to $B\subset\ZZ/(2m+1)\ZZ$.

In either case, \eqref{eq:A-sparse-in-S} holds. Combining this with \eqref{eq:A-dense-in-S} and \eqref{eq:S-large} gives
\[\log^{1/12}\abs{S}\lesssim\log^62K\implies\log 2K\gtrsim\log^{1/72}\abs{S}\gtrsim\log^{1/144}\abs{A}.\]
The result follows.
\end{proof}

\begin{proof}[Proof of \cref{lem:corner-abelian}] We are given some sets $B\subset G_0$ and $A\subset B\times B$ such that $A$ has no nontrivial corners and we wish to show that either $B$ has large doubling or $A$ is sparse in $B\times B$. Let $K:=\abs{B-B}/\abs{B}$ and $\eps:=\abs{A}/\abs{B}^2$. We may assume that $K\leq\exp(\log^{1/12}\abs{B})$, as otherwise the trivial bound $\abs{A}\leq\abs{B}^2$ is enough to conclude the result. We can thus apply \cref{lem:dense-on-nice-set} to $B\subset G_0$ to find some set $I\subset G_0$, either (a) an elementary $p$-group or (b) an interval in $G_0$, such that $B\subset Z+I$ for some set $Z$ of size at most $\exp(O(\log^62K))\cdot\frac{\abs{B}}{\abs{I}}$, and moreover $\abs{I}\geq\exp(\Omega(\log^{1/2}\abs{B}))$. Since $A\subset B\times B\subset (Z\times Z)+(I\times I)$,
we can find some $z_1,z_2$ for which
\begin{equation}\label{eq:A'-dense-in-I}
\frac{\abs[\big]{A\cap ((z_1+I)\times(z_2+I))}}{\abs{I}^2}\geq\frac{\abs{A}}{\abs{I}^2\abs{Z}^2}=\eps\cdot\pr*{\frac{\abs{B}}{\abs{I}\abs{Z}}}^2\geq\eps\cdot\exp\pr*{-O(\log^62K)}.
\end{equation}
Let $A':=(A-(z_1,z_2))\cap(I\times I)$, so that $A'$ is relatively dense in $I\times I$. Since $A'$ is contained in a translate of a corner-free set, $A'$ is also corner-free. If $I$ is an elementary $p$-group, then we can use \cref{thm:corners} to conclude directly that
\begin{equation}\label{eq:A'-sparse-in-I}
\abs{A'}\leq \abs{I}^2\exp\pr*{-\Omega\pr*{\log^{1/600}\abs{I}}}
\end{equation}
If $I$ is an interval of length $m$, we can (as in the proof of \cref{lem:3-AP-abelian} above) embed $A'$ as a corner-free subset of $(\ZZ/(2m+1)\ZZ)^2$ and obtain \eqref{eq:A'-sparse-in-I}. In either case, combining \eqref{eq:A'-sparse-in-I} with \eqref{eq:A'-dense-in-I} and our size assumption on $I$ gives
\[\log^{1/1200}\abs{B}\lesssim\log^{1/600}\abs{I}\lesssim\log\pr*{\frac{\abs{I}^2}{\abs{A'}}}\lesssim\log^62K+\log\eps^{-1}\lesssim\log^6(2K\eps^{-1}).\]
Since $\frac{\abs{B-B}^2}{\abs{A}}=K^2\eps^{-1}$, we conclude the result.
\end{proof}

\subsection{Finding patterns in \texorpdfstring{$A-A$}{A-A}}

In contrast to three-term arithmetic progressions and corners, the somewhat peculiar pattern specified in \cref{lem:find-multiples-general} has (to our knowledge) not appeared before in the literature. Fortunately, the pattern we must deal with is a variant of the patterns treated in \cite[Lemma~4.5]{CFPY}, and can, more or less, be found with an application of the pigeonhole principle. (See also \cite{CRS}, where a similar argument is used to find long arithmetic progressions in difference sets.) In the integers, this can be done directly by picking the relevant prime $q$ in advance. 

\begin{lemma}\label{lem:find-multiples-Z} Let $Q$, $R$, and $N$ be positive integers. Let $A\subset[N]$ be of size at least $\eps N$. For each prime $q<Q$, designate a nonnegative integer $s_q$, and suppose that not all $q^{s_q}$ exceed $Q$. 

Suppose that $N>8\eps^{-12}Q^{10}R$. Then there exists some prime $q<Q$ and at least $R$ values of $t\in\ZZ\setminus\{0\}$ such that
\begin{enumerate}
    \item $q^{s_q}\mid t$ while $(q')^{s_{q'}}\nmid t$ for any $q'<q$, and
    \item $t,qt,q^2t,q^3t\in A-A$.
\end{enumerate} 
\end{lemma}
\begin{proof} Let $\Delta:=2Q^7\eps^{-8}$ and let $q$ be the smallest prime for which $q^{s_q}\leq\Delta$. (Such a prime is guaranteed to exist since some $q^{s_q}$ is at most $Q$.) Consider the map $\phi\colon A\times A\times A\times A\to (A+qA)\times (A+qA)\times (A+qA)$ given by
\[(a_0,a_1,a_2,a_3)\mapsto (qa_0+a_1,qa_1+a_2,qa_2+a_3).\]
We have
\[\frac{\abs{A}^4}{\abs{A+qA}^3}\geq\frac{(\eps N)^4}{((q+1)N)^3}\geq\eps^4Q^{-3}N,\]
so there exists some $(u_1,u_2,u_3)\in(A+qA)\times(A+qA)\times(A+qA)$ which can be written as $\phi(a_0,a_1,a_2,a_3)$ for at least $\eps^4Q^{-3}N$ quadruples $(a_0,a_1,a_2,a_3)$. No two of these quadruples have the same first coordinate, as we have $a_1=u_1-qa_0$ and $a_2=u_2-qa_1$ and $a_3=u_3-qa_2$ for each such quadruple. Let $A_0$ be the set of first coordinates of these quadruples, so that $\abs{A_0}\geq\eps^4Q^{-3}N$. If $a_0,a_0'\in A_0$ then there exist $a_1,a_2,a_3,a_1',a_2',a_3'\in A$ for which, for each $1\leq i\leq 3$,
\[q^i(a_0-a_0')=q^ia_0-q^ia_0'=(u_i-a_i)-(u_i-a_i')=a_i'-a_i\in A-A.\]
So, any $t:=a_0-a_0'$ for $a_0,a_0'\in A_0$ satisfies condition (2). 

Now, for any integer $m$, write $f(m)$ for the number of ordered pairs of distinct elements of $A_0$ which are equivalent modulo $m$. We have
\begin{align*}
f(q^{s_q})&\geq q^{-s_q}\abs{A_0}^2-\abs{A_0};\\
f\left((q')^{s_{q'}}q^{s_q}\right)&\leq (q')^{-s_{q'}}q^{-s_q}N^2\leq\Delta^{-1}q^{-s_q}N^2\quad\text{for each}\quad q'<q,
\end{align*}
where we have used the Cauchy--Schwarz inequality on the first line and the fact that $A_0\subset[N]$, as well as our choice of $q$, on the second line. As a result, we have
\begin{align*}
\#\left\{(a_0,a_0')\in A_0:\begin{matrix}q^{s_q}\mid a_0-a_0',\\a_0=a_0'\text{ or }\\(q')^{s_{q'}}\nmid a_0-a_0'\text{ for all }q'<q\end{matrix}\right\}
&\geq \abs{A_0}+f(q^{s_q})-\sum_{q'<q}f((q')^{s_{q'}}q^{s_q})\\
&\geq q^{-s_q}\left(\abs{A_0}^2-\frac Q\Delta N^2\right)\\
&\geq \Delta^{-1}\left(\eps^4Q^{-3}\abs{A_0}N-\Delta^{-1}\eps^{-4}Q^4\abs{A_0}N\right)\\
&=\frac14Q^{-10}\eps^{12}\abs{A_0}N.
\end{align*}
where we have used $q<Q$ and $q^{s_q}\leq\Delta$, as well as our lower bound $\abs{A_0}\geq\eps^4Q^{-3}N$, to simplify. Since $N\geq 8\eps^{-12}Q^{10}R\geq 4\eps^{-12}Q^{10}(R+1)$, the set above is of size at least $(R+1)\abs{A_0}$. In particular, there exists some $a_0\in A_0$ for which there exist at least $R$ values $a_0'\in A_0$ which are equivalent to $a_0$ modulo $q^{s_q}$ but not equivalent to $a_0$ modulo $(q')^{s_{q'}}$ for each $q'<q$. Taking our values of $t$ to be $a_0-a_0'$ for each such $a_0'$, the result follows.
\end{proof}

We remark that, for $Q\lesssim\log\log n$, \cref{lem:find-multiples-Z} essentially follows from \cite[Lemma~4.5]{CFPY}, while for $Q\lesssim\log^{1/3}n$ this follows from a result of Croot, Ruzsa, and Schoen \cite{CRS} (both of these arguments have a similar flavor). Arguments along both these lines proceed by finding (using the pigeonhole principle) many $t$ for which $t,2t,3t,\ldots,Q^3t$ are all in $A-A$. So, condition (3) is certainly satisfied regardless of what value of $q$ is forced by condition (2). Our argument improves upon \cite{CRS} by an exponential factor by pinpointing exactly which value of $q$ is forced.

Using \cref{lem:find-multiples-Z} in the integer setting, we are able to prove \cref{lem:find-multiples-general}.

\begin{proof}[Proof of \cref{lem:find-multiples-general}] Since $G_0$ is abelian, it is the direct product of its $p$-parts, each of which has a cyclic subgroup $\ZZ/p\ZZ$. We conclude from the fact that $G_0$ has no cyclic subgroup of order $\lcm(1,\ldots,Q)$ that there exists some prime $q_0\leq Q$ for which the $q_0$-part of $G_0$ has no element of order exceeding $Q$. 

Let $S$ be the subset of $G_0$ given from \cref{lem:dense-on-nice-set} arising from $A\subset G_0$, so that $S$ is a translate of either (a) an elementary $p$-subgroup of $G_0$ or (b) a centered interval in a cyclic subgroup of $G_0$, and moreover
\begin{align}
\frac{\abs{A\cap S}}{\abs{S}}&\geq\exp\pr*{-O(\log^62K)}\label{eq:A-dense-in-S-fm}\\
\abs{S}&\geq\exp\pr*{\Omega(\log^{1/2}n)}.\label{eq:S-big}
\end{align}
We treat the cases in turn. 

First consider case (a). Let $S=x+H$ with $H$ a $p$-group for some prime $p$. By \eqref{eq:S-big} and our assumption that $G_0$ has no subgroup isomorphic to $(\ZZ/2\ZZ)^r$ or $(\ZZ/3\ZZ)^r$, we have $p\not\in\{2,3\}$. Let $q=2$, so that every element of $H$ has order relatively prime to $q$. We now essentially repeat part of the proof of \cref{lem:find-multiples-Z}. Let $A':=A\cap S$ and let $\phi\colon A'\times A'\times A'\times A'\to 3A'\times 3A'\times 3A'$ be defined by
\[\phi(a_0,a_1,a_2,a_3)=(2a_0+a_1,2a_1+a_2,2a_2+a_3).\]
Since $\abs{3A'}\leq \abs{3S}=\abs{S}$ and $\abs{A'}^4>\abs{S}^3$ by \eqref{eq:A-dense-in-S-fm}, \eqref{eq:S-big}, and our assumption $K\leq\exp(\log^{1/13}n)$, the map $\phi$ cannot be injective. We conclude that there exist $a_0,\ldots,a_3,a_0',\ldots,a_3'\in A'$ satisfying $2a_{i-1}+a_i=2a_{i-1}'+a_i'$ for each $1\leq i\leq 3$ and $(a_0,a_1,a_2,a_3)\neq(a_0',a_1',a_2',a_3')$. Let $t:=a_0-a_0'\neq 0$ so that $2^it=a_i'-a_i$ lies in $A-A$ for each $1\leq i\leq 3$; hence $t$ satisfies condition (3). Since the order of $t$ is a (nontrivial) power of $p\geq 5$, conditions (1) and (2) are also satisfied. This finishes the proof in case (a).

We now treat case (b). Let $m:=\abs{S}$, and let $g_1,g_2\in G$ be such that $S=\{g_1+ig_2:1\leq i\leq m\}$, where $g_2$ has order at least $m$. Write $H:=\langle g_2\rangle$ for the cyclic subgroup of $G_0$ generated by $H$ and $N$ for the order of $H$. Let $\pi$ be the projection map $\ZZ\to H$ given by $1\mapsto g_2$ and let $\pi_0$ be its restriction to $[m]$, so that $\pi_0$ is a bijection between $[m]$ and $S-g_1$. Let $A':=\{\pi_0^{-1}(a-g_1):a\in A\}\subset[m]$ and $\eps:=\frac{\abs{A\cap S}}{\abs{S}}=\abs{A'}/m$. For each prime $q<Q$, let $q^{s_q}$ be the largest power of $q$ dividing $N$. Let $R:=49$.

We now apply \cref{lem:find-multiples-Z} to $A'\subset[m]$ and $(s_q)_q$, noting that $q_0^{s_{q_0}}\leq Q$ and
\[m=\abs{S}\geq\exp\pr*{\Omega(\log^{1/2}n)}\geq 400\exp\pr*{\Omega(\log^62K)+10\log^{1/3}n}\geq 8\eps^{-12}Q^{10}R\]
by \eqref{eq:A-dense-in-S-fm}, so the hypotheses of \cref{lem:find-multiples-Z} are satisfied. We obtain the existence of some $q<Q$ and at least $49$ elements $t\in\ZZ\setminus\{0\}$ for which $t,qt,q^2t,q^3t\in A'-A'$ and $q^{s_q}\mid t$ and $(q')^{s_{q'}}\nmid t$ for any $q'<q$. For $t\in\ZZ$, a prime $q$ divides the order of $\pi(t)$ if and only if $q^{s_q}\nmid t$. Since the fibers of $\pi$ on $A'-A'\subset[-(m-1),m-1]$ are of size at most $2$, there are at least $25$ elements $u$ of $H$ for which $u,qu,q^2u,q^3u\in\pi(A')-\pi(A')=A-A$ and $q$ does not divide the order of $u$ and every prime $q'<q$ divides the order of $u$. Conditions (2) and (3) of our conclusion are thus satisfied for any of these $25$ values of $u$. The result follows from the fact that $H$, as a cyclic group, has at most $24$ elements whose orders divide $24$.
\end{proof}

\section{Conclusion and open problems}\label{sec:conclusion}

We conclude with some conjectures and some comments on the limitations of our methods.

\subsection{Abelian substructures}\label{sec:conc-ab}

Perhaps the most natural question arising from our main results is to what extent our assumption on $G$ in \cref{thm:gen-abelian-struct} (i.e., that $G$ avoids subquotients isomorphic to $\A_d$) is necessary. While our proof makes crucial use of this assumption (to control the dimension of the ambient general linear group in our application of \cref{thm:GL}), we conjecture that it is in fact unnecessary.

\begin{conjecture}\label{conj:gen-abelian-struct} There exists an absolute constant $c>0$ for which the following holds. If $G$ is a group and $A$ is a finite symmetric set satisfying $\abs{A^3}\leq K\abs{A}$, then there exists a commuting set $T\subset A^4$ satisfying
\[\abs{T}\geq\exp\pr*{\Omega\pr*{\frac{\log^{c}\abs{A}}{\log2K}}}.\]
\end{conjecture}

\cref{conj:gen-abelian-struct} would allow us to remove all conditions on the ambient group $G$ in \cref{thm:local-3-AP}. It would also allow a quantitative strengthening of \cref{thm:ramsey-cayley,thm:ramsey-cayley-general}, although our methods are not strong enough (as we detail in the following subsection) to remove any of the ``subgroup-avoidance'' conditions in \cref{thm:ramsey-cayley-general} entirely. 

Specialized to the $K=1$ case, \cref{conj:gen-abelian-struct} would give a weak form of \cref{thm:Pyber}. To our knowledge, no proof of \cref{thm:Pyber} (or any super-logarithmic bound therein) is known which avoids the use of the chief series and the classification of finite simple groups. So, a proof of \cref{conj:gen-abelian-struct} would likely require a close quantitative connection between approximate groups and genuine groups, a generalization of chief series and the classification of finite simple groups to approximate groups, or a new, classification-free proof of \cref{thm:Pyber}.

We give two remarks on weaker versions of \cref{conj:gen-abelian-struct} which can be proven with existing tools.

\begin{remark} If the ambient group $G$ is given to be torsion-free, then Sanders' result \cref{prop:sanders-quasipoly} is enough to prove an analogue of \cref{conj:gen-abelian-struct} with $\abs{T}\geq\ c\log^{1/2}\abs{A}/\log 2K$. Indeed, \cref{prop:sanders-quasipoly} gives for $t=c\log^{1/2}\abs{A}/\log 2K$ the existence of some nontrivial set $S$ with $S^t\subset A^4$. If $G$ is torsion-free then taking $T=\langle s\rangle\cap A^4$ for any $s\in S\setminus\{1\}$ gives the result.    
\end{remark}

\begin{remark} If one is willing to sacrifice the ability of $K$ to grow with $\abs{A}$, one can prove an analogue of \cref{conj:gen-abelian-struct}. Specifically, the Breuillard--Green--Tao structure theorem for approximate groups \cite[Theorem~1.6]{BGTGeneral} can be used, along with our work, to give a proof of \cref{conj:gen-abelian-struct} in the case where $\abs{A}\geq F(K)$ for some (ineffective) growth function $F(K)$.

We give a rough sketch of the argument. First, as in our proofs, we find some approximate group $B$ of size comparable to $A$ for which $B^t\subset A^4$ for some $t\geq\log^{\Omega(1)}\abs{A}$. Next, \cite[Theorem~1.6]{BGTGeneral} gives the existence of subgroups $H\unlhd G_0\leq G$ of the ambient group for which $B^4\supset H$ and $G_0/H$ is nilpotent and $B$ is covered by $O_K(1)$ left-translates of $G_0$. If $\abs{A}$ is large enough in terms of $K$, we have $\abs{B^2\cap G_0}\geq\abs{A}^{1/2}$. If $\abs{H}=\abs{B^4\cap H}=\abs{A^4\cap H}$ is of at least quasi-polynomial size in $\abs{A}$, we can use Pyber's result to conclude. Otherwise, we apply \cref{lem:commutable-for-nice-families}(1) and the fact that $G_0/H$ is nilpotent (and thus solvable) to find some commuting dissociated set $D\subset B^2\cap G_0$ of size polylogarithmic in $\abs{A}$. We can then conclude from the argument in the proof of \cref{lem:commutable-to-abelian-substruct}.
\end{remark}

We also believe that it is possible to improve \cref{thm:GL} to require only polynomial, rather than quasi-polynomial, dependence on the dimension of the ambient group.

\begin{conjecture}\label{conj:GL} There exist some polynomially bounded functions $M_1,M_2\colon\ZZ_{>0}\to\RR$ for which the following holds. Let $\FF$ be any field, let $d\geq 1$, and let $A\subset\GL_d(\FF)$ be a $K$-approximate group for some $K\geq 1$. Then there exists some abelian subgroup $H\leq\GL_d(\FF)$ for which
\[\abs{A^2\cap H}\geq\frac1{(2K)^{M_1(d)}}\abs{A^2}^{1/M_2(d)}.\]
\end{conjecture}

\begin{remark} The function $M_2$ in \cref{conj:GL} must grow at least polynomially. In particular, it is not possible for $M_2(d)$ to grow more slowly than $d^{1/3}$. This can be seen using the construction of Ol'shanskii in \cite{Olshanskii} of groups of order $p^n$ with no abelian subgroups of order exceeding $p^{c\sqrt n}$. We sketch the argument:

Let $n=m+k$ with $k$ as small as possible subject to $k^2>2m$. Ol'shanskii selects skew-symmetric bilinear forms $\alpha_1,\ldots,\alpha_k$ on $(\ZZ/p\ZZ)^m$ uniformly at random, and constructs a group $G_{m,k}$ of order $p^n$ as follows: the underlying set of $G_{m,k}$ is $(\ZZ/p\ZZ)^k\times(\ZZ/p\ZZ)^m$, and the group operation is
\[(w,v)\cdot(w',v')=\left(w+w'+\left(\alpha_i(v,v')\right)_{1\leq i\leq k},v+v'\right).\]
By \cite[Lemma~2]{Olshanskii}, such a group has (with high probability) no abelian subgroup of order $p^{2k}$.

Let $B_i$ be a matrix satisfying $\alpha(v,v')=v^\intercal B_iv'$. The group $G_{m,k}$ embeds in $\GL_d(\FF_p)$ for $d:=(m+2)k=\Theta(n^{3/2})$ via the map
\[(w,v)\mapsto\bigoplus_{i=1}^k\begin{pmatrix}1&v^\intercal&w_i\\&I_m&B_iv\\&&1\end{pmatrix}.\]
However, for any abelian subgroup $H$ of $\GL_d(\FF_p)$ we have
\[\abs{G_{m,k}^2\cap H}\leq p^{2k}\leq p^{4\sqrt n}\leq p^{O(d^{1/3})}.\]
Observe $\abs{G_{m,k}}=p^n=p^{\Theta(d^{2/3})}$. Taking $p$ large, we conclude that, if \cref{conj:GL} is true for some $M_1,M_2$, then we must have $M_2(d)\geq cd^{1/3}$.
\end{remark}

There are a handful of aspects of our argument that prevent us from being able to obtain a bound as strong as \cref{conj:GL}. One is the presence of the $\log\log d$ factor in \cref{prop:reg-el}; if we were able to replace this factor by a constant, we would correspondingly save a factor of $\log\log d$ in the exponent of $M_2(d)$. (In fact, it is plausible that a quantitatively much stronger result may hold; see \cref{rmk:reg-el}.) The more significant issue is the fact that we have only worked with varieties of degree one. (For example, the variety $\Conj(a)$ always has codimension at least $d$, but we never consider this variety directly, always replacing it by the affine subspace $\{x:\tr x=\tr a\}$ in which it is contained.) This has the benefit of avoiding quantitative penalties which scale exponentially in $d$, as are encountered in \cite{BDH} (in which the authors prove strong quantitative bounds on the product theorem in linear groups). However, it means that we are unable to fully exploit the ``$a$ relatively irregular'' case. If $a$ is regular semisimple, then $\Conj(a)$ has codimension exactly $d$, and we do not lose so much by passing to $\{x:\tr x=\tr a\}$. However, if $a$ has an eigenspace of dimension on the order of $d$, then $\Conj(a)$ has codimension $\Omega(d^2)$, and the Larsen--Pink inequality tells us that (in the absence of additional structure) $\abs{A^3\cap\Conj(a)}$ should not be much larger than $\abs{A}^{1-\delta}$ for some constant $\delta$. We were not able to exploit the structure of the small variety $\Conj(a)$ to obtain a strong enough quantitative version of this statement for our purposes.

\subsection{Ramsey Cayley graphs}\label{sec:conc-ramsey}

The overarching question surrounding our results \cref{thm:ramsey-cayley,thm:ramsey-cayley-general} on the existence of Ramsey Cayley graphs is, of course, Alon's conjecture (\cref{conj:alon}) --- that, for some absolute constant $C$, every finite group has a $C$-Ramsey Cayley graph. We have proven this conjecture for all groups which possess no large subgroup of the form $(\ZZ/2\ZZ)^r$ or $(\ZZ/3\ZZ)^r$ or $\ZZ/\lcm(1,\ldots,Q)\ZZ$. We briefly explain (qualitatively) why our methods are unable to handle these subgroups, even in the abelian setting.

The difficulty inherent in elementary abelian groups was first noticed by Green \cite{Green}. For $p$ a fixed prime and $G=(\ZZ/p\ZZ)^n$, a uniformly random Cayley graph on $G$ has cliques of size $\Theta(\log\abs{G}\log\log\abs{G})$ with high probability; moreover, such a clique can be taken to be a subgroup of $G$. (For $p=2$ this is \cite[Theorem~9]{Green}.) This difficulty was reiterated by Conlon, Fox, Pham, and Yepremyan in \cite{CFPY}; their main combinatorial results on Cayley graphs of $G$ are tight in the $G=(\ZZ/p\ZZ)^n$ setting. For $p>3$, the construction in \cite{CFPY} of a non-uniform random Cayley graph presents a way to deterministically ensure that no nontrivial subgroup of $(\ZZ/p\ZZ)^n$ can be a clique or independent set. However, for $p\in\{2,3\}$, there are no clear ways to modify the uniformly random graph to avoid the appearance of large subgroups as cliques. Our slightly more general \cref{constr:ramsey-cayley-general} presents no way to surmount this obstacle, and so we encounter the same issue. A necessary, but perhaps not sufficient, condition for the groups $G=(\ZZ/2\ZZ)^n$ to have Ramsey Cayley graphs was proposed in \cite[Conjecture~5.1]{CFPY}.

The avoidance of $\ZZ/\lcm(1,\ldots,L)\ZZ$ may be more particular to our methods. To make a construction like \cref{constr:ramsey-cayley-general} work for some group $G$, it must hold that, for many $g\in G$, there exists some small prime $q$ which does not divide the order of $g$. This fails to hold for the group $\ZZ/\lcm(1,\ldots,L)\ZZ$. More specifically, the $\ZZ/\lcm(1,\ldots,L)\ZZ$-avoidance assumption is used in \cref{lem:find-multiples-Z}, wherein the resulting condition (that some $q^{s_q}$ is not too large) is essential. \cref{lem:find-multiples-Z} is applied once one has already restricted oneself to a coset of an abelian subgroup; this is why we must take the parameter $L$ to be quasi-polynomial, rather than polynomial, in $\log N$.

In some specific cases, this issue can be remedied by choosing a slightly different construction. In the groups $G=(\ZZ/\lcm(1,\ldots,n)\ZZ)^k$ \emph{for bounded $k$}, one can use a construction related to \cref{constr:ramsey-cayley}. For any $x\in G$, there are at most $2^k$ elements $y\in G$ for which $2y=x$. This allows for a randomized construction wherein (a slightly weaker version of) \cref{lem:constr-1248-avoid} holds, at the cost of imposing, for many $x$, that those elements $y$ satisfying $2y=x$ deterministically either all lie inside $S$ (if $x\not\in S$) or outside $S$ (if $x\in S$). We can apply \cref{lem:ramsey-cayley-general} to this construction, but we must take $c_0=\Theta(2^{-k})$. We thus find that $G$ has a $O(2^k)$-Ramsey Cayley graph. 

Based on this discussion, we highlight two particular cases of Alon's conjecture. The first case, in which the ambient group is abelian, would follow from \cref{thm:ramsey-cayley-general} if the $\ZZ/\lcm(1,\ldots,Q)\ZZ$-avoidance condition could be removed.

\begin{question}\label{qn:unbalanced-abelian-rc} Does the family of abelian groups
\[G_n:=\left(\ZZ/\lcm\big(1,\ldots,2^{2^n}\big)\ZZ\right)^n\]
have Ramsey Cayley graphs?
\end{question}

The second is the case of symmetric groups. Symmetric groups pose a natural barrier to our methods in every sense: they possess large elementary abelian subgroups, large cyclic subgroups, and large composition rank. We expect that proving Alon's conjecture for symmetric groups will require entirely new ideas.

\begin{question}\label{qn:symmetric-rc} Does the family of symmetric groups have Ramsey Cayley graphs?
\end{question}

\subsection{Miscellany}\label{sec:conc-misc}

There are two intermediate statements appearing in the body of this work which we believe may be of independent interest and which we wish to highlight in this concluding section.

The first is \cref{prop:reg-el}. A natural generalization of \cref{prop:reg-el}, removing the $\log\log d$ factor, is the following.

\begin{conjecture}\label{conj:reg-el} There exists some constant $c>0$ for which the following holds. Let $\FF$ be an algebraically closed field and let $A\subset\GL_d(\FF)$ be a $K$-approximate group. Suppose that every element $a\in A^2$ has an eigenspace of dimension $(1-c)d$. Then there exist two subspaces $V_1\subset V_2\subset\FF^d$ with $\dim V_2-\dim V_1>cd$ for which at least $(2K)^{-O(d^{O(1)})}\abs{A^2}$ elements $a$ of $A^2$ satisfy the following property: there exists some $\lambda\in\FF$ for which $(a-\lambda)V_2\subset V_1$.
\end{conjecture}

No conclusion of the above form can be obtained if we relax our assumption to the statement that every element of $A^2$ has an eigenspace of dimension $m$ for any $m<d/2$. This can be seen by considering the diagonal embedding
\[\GL_2(\FF_q)\hookrightarrow\bigoplus_{i=1}^k\GL_2(\FF_q)\leq\GL_{2k}(\FF_q).\]

We also ask a question directly motivated by what we need to prove \cref{thm:GL}. While \cref{prop:reg-el} presents (qualitatively) a classification for when an approximate groups fails to contain a not-too-irregular element, we do not use its full strength. All we need from its conclusion is that $(2K)^{-O(d^{O(1)})}\abs{A^2}$ many elements of $A^2$ lie in some proper affine subspace of $\Mat_d(\FF)$. 

\begin{question}\label{qn:reg-el} Given $d$, what is the smallest $m$ for which the following holds? 

Let $\FF$ be an algebraically closed field and let $A\subset\GL_d(\FF)$ be a $K$-approximate group. Suppose that every element $a\in A^2$ has an eigenspace of dimension at least $m$. Then there exists some affine subspace $\bU$ of $\Mat_d(\FF)$ containing at least $(2K)^{-O(d^{O(1)})}\abs{A^2}$ elements of $A^2$.
\end{question}

This question may be interesting even in the case of genuine groups. (Here, Burnside's theorem more or less requires such an affine subspace to be contained in the stabilizer of some nontrivial flag on $\FF^d$.) 

\begin{remark}\label{rmk:reg-el} For some constant $c$ and infinitely many $d$, such an $m$ in \cref{qn:reg-el} must exceed $d\cdot e^{-c\sqrt{\log d}\log\log d}$, even for $K=1$. We sketch a proof. Let $n$ be an odd positive integer, and consider the irreducible representation $\rho$ of $\GL_n$ associated to the weight $(n-1,n-2,\ldots,1,0)$. This representation has dimension $d:=2^{\binom n2}$. Let $q$ be a large prime power, and let $G\leq\GL_d(\FF_q)$ be the image of $\GL_n(\FF_q)$ in this representation. One can compute using Schur polynomials that, if $t_1,\ldots,t_n$ are the eigenvalues of some $a\in\GL_n(\FF_q)$, then the associated $\rho(a)\in G$ has eigenvalue $t_1^{e_1}t_2^{e_2}\cdots t_n^{e_n}$ with multiplicity equal to the coefficient of $x_1^{e_1}\cdots x_n^{e_n}$ in the expansion of 
\[\prod_{1\leq i<j\leq n}(x_i+x_j).\]
Let $m$ be the coefficient of $x_1^{(n-1)/2}\cdots x_n^{(n-1)/2}$ in this expansion. Every $\rho(a)$ has an eigenvalue of multiplicity $m$, and in fact one can show (by the lower semicontinuity of matrix rank) that every $\rho(a)$ in fact has an \emph{eigenspace} of dimension $m$. The number $m$ counts the number of regular tournaments on $n$ vertices; the asymptotic $m=d\cdot e^{-\Theta(\sqrt{\log d}\log\log d)}$ then follows from a theorem of McKay \cite[Theorem~3.1]{McKay}.

What remains is to show that, for $q$ sufficiently large, a vanishing proportion of $G$ lives in any affine subspace in $\Mat_d(\overline \FF_q)$. Suppose the contrary holds. Then there is some affine-linear equation in the entries of $\rho(a)$ which holds for a positive proportion of $a$. This is a polynomial equation $f(a_{11},\ldots,a_{nn})=0$ in the entries of $a$ (as can be seen by any explicit description of the representations of $\GL_n$). The equation $f=0$ either holds identically or cuts out some variety in $\overline\FF_q^{n\times n}$ of positive codimension. In the latter case, the Schwarz--Zippel lemma implies that the proportion of $a\in\Mat_n(\FF_q)$ satisfying $f=0$ is at most $C_n/q$ for some constant $C_n$ depending only on $n$. In the former case, the entirety of $G$ lives in some affine subspace of $\Mat_d(\overline\FF_q)$. Burnside's theorem, applied to the subalgebra generated by $\{b-1:b\in G\}$, then implies that $G$ has some invariant subspace. For large $q$, this contradicts the irreducibility of $\rho$.
\end{remark}

Lastly, we present a question in pure group theory which arises naturally from the statement of \cref{lem:sectional-p-rank}. Despite the simplicity of this statement, it has, to our knowledge, not yet appeared in the literature.

\begin{question}\label{qn:quot-to-subgp} For which pairs $(H_1,H_2)$ of finite groups does the following statement hold: every finite group $G$ with a quotient isomorphic to $H_1$ has a subgroup isomorphic to $H_2$?
\end{question}

Write $H_1\implies H_2$ if the condition in \cref{qn:quot-to-subgp} holds true for $(H_1,H_2)$. 

\cref{qn:quot-to-subgp} can be motivated by considering the category of abelian groups. The fundamental theorem of finite abelian groups implies that every quotient of such a group is isomorphic to a subgroup of the same group. In other words, if all groups (including the ambient group $G$) in \cref{qn:quot-to-subgp} were required to be abelian, then a pair $(H_1,H_2)$ would satisfy the stated property if and only if $H_2$ is isomorphic to a subgroup of $H_1$. \cref{qn:quot-to-subgp} asks how poorly this property of finite abelian groups fails to transfer to non-abelian groups. As \cref{qn:quot-to-subgp} asks about the subgroup structure of arbitrary finite extensions of $H_1$, it is not surprising that it seems a fairly difficult problem in general.

Indeed, we believe \cref{qn:quot-to-subgp} to be interesting even when $H_1$ and $H_2$ are abelian, and indeed essentially all we know about this question is in that restrictive setting. This limited knowledge is enumerated in the following list.

\begin{itemize}
    \item For any positive integer $m$ it holds that $\ZZ/m\ZZ\implies\ZZ/m\ZZ$. Indeed, if $G/N\cong\ZZ/m\ZZ$ then any lift of any generator of $G/N$ to $G$ has order a multiple of $m$, and thus has a power whose order is exactly $m$. 

    \item For $p$ prime, $(\ZZ/p\ZZ)^m\implies(\ZZ/p\ZZ)^{m'}$ if $m'<\sqrt{m/2}$. This is the statement of \cref{lem:sectional-p-rank}.
    
    \item On the other hand, for $p$ prime, $(\ZZ/p\ZZ)^m\nimplies(\ZZ/p\ZZ)^{m'}$ for $m'>\sqrt{8m+9}$. This is the main result of \cite{Olshanskii}. Following Ol'shanskii's methods, Halasi, Podoski, Pyber, and Szab\'o \cite[Theorem~1.9]{HPPS} show a stronger result where $p$ may be replaced by any prime power.
\end{itemize}

In particular, we know no cases in which $H_1\implies H_2$ for any non-abelian group $H_2$. It seems to us plausible, albeit unlikely, that no such examples exist. We present the following conjecture to the contrary.

\begin{conjecture}\label{conj:quot-to-subgp} For every positive integer $n$, there exists some positive integer $N$ for which $\operatorname S_N\implies\operatorname S_n$. 
\end{conjecture}

\noindent Indeed, \cref{conj:quot-to-subgp} is equivalent to, for every $H_2$, the existence of some $H_1$ satisfying $H_1\implies H_2$.

\section*{Acknowledgements}

We thank Jacob Fox for suggesting this area of study, for many helpful conversations throughout this work, and for numerous suggestions which greatly improved the final product. We also thank Daniel Altman, Amichai Lampert, and Jonathan Tidor for interesting discussions concerning \cref{sec:not-too-irreg}; Tristan Shin for pointing us to \cite[Theorem~21]{BloomSisask}; and Benjamin Church, Spencer Dembner, and Vaughan McDonald for providing algebro-geometric support (including, but not limited to, the proof of \cref{lem:GL-finiteness}).

The author is supported by the National Science Foundation Graduate Research Fellowship Program under Grant No.~DGE-2146755. Additionally, this material is partially based upon work supported by the National Science Foundation under Grant No.~DMS-1928930 while the author was in residence at the Simons Laufer Mathematical Sciences Institute in Berkeley, California during the semester of Spring 2025.

\bibliographystyle{alpha}
\bibliography{bib.bib}

\appendix

\section{Understanding subquotients: proofs of \texorpdfstring{\cref{lem:comp-rank-alternating,lem:sectional-p-rank}}{Lemmas 2.7 and 2.8}}\label{sec:appendix-cfsg}

In this section we prove the two purely group-theoretic facts we have needed in the body of the article: \cref{lem:comp-rank-alternating,lem:sectional-p-rank}.

We begin with a definition and a couple simple properties thereof. For a prime $p$ and a finite group $G$, we let $r_p(G)$ (resp.\ $s_p(G)$) denote the largest integer $t$ for which $G$ has a subgroup (resp.\ subquotient) isomorphic to $(\ZZ/p\ZZ)^t$. Following Guralnick and Tiep \cite{GuralnickTiepI}, we call $s_p(G)$ the \emph{sectional $p$-rank} of $G$. 

The sectional $p$-rank of a group has applications outside of this work. It is an important invariant in group cohomology, where it controls the dimension and complexity of group cohomology; see for example \cite{GuralnickTiepI,GuralnickTiepII,Symonds}. \cref{lem:sectional-p-rank} allows the results of the aforementioned papers to be stated in terms of existence of subgroups only; it is not clear to us whether $r_p(G)$ is also a natural invariant in that setting.

\subsection{Understanding the sectional \texorpdfstring{$p$}{p}-rank}

We begin with two simple properties of sectional $p$-rank.

\begin{lemma}\label{lem:srank-in-quotients} Let $G$ be a finite group and $N$ a normal subgroup. Then $s_p(G)\leq s_p(N)+s_p(G/N)$.    
\end{lemma}
\begin{proof} Let $H/K$ be a subquotient of $G$ isomorphic to $(\ZZ/p\ZZ)^s$, where $s:=s_p(G)$. Since $H\cap N\unlhd H$ and $H/(H\cap N)\leq G/N$, we have $s_p(N)+s_p(G/N)\geq s_p(H\cap N)+s_p(H/(H\cap N))$. It thus suffices to prove the result with $(G,N)$ replaced by $(H,H\cap N)$; in other words, we may assume $H=G$.

Now, the inclusion $N\unlhd G$ induces an inclusion $N/(K\cap N)\unlhd G/K$. Since $G/K$ is an elementary abelian group, both $N/(K\cap N)$ and
\[(G/K)/(N/(K\cap N))\cong (G/N)/(K/(K\cap N))\]
are elementary abelian groups, with ranks summing to $s$. Therefore
\[s_p(G)=\rk G/K=\rk\pr[\big]{N/(K\cap N)}+\rk\pr[\big]{(G/N)/(K/(K\cap N))}\leq s_p(N)+s_p(G/N).\qedhere\]
\end{proof}

\begin{lemma}\label{lem:reduce-to-p} For any finite group $G$, we have $r_p(G)=r_p(\Syl_p(G))$ and $s_p(G)=s_p(\Syl_p(G))$.
\end{lemma}
\begin{proof} The first identity follows from the fact that every $p$-subgroup of $G$ is contained within some Sylow $p$-subgroup of $G$. For the second, since $\Syl_p(G)\leq G$, it is obvious that $s_p(\Syl_p(G))\leq s_p(G)$. For the reverse inequality, let $s=s_p(G)$, and let $H/K$ be a subquotient of $G$ isomorphic to $(\ZZ/p\ZZ)^s$. Let $S_K$ be a Sylow-$p$ subgroup of $K$, and let $S_H$ be a Sylow-$p$ subgroup of $H$ containing $S_K$, so that $S_K=S_H\cap K$. Then we have
\[S_H/S_K=S_H/(S_H\cap K)\leq H/K.\]
Moreover, $\abs{S_H}/\abs{S_K}$ is the $p$-part of $\abs{H}/\abs{K}=\abs{H/K}$; since $H/K$ is a $p$-group, we conclude $S_H/S_K=H/K$. Therefore $H/K$ is a quotient of $\Syl_p(H)$, and thus $s_p(\Syl_p(G))\geq s_p(\Syl_p(H))\geq s=s_p(G)$.
\end{proof}

Given \cref{lem:reduce-to-p}, \cref{lem:sectional-p-rank} follows transparently from the following result in $p$-groups.

\begin{lemma}[Thompson--Mann; see {\cite[Theorems~A~and~B]{Mann}}]\label{lem:sectional-for-p-groups} For any finite $p$-group $G$ we have $s_p(G)\leq 2r_p(G)^2$.
\end{lemma}

In fact, something stronger is true: $r_p(G)$ can be replaced by the largest rank of an abelian \emph{normal} subgroup of $G$. In the case of $p$ odd, this is a classical result of Thompson (see \cite[Theorem~III.12.3]{Huppert}), who proved that if every abelian normal subgroup of a $p$-group $G$ can be generated by at most $r$ elements then every subgroup of $G$ can be generated by at most $r(r+1)/2\leq r^2$ elements. For $p=2$, this follows from a result of Mann \cite[Theorem~B]{Mann}, and seems also to have been discovered by Ol'shanskii \cite{Olshanskii}, based on an earlier result of Merzljakov \cite[Theorem~1]{Merzljakov}. We remark that Ol'shanskii also provides a construction of $p$-groups with $s_p(G)=\Omega(r_p(G)^2)$. So, \cref{lem:sectional-for-p-groups} (and thus \cref{lem:sectional-p-rank}) are tight up to a constant factor.

\subsection{Alternating subquotients}

We now prove \cref{lem:comp-rank-alternating}; we begin by proving part (1). We first need the following simple lemma.

\begin{lemma}\label{lem:simple-subquot} Let $G$ be a finite group and let $S$ be a finite simple group. If $G$ contains $S$ as a subquotient, then some composition factor of $G$ contains $S$ as a subquotient.
\end{lemma}
\begin{proof} Let $G=G_0\unrhd G_1\unrhd \cdots\unrhd G_s=\{1\}$ be a composition series of $G$. Let $H\leq G$ and $N\unlhd H$ be such that $S\cong H/N$. The groups $(H\cap G_i)/(N\cap G_i)$ form a subnormal series of $H/N\cong S$, and each successive quotient of this series is a subquotient of the corresponding composition factor of $G$. The result follows from the fact that $S$ is simple: any subnormal series of $S$ contains $S$ as a successive quotient.
\end{proof}

Now, \cref{lem:comp-rank-alternating}(1) follows from the following.

\begin{lemma}\label{lem:no-Ad-1} Let $\FF$ be any field and $G\leq\GL_r(\FF)$ be a finite group. If $G$ contains $\A_d$ as a subquotient, then $r\geq d/7-1$.
\end{lemma}

A statement similar to \cref{lem:no-Ad-1} has appeared (for $\GL_r(\ZZ/m\ZZ)$ with arbitrary $m$, as opposed to $\GL_r(\FF)$ for arbitrary $\FF$) in a question asked by Pete L.\ Clark and answered by Peter Mueller on MathOverflow \cite{mathoverflow}. We present a self-contained proof of \cref{lem:no-Ad-1}, as we are able to improve quantitatively on Mueller's argument.

Our first step is to argue that we can assume $\FF$ is finite.

\begin{lemma}\label{lem:GL-finiteness} Let $G$ be a finite group which embeds in $\GL_r(\FF)$ for some field $\FF$. Then there exists a finite field $\FF_q$ for which $G$ embeds in $\GL_r(\FF_q)$.
\end{lemma}

\cref{lem:GL-finiteness} follows from combining \cite[Theorems~VI--VIII]{Malcev}. For the reader's convenience, we sketch a proof using ring theory provided to us by Spencer Dembner and Benjamin Church.

\begin{proof} Let $x_{g,i,j}$ for $g\in G$ and $1\leq i,j\leq r$ be variables and define matrices $(m_g)_{g\in G}$ by $(m_g)_{ij}=x_{g,i,j}$. Let $A$ be the finite type $\ZZ$-algebra generated by $\{x_{g,i,j}\}$ with relations $x_{1,i,j}=1$ if $i=j$ and $0$ if $i\neq j$, as well as $m_gm_h=m_{gh}$ for each $g,h\in G$. 

We are given the existence of an embedding $G\hookrightarrow\GL_r(\FF)$ for some field $\FF$. This embedding gives us a ring homomorphism $\varphi\colon A\to\FF$. Let $I:=\ker\varphi$ and $B:=A/I$. The ring $B$ embeds into $\FF$ and is thus an integral domain. Moreover, since $\varphi(m_g)\neq\varphi(m_1)$ for each $g\in G\setminus\{1\}$, we can find, for each such $g$, some $(i_g,j_g)$ for which $x_{g,i_g,j_g}-x_{1,i_g,j_g}\not\in I$. Consider the localization $B'$ of $B$ at all $x_{g,i_g,j_g}-x_{1,i_g,j_g}$ for all $g\in G\setminus\{1\}$. The ring $B'$ possesses some maximal ideal $\mathfrak m$; since $B'$ is a finite-type $\ZZ$-algebra, the quotient $B'/\mathfrak m$ is a finite field  (see \cite[Corollary~5.3.4.1]{Bourbaki}). The quotient of $B$ by the preimage $\mathfrak p$ of $\mathfrak m$ embeds into $B'/\mathfrak m$, and is thus isomorphic to some finite field $\FF_q$. The quotient map $\pi\colon B\to B/\mathfrak p\cong\FF_q$ gives us the existence of a homomorphism $\varphi'\colon G\to\GL_r(\FF_q)$. Since $\mathfrak p=B\cap\mathfrak m$ for an ideal $\mathfrak m$ of $B'$, we have $\pi(x_{g,i_g,j_g})\neq \pi(x_{1,i_g,j_g})$, and thus $\varphi'(g)\neq\varphi'(1)$, for every $g\in G$. Therefore the homomorphism $\varphi'$ is injective, as desired.
\end{proof}

\begin{proof}[Proof of \cref{lem:no-Ad-1}] By \cref{lem:GL-finiteness}, we have $G\leq\GL_r(\FF_q)$ for some prime power $q=p^k$. Choose some prime $\ell\in\{3,5\}$ distinct from $p$; by extending our field if necessary, we can assume $\ell\mid q-1$. Our hypothesis that $G\leq\GL_r(\FF_q)$ has an $\A_d$-subquotient implies that $\GL_r(\FF_q)$ has an $\operatorname S_{d-2}$ subquotient. In particular, since $\operatorname S_{d-2}$ has a subgroup isomorphic to $(\ZZ/\ell\ZZ)^{\lfloor{(d-2)}/\ell\rfloor}$, we conclude
\begin{equation}\label{eq:s-rank-GL}
s_\ell\pr*{\GL_r(\FF_q)}\geq\floor*{\frac{d-2}\ell}.
\end{equation}

For an integer $n$, we write $\nu_\ell(n)$ for the largest $\nu$ for which $\ell^\nu\mid t$. Consider the subgroup $H\cong(\FF_q^\times)^r\rtimes\operatorname S_r$ of $\GL_r(\FF_q)$ consisting of products of permutation matrices and diagonal matrices. The index of $H$ in $\GL_r(\FF_q)$ is
\[\prod_{t=1}^r\frac{q^{r-t}(q^t-1)}{t(q-1)}.\]
Since $\ell\mid q-1$, the lifting-the-exponent lemma implies that $\nu_\ell(q^t-1)=\nu_\ell(q-1)+\nu_\ell(t)$. We conclude that the index of $H$ in $\GL_r(\FF_q)$ is relatively prime to $\ell$. In particular, $H$ contains a Sylow-$\ell$ subgroup of $\GL_r(\FF_q)$. We conclude from \cref{lem:reduce-to-p} and \eqref{eq:s-rank-GL} that
\begin{equation}\label{eq:s-rank-H}
s_\ell(H)=s_\ell\pr*{\GL_r(\FF_q)}\geq\floor*{\frac{d-2}\ell}.
\end{equation}
Now, since $H$ has a normal subgroup $(\FF_q^\times)^r\cong(\ZZ/(q-1)\ZZ)^r$ with quotient $\operatorname S_r$, \cref{lem:srank-in-quotients} implies that
\[s_\ell(H)\leq s_\ell\pr*{(\ZZ/(q-1)\ZZ)^r}+s_\ell(\operatorname S_r)\leq r+\nu_\ell(r!)\leq\frac{\ell r}{\ell-1}.\]
We conclude from \eqref{eq:s-rank-H} and our choice of $\ell\in\{3,5\}$ that
\[d\leq \ell\cdot s_\ell(H)+\ell+1\leq \frac{\ell^2}{\ell-1}r+\ell+1\leq 7r+6.\qedhere\]
\end{proof}

To prove \cref{lem:comp-rank-alternating}(2), we need to use the classification of finite simple groups and some properties of said groups. Our source for this is the book of Wilson \cite{Wilson}. 

\begin{theorem}[Classification of finite simple groups, sourced from {\cite[Section~1.2]{Wilson}}]\label{thm:CFSG} Let $S$ be a finite simple group. Then one of the following holds:
\begin{enumerate}[(a)]
    \item $S$ is a cyclic group of prime order.
    \item $S$ is an alternating group.
    \item $S$ is a classical group, either linear $\operatorname{PSL}_n(q)$, unitary $\operatorname{PSU}_n(q)$, symplectic $\operatorname{PSp}_{2n}(q)$, or orthogonal $\operatorname{P\Omega}_n^{\pm}(q)$.
    \item $S$ is an exceptional group of Lie type.
    \item $S$ is one of the $26$ sporadic simple groups.
\end{enumerate}
\end{theorem}

We mainly need properties of the classical and exceptional groups. We consider first the classical groups.

\begin{lemma}\label{lem:cfsg-classical} If $S$ is a classical group in $\{\operatorname{PSL}_n(q),\operatorname{PSU}_n(q),\operatorname{PSp}_n(q),\operatorname{P\Omega}_n^{\pm}(q)\}$, then
\begin{enumerate}
    \item $S$ can be embedded into $\GL_r(\FF_{q^2})$ for $r=n^2$, and
    \item $S$ contains a subgroup isomorphic to $\A_d$ for $d=\lfloor n/2\rfloor$.
\end{enumerate}
\end{lemma}

\begin{proof} To prove (1), we first note that $S$ is a subgroup of $\PGL_n(\FF_{q^2})$, the quotient of $\GL_n(\FF_{q^2})$ by the scalar matrices:
\begin{itemize}
    \item The group $\operatorname{PSL}_n(q)$ is the quotient of $\SL_n(q)\leq\GL_n(\FF_q)$ by scalar matrices. 

    \item The group $\operatorname{PSU}_n(q)$ is defined as the quotient by the scalar matrices of $\operatorname{SU}_n(q)$, which is defined as a subgroup of $\GL_n(\FF_{q^2})$ \cite[Section~3.6]{Wilson}.
    
    \item Let $n=2m$. The group $\operatorname{PSp}_{2m}(q)$ is defined as the quotient by the scalar matrices of $\operatorname{Sp}_{2m}(q)$, which is defined as a subgroup of $\GL_{2m}(\FF_q)$ \cite[Section~3.5]{Wilson}. It is thus a subgroup of $\PGL_n(\FF_{q^2})$.

    \item Regardless of the parity of $n$ and the sign of the $\pm$, the group $\operatorname{P\Omega}_n^{\pm}(q)$ is defined as the quotient by the scalar matrices of $\operatorname{GO}(\FF_q,f)$ for some non-singular symmetric bilinear form $f$, which itself is defined as a subgroup of $\GL_n(\FF_q)$ \cite[Sections~3.7~and~3.8]{Wilson}.
\end{itemize}
Property (1) now follows by embedding $\PGL_n(\FF_{q^2})$ into $\GL(\Mat_n(\FF_{q^2}))\cong\GL_{n^2}(\FF_{q^2})$ by the adjoint action $g\colon m\mapsto gmg^{-1}$.

We now prove (2).

\begin{itemize}
    \item The embedding $\operatorname{S}_n\hookrightarrow\GL_n(\FF_q)$ via the permutation matrices gives an embedding $\A_n\hookrightarrow\SL_n(q)$ which induces an embedding $\A_n\hookrightarrow\operatorname{PSL}_n(q)$.

    \item By choosing a basis, one can write $\operatorname{SU}_n(q)$ as the subgroup of $\SL_n(\FF_{q^2})$ consisting of matrices $g$ satisfying $\sigma(g)^\intercal g=I$, where $\sigma$ is the entrywise map $x\mapsto x^q$. With this description, each even permutation matrix lies in $\operatorname{SU}_n(q)$, and we thus get an embedding $\A_n\hookrightarrow\operatorname{PSU}_n(q)$.

    \item Let $n=2m$. By choosing a basis, one can write
    \[\operatorname{Sp}_{2m}(q)=\set*{\begin{pmatrix}A&B\\C&D\end{pmatrix}:A,B,C,D\in\Mat_m(\FF_q),\begin{matrix}A^\intercal C=C^\intercal A,\\B^\intercal D=D^\intercal B,\\A^\intercal D-C^\intercal B=I\end{matrix}}.\]
    We can thus embed $\SL_m(q)$ into $\operatorname{Sp}_{2m}(q)$ by mapping $g$ to $\operatorname{diag}(g,(g^{-1})^\intercal)$. This induces an embedding $\operatorname{PSL}_m(q)\hookrightarrow\operatorname{PSp}_{2m}(q)$, and thus an embedding of $\A_m=\A_{n/2}$ into $\operatorname{PSp}_n(q)$.

    \item The case of orthogonal groups is a bit messier. By \cite[Equation~3.34]{Wilson}, there are inclusions $\operatorname{GO}_{2m+1}(q)\leq\operatorname{GO}_{2m+2}^\pm(q)$ when $q$ is odd and $\operatorname{Sp}_{2m}(q)\cong\operatorname{GO}_{2m+1}(q)\leq\operatorname{GO}_{2m+2}^\pm(q)$ when $q$ is even \cite[Section~3.8.2]{Wilson}. Combined with the previous bullet point, this is enough to prove (2) when $q$ is even. When $q$ is odd, it is enough to show that $\A_{2m+1}\leq\operatorname{P\Omega}_{2m+1}(q)$. To do this, we may simply choose a basis in which $\operatorname{GO}_{2m+1}(q)$ consists exactly of those matrices $g$ satisfying $g^\intercal g=I$. The even permutation matrices thus all lie in $\operatorname{SO}_{2m+1}(q)$ under this basis, and thus project down to $\operatorname{P\Omega}_{2m+1}(q)$. \qedhere
\end{itemize}
\end{proof}

\begin{lemma}\label{lem:cfsg-exceptional} If $S$ is an exceptional group of Lie type, then $S$ can be embedded in $\GL_{10^4}(\FF_q)$ for some prime power $q$.
\end{lemma}
\begin{proof} We treat each of the classes of exceptional groups of Lie type in turn.
\begin{itemize}
    \item The group $G_2(q)$ is contained in $\operatorname{SO}_7(q)\leq\GL_7(\FF_q)$ \cite[Section~4.3.2]{Wilson}.

    \item For every prime power $q$, the group $F_4(q)$ can be defined as a subgroup of $\GL_{26}(\FF_q)$ \cite[Section~4.8.4]{Wilson}.

    \item The group $E_6(q)$ is defined as the quotient of a subgroup of $\GL_{27}(\FF_q)$ by its scalars \cite[Section~4.10.1]{Wilson}. Using the adjoint action as in the proof of \cref{lem:cfsg-classical}, $E_6(q)$ can be realized as a subgroup of $\GL_{729}(\FF_q)$.

    \item The group $^2E_6(q)$ is defined as the quotient of a subgroup of $E_6(q^2)$ by scalars, and is thus contained in $\GL_{729}(\FF_q)$ \cite[Section~4.11]{Wilson}. 
    
    \item The group $^3D_4(q)$ can be defined as a subgroup of $\operatorname{P\Omega}_8^+(q^3)$ \cite[Section~4.6.1]{Wilson}, and thus is a subgroup of $\GL_{64}(\FF_{q^6})$ by \cref{lem:cfsg-classical}(2).

    \item The group $E_7(q)$ can be defined as the quotient of a subgroup of $\GL_{56}(\FF_q)$ by its scalars \cite[Section~4.12.3]{Wilson}. Using the adjoint action as in the proof of \cref{lem:cfsg-classical}, $E_7(q)$ can be realized as a subgroup of $\GL_{56^2}(\FF_q)$.

    \item The group $E_8(q)$ is defined as the automorphism group of a $248$-dimensional Lie algebra \cite[Section~4.12.1]{Wilson}. As a result, it is a subgroup of $\GL_{248}(\FF_q)$.

    \item The group $^2B_2(2^{2n+1})=\operatorname{Sz}(2^{2n+1})$ is defined as a subgroup of $\operatorname{Sp}_4(2^{2n+1})$ \cite[Section~4.2.1]{Wilson}, which is itself a subgroup of $\GL_4(\FF_{2^{2n+1}})$.

    \item The group $^2G_2(3^{2n+1})$ is defined as a subgroup of $G_2(3^{2n+1})$ \cite[Section~4.5.1]{Wilson}, which we have observed to be a subgroup of $\GL_7(\FF_{3^{2n+1}})$.

    \item The group $^2F_4(2^{2n+1})$ is defined as a subgroup of $F_4(2^{2n+1})$ \cite[Section~4.9.1]{Wilson}, which we have observed to be a subgroup of $\GL_{26}(\FF_{2^{2n+1}})$.

    \item The Tits group $^2F_4(2)'$ is an index-$2$ subgroup of $^2F_4(2)$ \cite[Section~4.9.4]{Wilson}, and is thus a subgroup of $\GL_{26}(\FF_2)$. \qedhere
\end{itemize}
\end{proof}

We are now ready to prove \cref{lem:comp-rank-alternating}(2).

\begin{lemma}\label{lem:no-Ad-2} There exists an absolute constant $C$ such that, if $S$ is a finite simple group with no copy of $\A_d$ as a subquotient, then there exists some $r\leq Cd^2$ and some finite field $\FF$ for which $S\leq\GL_r(\FF)$.
\end{lemma}
\begin{proof} Let $C$ exceed the order of each of the $26$ sporadic simple groups. We use \cref{thm:CFSG}. 

If $S\cong\ZZ/p\ZZ$ is a cyclic group of prime order, then $S$ embeds into $\GL_1(\CC)\cong\CC^\times$ via the $p$th roots of unity.

If $S$ is an alternating group, then $S\cong\A_r$ for some $r<d$, and so $S$ embeds into $\GL_{d-1}(\FF_2)$ via the permutation matrices. 

If $S$ is a classical group, then \cref{lem:cfsg-classical}(2) implies that
\[S\in\{\operatorname{PSL}_n(q),\operatorname{PSU}_n(q),\operatorname{PSp}_n(q),\operatorname{P\Omega}_n^{\pm}(q)\}\]
for some $n\leq 2d+1$ and some prime power $q$. Now \cref{lem:cfsg-classical}(1) implies that $S$ can be embedded ingo $\GL_r(\FF_{q^2})$ for some $r\leq n^2\leq (2d+1)^2\leq Cd^2$.

If $S$ is an exceptional group of Lie type, then \cref{lem:cfsg-exceptional} implies that $S$ embeds into $\GL_{10^4}(\FF_q)\leq\GL_{Cd^2}(\FF_q)$ for some prime power $q$. 

If $S$ is a sporadic group, then $S$ embeds into $\GL_{\abs{S}}(\FF_2)\leq\GL_C(\FF_2)\leq\GL_{Cd^2}(\FF_2)$ via the right regular representation.
\end{proof}
\end{document}